\documentclass[10pt]{amsart}
\usepackage{chngpage, float, graphicx, tabularx}
\title{}
\author{}
\newcommand{\pt}{\on{pt}}
\newcommand{\IC}{\on{IndCoh}}
\newcommand{\ICn}{\on{IndCoh}_{\on{nilp}}}
\newcommand{\fgc}{\widehat{\fg}_{\kappa_c}}
\newcommand{\cB}{\check{B}}
\newcommand{\cN}{\check{N}}
\newcommand{\fgk}{\widehat{\fg}_\kappa}
\newcommand{\Op}{\mathrm{Op}_{\cG}}
\newcommand{\Fp}{\mathbb{F}_p}
\newcommand{\gkc}{\gk_{\kappa_c}}
\newcommand{\TwoCoh}{\on{2-QCoh}}
\newcommand{\TwoIC}{\on{2-IndCoh}_{\on{nilp}}}
\newcommand{\LS}{\on{LocSys}_{\cG}}
\newcommand{\Fqt}{\Fq(\!(t)\!)}

\newcommand{\udHom}{\on{Hom}_{\DGCat}}
\newcommand{\sB}{\on{Fl}_G}

\newcommand{\sM}{\mathscr{M}}
\newcommand{\fg}{\mathfrak{g}}
\newcommand{\fh}{\mathfrak{h}}
\newcommand{\on}{\operatorname}
\renewcommand{\mod}{\on{-mod}}
\newcommand{\Dmod}{\on{D-mod}}
\newcommand{\Fq}{\mathbb{F}_q}
\newcommand{\Fun}{\on{Fun}}
\newcommand{\Vect}{\on{Vect}}
\newcommand{\rmod}{\on{mod-}}

\theoremstyle{remark}
\newtheorem{defn}[subsubsection]{Definition}
\newtheorem{re}[subsubsection]{Remark}

\newtheorem{ex}[subsubsection]{Example}

\usepackage{amssymb,amsmath,amscd, amsthm,xspace,color,tocvsec2, array,
	tikz-cd, aligned-overset}
\usepackage[all,cmtip]{xy}
\usepackage[mathscr]{euscript}
\usepackage[breaklinks=true]{hyperref}
\usepackage[myheadings]{fullpage}

\newcommand{\fk}{\mathfrak{k}}

\newcommand{\sL}{\mathscr{L}}
\newcommand{\sA}{\mathscr{A}}
\newcommand{\sT}{\mathscr{T}}

\newcommand{\sP}{\mathscr{P}}

\newcommand{\sS}{\mathscr{S}}

\newcommand{\Sp}{\on{Spec }}

\newcommand{\sE}{\mathscr{E}}

\newcommand{\End}{\operatorname{End}}
\newcommand{\Hom}{\operatorname{Hom}}
\newcommand{\CC}{\mathbb{C}}

\newcommand{\sN}{\mathscr{N}}

\newcommand{\Gr}{\operatorname{Gr}}

\newcommand{\fsl}{\mathfrak{sl}}

\newcommand{\ft}{\mathfrak{t}}

\newcommand{\fz}{\mathfrak{z}}

\newcommand{\id}{\text{id}}

\newcommand{\sF}{\mathscr{F}}
\newcommand{\sG}{\mathscr{G}}

\newcommand{\tr}{\operatorname{tr}}

\newcommand{\Bun}{\operatorname{Bun}}
\newcommand{\Fl}{\operatorname{Fl}}

\newcommand{\C}{\ensuremath{\mathbb{C}}}

\newcommand{\Z}{\ensuremath{\mathbb{Z}}}

\newcommand{\comment}[1]{}

\newcommand{\OO}{\mathscr{O}}

\newcommand{\Spec}{\text{Spec }}
\newcommand{\Sym}{\operatorname{Sym}}

\newcommand{\Ext}{\operatorname{Ext}}

\newcommand{\fb}{\mathfrak{b}}
\newcommand{\gk}{\widehat{\fg}}
\newcommand{\fn}{\mathfrak{n}}
\newcommand{\fm}{\mathfrak{m}}

\newcommand{\sD}{\mathscr{D}}

\newtheorem{theo}[subsubsection]{Theorem}
\newtheorem*{theo*}{Theorem}

\numberwithin{equation}{section}
\newcommand{\fgl}{\mathfrak{gl}}

\newcommand{\cfg}{\check{\fg}}

\newcommand{\ck}{\underline{k}}

\newcommand{\sO}{\mathscr{O}}

\newcommand{\cG}{\check{G}}

\newcommand{\sC}{\mathscr{C}}

\newcommand{\msf}{\operatorname}

\newcommand{\sW}{\EuScript{W}_\kappa}

\newcommand{\Qlb}{\overline{\mathbb{Q}}_\ell}

\renewcommand{\mod}{\operatorname{-mod}}

\newcommand{\QCoh}{\on{QCoh}}

\newcommand{\DGCat}{\msf{DGCat}_{\on{cont}}}

\begin{document}
	
	
	\title[Categorical representations and local geometric Langlands]{An 
	informal introduction to categorical representation theory and 
	the local geometric Langlands program}

	\author{Gurbir Dhillon}
	\address{Yale University \\ New Haven, CT
		06511, USA}

	\email{gurbir.dhillon@yale.edu}

	\date{\today}

	\begin{abstract}  We provide a motivated introduction to the theory of 
	categorical actions of groups and the local geometric Langlands 
	program. Along the way we emphasize applications, old and new, to 
	the usual representation theory of reductive and affine Lie algebras.  
	\end{abstract}
	
	\maketitle

\setcounter{tocdepth}{1}
\tableofcontents

	\section{Introduction}
	\subsection{} The theory of categorical actions of groups and the local 
	geometric Langlands correspondence are two important actors in modern 
	representation theory which have undergone rapid development in the past 
	twenty years, with decisive contributions by Beilinson, Drinfeld, Frenkel, and Gaitsgory, and more recently Beraldo, Raskin, and Yang. 
	
	Both categorical actions and local geometric 
	Langlands are deeply intertwined with the representation theory of Lie 
	algebras in their history and major 
	applications, and we will 
	emphasize these connections throughout.

	On the one hand, the theory of categorical actions provides an efficient 
	framework for treating the representation theory of reductive Lie algebras. 
	This is a convenient entry point into the 
	theory of categorical actions 
	and is how we shall first meet them in the present survey. 
	
	On the other hand, and more significantly, the formalism of categorical 
	actions and perspectives from local geometric Langlands have proven 
	essential to the development of the representation theory of 
	affine Lie algebras and related vertex algebras, providing tools to 
	resolve old conjectures and generating new ones. We will reach these developments in the 
	latter half of this survey.

	\subsection{} 	
	
	A wonderful feature of this circle of ideas is that, despite the appearance 
	of fairly sophisticated tools from modern algebra and geometry, the basic 
	results, techniques, and perspectives can be quickly 
	absorbed by a nonspecialist if presented in the right light. It is our 
	hope that the present article can provide such an entry point for 
	the nonexpert reader.

	To this end, our orienting goal has been to 
	survey this area in 
	a manner which is accessible and complements what 
	is already available in the literature. In view of this, we 
	have chosen to emphasize analogies, motivations, and historical precursors 
	over providing a fully rigorous development of the treated material. 
	
   \subsection{} Two immediate 
	problems 
	 arise when 
	pursuing such a goal. 
	
     The first problem is how much to presume as known by the 
	reader and how much to explain in some detail. We will assume the reader 
	has some 
	familiarity with the basics of the representation theory of simple Lie groups 
	and Lie algebras, algebraic geometry, and homological algebra but not more 
	than what would be covered in an introductory course on these subjects. 
	Certainly these may be further substituted by some mathematical maturity 
	or at least the willingness to take a 
	few assertions on faith.  In particular, we will not assume any great 
	familiarity with the tools of geometric representation theory, 
	e.g., D-modules and perverse sheaves, higher algebra, and derived and 
	infinite 
	dimensional algebraic geometry. Instead, we will develop what we need as we 
	go along. In doing so, we have again attempted to provide treatments which 
	convey the basic ideas with minimal fuss, in a way which hopefully provides some 
	orientation to then tackle the relevant literature. 
	
		The second problem is that the choice of which ideas, methods, 
and 
	intuitions to emphasize is necessarily constrained by space as well as the 
	biases and the 
	limitations of the surveyor. Accordingly, we apologize to the many 
	mathematicians and physicists whose contributions were omitted 
	or underemphasized. As a first indication, the extremely rich connections 
	with three- 
	and four-dimensional supersymmetric gauge theory are only touched on 
	briefly towards the end.

	\subsection{} We have endeavored to organize the survey in a such a way 
	that the complexity of the material increases monotonically throughout. In 
	particular, the reader may wish to start at the beginning, skip to where 
	things seem unfamiliar or interesting, and read from there onwards. We hope 
	that the first sections, particularly Section \ref{s:bbloc}, can serve as a 
	light introduction to geometric representation theory for nonspecialists. 
	We also expect that to appreciate every last detail of the final sections, 
	particularly 
	Section \ref{s:affbbloc}, would require some degree of prior familiarity.

	The precise organization of the article is as follows. 
	
	In Section \ref{s:repthryfingroup} we recall some basic aspects of the 
	representation theory of finite groups to set the stage for their 
	categorical counterparts, which we will encounter later. 
	
	In Section \ref{s:bbloc} we review the localization theorem of 
	Beilinson--Bernstein and see how the resulting correspondence between 
	equivariant D-modules and Harish-Chandra modules for the Lie algebra leads 
	to the idea of a categorical group action.
	
	In Section \ref{s:catrepgrp} we then discuss the basics of the theory of 
	categorical group actions for finite dimensional algebraic groups and take 
	some time to explain some of the 
	 tools and constructions from modern homological algebra needed to get up 
	 and running.
	 
To get to local geometric Langlands, we must pass from the setting of algebraic 
groups to loop groups. Accordingly, in Section \ref{s:catloop}, we then set up 
the basics of the categorical representation theory of loop groups, i.e., review what changes when passing to the infinite dimensional setting.

	Finally, in Section \ref{s:affbbloc} we meet the local geometric Langlands 
	correspondence, discuss known and conjectured aspects in some detail, and 
	highlight an important application to localization theory for affine Lie 
	algebras at critical level. 
	
	We also provide in Appendix \ref{s:frmfun} a concise introduction to 
	D-modules, the Riemann--Hilbert correspondence, and the functions-sheaves 
	correspondence for the benefit of nonspecialists.

	 \subsection{} In addition to the references to the literature provided throughout, the reader is encouraged to consult the wonderful earlier survey \cite{frenkellocalnotes} and book \cite{frenkelloop} on local geometric Langlands and the notes \cite{paris-notes}; see also \cite{frenkelglobalnotes} for an accessible introduction to the global story.

	\vspace{.5cm}
	
	\noindent {\bf Acknowledgments.} We thank Pramod Achar, Dima Arinkin, Dario Beraldo, 
	Roman 
	Bezrukavnikov, Alexander Braverman, Justin Campbell, Harrison Chen, Thomas Creutzig,  Davide 
	Gaiotto, 
	Tony Feng,  Edward 
	Frenkel, Igor Frenkel,  
	Dennis Gaitsgory, Ivan Losev, Sam Raskin, Daniel Sage, Germ\'an Stefanich, Jonathan Wang, 
	David Yang, Zhiwei Yun, and Xinwen Zhu for conversations which have enhanced our 
	understanding of this subject. We are also deeply grateful to Cora Bernard, Apoorva Khare, and the anonymous referee for their careful reading and helpful feedback on earlier versions of the present text. 
	
	The author was supported by an NSF Postdoctoral Fellowship under grant No. 
	2103387.


	\section{Representations of group algebras and convolution algebras} 
	\label{s:repthryfingroup}
	\subsection{Overview} 
	
	\subsubsection{} Given a split reductive group $G$ over the 
	integers,\footnote{The reader will not lose much by taking $G$ throughout to 
	be the group $GL_n$ of invertible $n \times n$ matrices or even just 
	$GL_2$.} there are instructive analogies between the complex 
	representation theory of its points over finite fields, local fields, and 
	the real numbers. Using these, one can transfer constructions and ideas 
	from one area to the other. Harish-Chandra even dubbed this 
	permeability between the subjects as a Lefschetz principle for 
	representation theory \cite{langlands1985harish}.\footnote{Recall the usual Lefschetz principle 
	ensures that facts which are true about any complex algebraic variety, e.g., 
	generic smoothness, hold for algebraic varieties over any algebraically 
	closed field of characteristic zero.}

	Similarly, a basic principle in the categorical theory, with its origins in 
	the work of Drinfeld and Laumon on the global geometric Langlands 
	correspondence and Lusztig on character sheaves, is that the categorical 
	representation theory of $G(\mathbb{C})$ should behave similarly to (a part 
	of) the complex representation theory of $G(\mathbb{F}_q)$ and that the 
	categorical representation thory of $G(\mathbb{C}(\!(t)\!))$ should behave 
	similarly to the complex representation theory of 
	$G(\mathbb{F}_q(\!(t)\!))$.

	\subsubsection{} With this in mind, before moving to the categorical 
	setting, we first review some salient facts about the usual representations 
	of finite groups and then specialize the discussion to $G(\mathbb{F}_q)$.

	\subsection{The group algebra via convolution} 
	
	\subsubsection{} Suppose that $H$ is a finite group and $k$ is an 
	algebraically closed field of characteristic zero. Let us write $\on{Rep}(H)^\heartsuit$ for the abelian category 
	of representations of $H$ on $k$ vector spaces.
	
	\begin{re}The superscript on $\on{Rep}(H)^\heartsuit$ is to emphasize we are speaking of the abelian category, thought of as the heart of the standard $t$-structure on the corresponding derived category $\on{Rep}(H)$. Similarly, in what follows, categories of modules and sheaves are derived unless otherwise specified, and the corresponding abelian categories will be denoted with a superscript `$\heartsuit$'. 
	\end{re}

	\subsubsection{} Recall that to give a representation of $H$ on a vector 
	space $V$ is the data of a map
	\begin{equation} \label{e:hom}
	\alpha: H \rightarrow \on{End}(V).
	\end{equation}
	Here $\on{End}(V)$ denotes the $k$ linear endomorphisms of $V$, and $\alpha$ is 
	unital and intertwines multiplication in $H$ with composition in 
	$\on{End}(V)$. More generally, for any $k$-algebra $A$, consider the 
	collection of maps of semigroups 
	\begin{equation}  \label{e:homalg}
	H \rightarrow A,
	\end{equation}
	where $A$ is viewed as a semigroup under multiplication. A crucial 
	observation is that, since one is forgetting structure on $A$, namely addition, this extends 
	to a map from a $k$-algebra built from $H$, namely the group algebra. Let 
	us phrase the construction of the group algebra in a way convenient for the categorical 
	variant we shall meet later.

	\subsubsection{} \label{ss:opsfinitesets}For a finite set $X$, let us write $\on{Fun}(X)$ for the 
	vector space of $k$-valued functions on $X$. Let us recall two basic 
	properties of this assignment. First, a factorization of $X$ as a product 
	$X = X_1 \times X_2$ yields a canonical tensor product factorization of its 
	space of functions
	$$ \Fun(X) \simeq \Fun(X_1) \otimes \Fun(X_2).$$Second, given a map $f: X 
	\rightarrow Y$ of maps of finite sets, one has pushforward and pullback maps
	$$
	 f^!: \on{Fun}(Y) \rightarrow \Fun(X) \quad \text{and} \quad   
	f_*: \Fun(X) \rightarrow \Fun(Y). 
	$$
	Explicitly, $f^!$ sends a map $\phi: Y \rightarrow k$ to the 
	composition $\phi \circ f$, and  $f_*$ integrates a function $\psi$ in 
	$\Fun(X)$ along the fibres of $f$. I.e.,  for any element $y$ of $Y$ one has
	\[
	(f_*\psi)(y) := \underset{x \in f^{-1}(y)}{\Sigma} \hspace{.5mm} \psi(x). 
	\]
	%
	\subsubsection{} Combining these properties, the group structure on 
	$H$ endows $\Fun(H)$ with the structure of a $k$-algebra. Namely, if we 
	write 
	$
	m: H \times H \rightarrow H
	$
	for the multiplication map, this yields a convolution map on functions 
	\begin{equation} \label{e:convfun}
	\Fun(H) \otimes \Fun(H) \simeq \Fun(H \times H) \xrightarrow{ m_*}  
	\Fun(H). 
	\end{equation}

	Explicitly, if for an element $h$ of $H$ we denote by $\delta_h$ the 
	corresponding delta function in $\Fun(H)$, i.e., 
	\[
	\delta_h(h') := \begin{cases}   1 & \on{if } h = h', \\ 0 & \on{otherwise,} 
	\end{cases} 
	\]
	then the $\delta_h$, for $h \in H$, form a basis of $\Fun(H)$, and multiply 
	by the rule
	\[
	\delta_{h_1} \star \delta_{h_2} = \delta_{h_1h_2},
	\]
	where $h_1h_2$ denotes the product in $H$. The assignment $h \mapsto 
	\delta_h$ yields a map $H \rightarrow \Fun(H)$, with the property that any 
	map of the form \eqref{e:homalg} extends uniquely to a map of $k$-algebras 
	\[
	\Fun(H) \rightarrow A. 
	\]
	In particular, this yields the familiar equivalence between representations 
	of $H$ and left modules for $\Fun(H)$.
	
	Plainly, one forms $\Fun(H)$ by taking linear combinations of points of 
	$H$, and it acts in a representation of $H$ on a vector space by taking 
	the corresponding linear combination of endomorphisms.

	\subsection{The principal series representation} 
	
	\subsubsection{} In the remainder of this section, we will recall a basic 
	representation of $\mathbb{F}_q$-points of $G$, namely the space of 
	functions on its flag manifold. To do so, we first begin with some 
	preliminaries.

	\subsubsection{}  \label{ex:actspace}As before, $H$ denotes a finite group. 
	If $H$ acts on a finite set $X$, this induces a natural representation of 
	$H$ on $\Fun(X)$. Namely, given an element $h$ of $H$ and a function $\phi$ 
	in $\Fun(X)$, the action of $h$ on $\phi$ yields the element $h \cdot \phi$ 
	given by
	\[
	(h \cdot \phi)(x) := \phi(h^{-1} x), \quad \text{for } x \in X. 
	\]
	The corresponding left module structure 
	\begin{equation} \label{e:leftmod}
	\Fun(H) \otimes \Fun(X) \rightarrow \Fun(X)
	\end{equation}
	may be described as follows. If we denote the action by $a: H \times X 
	\rightarrow X$, then \eqref{e:leftmod} is given by the composition 
	\begin{equation} \label{e:actionfun}
	\Fun(H) \otimes \Fun(X) \simeq \Fun(H \times X) \xrightarrow{a_*} \Fun(X). 
	\end{equation}

	\subsubsection{}\label{s:funselfdual} We recall that $\Fun(X)$ is 
	canonically self dual. Namely, 
	writing $\Delta: X \rightarrow X \times X$ for the diagonal and $\pi: X 
	\rightarrow \on{pt}$ for the projection, one has the perfect pairing 
	\[
	\Fun(X) \otimes \Fun(X) \simeq \Fun(X \times X) \xrightarrow{\Delta^!} 
	\Fun(X) \xrightarrow{\pi_*} \Fun(\on{pt}) \simeq k,
	\]
	Explicitly, this is simply the `$L^2$' inner product 
	\[
	f \otimes g \mapsto \underset{x}{\Sigma} \hspace{1mm} f(x) \cdot g(x). 
	\]

	\subsubsection{} \label{sss:endos}
	
	In particular, given finite $H$-sets $X$ and $Y$, one has a canonical 
	identification 
	\[
	\on{Hom}_{\on{Rep}(H)^\heartsuit}( \Fun(X), \Fun(Y)) \simeq ( \Fun(X)^* \otimes 
	\Fun(Y))^{H} \simeq \Fun(X \times Y)^H \simeq \Fun(X \times Y/ H),
	\]
	where the superscript `$H$' denotes the subspace of $H$-invariants, and $X 
	\times Y / H$ denotes the orbit set of $X \times Y$ with respect to the 
	diagonal action of $H$. 
	
	Explicitly, given a function $K \in \Fun(X \times Y)$, the corresponding homomorphism $\Fun(X) \rightarrow \Fun(Y)$ is given by the integral transform 
	\[
	   \Fun(X) \xrightarrow{\pi_X^!} \Fun(X \times Y) \xrightarrow{ K \cdot - } \Fun(X \times Y) \xrightarrow{\pi_{Y, *}} \Fun(Y), 
	\]
	where $\pi_X: X \times Y \rightarrow X$ and $\pi_Y: X \times Y \rightarrow Y$ denote the projection maps, and $K \cdot -$ denotes pointwise multiplication of functions. Moreover, $K$ is $H$-invariant if and only if the corresponding integral transform is a map of $H$ representations.

	\subsubsection{} Let us specialize the preceding discussion to 
	$G(\mathbb{F}_q)$, where $\mathbb{F}_q$ denotes a finite field with $q$ 
	elements.

	Consider the flag manifold $\sB$, which parametrizes the Borel subgroups of $G$. This is a projective scheme over the 
	integers. The natural action of $G$ on $\sB$ by conjugation of Borel subgroups yields on 
	$\mathbb{F}_q$-points an action 
	\[
	G(\mathbb{F}_q) \times \sB(\mathbb{F}_q) \rightarrow \sB(\mathbb{F}_q). 
	\]
	In particular, as in Section \ref{ex:actspace} we obtain the unipotent principal series module 
	\[
	\Fun(G(\mathbb{F}_q)) \otimes \Fun(\sB(\mathbb{F}_q)) \rightarrow 
	\Fun(\sB(\mathbb{F}_q)).
	\]
	Moreover, its algebra of endomorphisms, the {\em Hecke algebra}, is canonically 
	identified as in Section \ref{sss:endos} with 
	\[
	\Fun( \sB(\mathbb{F}_q) \times \sB(\mathbb{F}_q) / G(\mathbb{F}_q)). 
	\]
	We emphasize that we are not considering the algebras of polynomial 
	functions on the corresponding $\Fq$-varieties, but simply the algebra of 
	$k$-valued functions on their finite sets of rational points. 
	
	\subsubsection{} Slightly more explicitly, recall that $\sB$ parametrizes 
	Borel subgroups of $G$. These subgroups are all conjugate and self-normalizing in 
	$G$. So, if one fixes a Borel $B$, one obtains an identification of the 
	flag manifold with the homogeneous space $G/B$, and moreover on $\Fq$-points
	\[
	\sB(\mathbb{F}_q) \simeq G(\Fq)/B(\Fq).
	\]
	Similarly, one obtains an identification of $(\sB \times \sB) /G$ with $B 
	\backslash G / B$ and hence an identification 
	\[
	\Fun( \sB(\mathbb{F}_q) \times \sB(\mathbb{F}_q) / G(\mathbb{F}_q)) \simeq 
	\Fun( B(\Fq) \backslash G(\Fq) / B(\Fq)). 
	\]
	This intertwines the algebra structure on the left with the algebra structure on the right given by convolution.

	\subsubsection{} For concreteness, let us discuss the smallest nontrivial 
	example in more detail. 
	
	\begin{ex} Suppose that $G = GL(V)$ for a free abelian group $V$ of rank 
	two. In this case, $\sB$ is the projective line $\mathbb{P}(V)$, which 
	parametrizes rank one quotients of $V$ with its natural action of $GL(V)$. 
	To see this, pick a basis $e_1, e_2$ of $V$, and note that the stabilizer 
	of the projection $V \rightarrow V / \langle e_2 \rangle$ is the Borel 
	subgroup of upper triangular matrices
		\[
		B =\begin{pmatrix} * & * \\ & * \end{pmatrix},
		\]
		which yields the desired identification with $G/B$. In particular, 
		$\sB(\mathbb{F}_q)$ has $q+1$ points.
		
		The module $\Fun(\sB(\mathbb{F}_q))$ decomposes into two simples as 
		follows. The projection $\pi: \sB \rightarrow \on{pt}$ yields a (split) 
		short exact sequence of $G(\mathbb{F}_q)$-modules
		\[
		0 \rightarrow \on{St} \rightarrow \Fun(\sB(\mathbb{F}_q)) \rightarrow 
		\Fun(\on{pt}(\mathbb{F}_q)) \simeq k \rightarrow 0,
		\]
		and the kernel $\on{St}$, the Steinberg module, is simple e.g., because the Hecke algebra is two dimensional. 
	\end{ex}

	\subsubsection{} The contents of this section may be summarized as follows. 
	First, we reviewed the construction of the group algebra of a finite group 
	as a convolution algebra. Second, for the finite group $G(\Fq)$, we 
	reviewed its natural action on functions on the flag manifold $\sB(\Fq)$.

	\section{Beilinson--Bernstein localization and hidden symmetries}
	\label{s:bbloc}
	
	\subsection{Overview} 
		\label{s:overview3}
	\subsubsection{} We considered in the previous section the 
	vector space of functions on the $\Fq$-points of the flag manifold
	\begin{equation} \label{e:funflgs}
	\Fun(\sB(\Fq)), 
	\end{equation}
	which naturally carried an action of $G(\Fq)$. 
	
	In this section, we will instead be interested in the derived category of D-modules on the 
	complex flag manifold
	\begin{equation} \label{e:dmodflg}
	\Dmod(\sB(\C)).
	\end{equation}
		We will begin by recalling the fundamental connection of this category \eqref{e:dmodflg} to the representation theory of Lie algebras following Beilinson--Bernstein. In reviewing several important properties of this correspondence, we will be led to the analogue of the action of $G(\Fq)$ on $\Fun(\Fl_G(\Fq))$. In this way we will meet the idea of categorical actions of groups.

	\subsection{The Localization Theorem}

	\subsubsection{} As the reader may have gleaned from our discussion in Section \ref{s:overview3}, one thinks of the category of D-modules \eqref{e:dmodflg} as a sort of categorification of \eqref{e:funflgs}. This analogy is absolutely fundamental in geometric representation theory. For this reason, we have provided a fairly detailed explanation of it, as well as an introduction to D-modules, in Appendix \ref{s:frmfun}. The nonspecialist reader may wish to turn there now and return to this section after picking up the basics.

\subsubsection{} In the remainder of this section, we work over an algebraically closed field $k$ of characteristic zero. For ease of notation, we redefine $\Fl_G$ to be the corresponding base change. That is, if we write $\Fl_{G, \mathbb{Z}}$ for the version over $\Spec \Z$ discussed in Section \ref{s:repthryfingroup}, we set 
\[
 \Fl_G := \Fl_{G, \Z} \underset {\Spec \Z} \times \Spec k.
\]
Of course, the reader may take $k$ to be the complex numbers but should bear in mind that we do not need anything particular about their metric topology.

	\subsubsection{} It turns out that the category of D-modules on the flag variety 
	$$\Dmod(\Fl_G)$$has a seemingly different realization, namely as a category of 
	Lie algebra representations. Let us review this remarkable correspondence, 
	which birthed geometric representation theory.

	\subsubsection{} Recall that a D-module $\sM$ on a smooth algebraic variety 
	$X$ is a quasicoherent sheaf equipped with a flat connection. Equivalently, 
	$\sM$ is equipped with an action of the sheaf of differential operators 
	$\sD_X$ extending the action of regular functions $\OO_X$. In particular, 
	its global sections will carry an action by global differential operators, 
	i.e., we have a tautological functor
	\begin{equation} \label{e:glob}
	\Gamma(X, -): \Dmod(X) \rightarrow \Gamma(X, \sD_X)\on{-mod.} 
	\end{equation}
	It is a remarkable theorem of Beilinson--Bernstein that, for $X = \sB$, this 
	functor is an equivalence.

	\subsubsection{} Let us explain why this is surprising. 
	
	First, if $X$ is affine, the analogous assertion is true but unsurprising. 
	Indeed, in this case it follows from the fact that the definition of 
	quasicoherent sheaves is rigged so that on an affine variety the functor
	\begin{equation} \label{e:globqcoh}
	\Gamma(X, -): \on{QCoh}(X) \rightarrow \Gamma(X, \OO_X)\mod
	\end{equation}
	is an equivalence. 
	
	By constrast, if $X$ is not affine, e.g., projective, the functor 
	\eqref{e:globqcoh} is hardly ever an equivalence. For example, for $\sB$ 
	its only global (derived) functions are scalars, and so \eqref{e:globqcoh} 
	reduces to the map
	\[
	\Gamma(\sB, -): \on{QCoh}(\sB) \rightarrow \on{Vect,} 
	\]
	where $\on{Vect}$ denotes the derived category of vector spaces. This is far from an equivalence, since $\sB$ is not a point! For this 
	reason, it is non-trivial that, upon further equipping quasicoherent 
	sheaves with flat connections, \eqref{e:glob} is an equivalence for the 
	flag manifold. 
	
	We also emphasize that, for a general projective variety, \eqref{e:glob} is 
	not an equivalence. For example, the reader familiar with algebraic curves 
	may wish to check that the only smooth projective curve for which 
	\eqref{e:glob} is an equivalence is $\mathbb{P}^1$, i.e., the smallest flag 
	manifold.

	\subsubsection{} Let us turn next to the description of global differential 
	operators 
	\[
	\Gamma(\sB, \sD_{\sB}). 
	\]
	Recall that for any smooth variety $X$, $\sD_{X}$ is generated as a sheaf 
	of algebras by the functions $\OO_X$ and vector fields $\sT_X$. Let us see what 
	they contribute to its global sections in the present case of $X = \sB$.

	As we mentioned above, the global functions on the flag manifold are simply 
	scalars 
	\[
	k \xrightarrow{\sim} \Gamma(\sB, \OO_{\sB}). 
	\]
	For the vector fields, the natural action of $G$ on $\sB$ yields 
	the infinitesimal action of its Lie algebra $\fg$, i.e., a homomorphism of 
	Lie algebras
	\begin{equation}\label{e:vecfields}
	\fg \rightarrow \Gamma(\sB, \sT_{\sB}). 
	\end{equation}
	To approximate the final answer, we might optimistically hope that this is 
	an isomorphism. 
	
	To account for higher order differential operators, we could similarly ask 
	that the resulting map from the enveloping algebra
	\begin{equation}\label{e:globdiff}
	U(\fg) \rightarrow \Gamma(\sB, \sD_{\sB})
	\end{equation}
	be an isomorphism. This would be a strong global analogue of the local 
	generation of sections of $\sD_{\sB}$ by functions and vector fields.

	\subsubsection{} As it turns out, both \eqref{e:vecfields} and 
	\eqref{e:globdiff} are not quite isomorphisms. Instead, they must be 
	corrected by considering the center.

	Let us begin with the case of vector fields. Note that the center $Z$ of 
	$G$ acts trivially on $\sB$. So, we need to quotient by the center $\fz$ of 
	$\fg$. After doing so, we do obtain an equivalence
	\begin{align*}
	\fg/\fz &\xrightarrow{\sim} \Gamma(\sB, \sT \hspace{.2mm}).
	\intertext{\hspace{1em}Similarly, when we pass to all differential 
	operators, we must impose a central quotient. Namely, if we write $Z(\fg)$ 
	for the center of $U(\fg)$, acting on the trivial representation of $\fg$ 
	yields a character $\chi: Z(\fg) \rightarrow k.$ Let us denote by 
	$U_0(\fg)$ the corresponding central quotient of the entire enveloping 
	algebra}
	U_0(\fg) &:= U(\fg) \underset{Z(\fg)} \otimes  k_\chi.
	\intertext{With this, as in the case of vector fields, \eqref{e:globdiff} 
	factors through the resulting quotient, and we obtain an equivalence}
	U_0(\fg) &\xrightarrow{\sim} \Gamma(\sB, \sD). 
	\end{align*}
	
	\label{s:casa}
	
	\begin{ex} Let us explicitly describe what $U_0(\fg)$ looks like in the 
	case of $\fg = \fgl_2$. We may decompose the Lie algebra as its center and 
	the traceless matrices
		\[
		\fgl_2 = \fz \oplus \fsl_2. 
		\]
		This induces decompositions of the enveloping algebra, and hence its center, as 
		\[
		U(\fgl_2) \simeq \on{Sym}(\fz) \otimes U(\fsl_2) \quad \text{and} \quad 
		Z(\fgl_2) \simeq \on{Sym}(\fz) \otimes Z(\fsl_2). 
		\]
		\label{s:cas}
		The center $Z(\fsl_2)$ is a polynomial ring in one variable, generated 
		by the Casimir element $\Omega$ of filtered degree two. Explicitly, if 
		we pick a basis element $z$ for $\fz$ and standard generators $f, h, 
		e$ for $\fsl_2$, we have that 
		\[
		Z(\fgl_2) \simeq k[ z, ef + fe + \frac{1}{2}h^2],
		\]
		where the second generator is $\Omega$. It is straightforward to see 
		both generators act by zero on the trivial representation. In 
		particular, $U_0(\fg)$ is the quotient of $U(\fgl_2)$ by the ideal generated by those two elements, i.e., 
		\[
		U_0(\fgl_2) = U(\fgl_2)/( \fz, \Omega) \simeq U(\fsl_2)/(\Omega).  
		\]

	\end{ex}

	\subsubsection{} 
	Let us denote the derived category of $U_0(\fg)$-modules, i.e., representations of 
	$\fg$ with trivial central character, by $\fg\on{-mod}_0$. Summarizing the 
	previous discussion, the theorem of Beilinson--Bernstein reads as follows. 
	
	\begin{theo} \label{t:bbloc}(Beilinson--Bernstein \cite{beilinson-bernstein81}) The functor of global sections yields 
	a $t$-exact equivalence 
		\[
		\Gamma(\sB, -): \Dmod(\sB) \rightarrow \fg\on{-mod}_0.       
		\]
	\end{theo}

	It is hard to overstate the significance of this theorem for the 
	development of representation theory. We confine ourselves here to a few 
	remarks. 
	
	\subsubsection{} First, by a version of Schur's lemma due to Quillen, every 
	simple $\fg$-module has a central character, i.e., the center $Z(\fg)$ acts 
	by scalars. Therefore, for the study of irreducible modules, working one 
	central character at a time is sufficient.
	
	The above theorem identifies those with the trivial central character and 
	simple D-modules on $\sB$. This allows one to translate problems in the study of such irreducible modules into geometric or topological problems. The latter problems may be approached 
	using the wealth of information known about the topology of algebraic 
	varieties, notably Hodge theory and its singular variants. 
	
	The original motivation for the theorem was such an application. Namely, 
	within the category of all $\fg$-modules, to any Borel subgroup $B$ of $G$ 
	one associates the category of $(\fg, B)$-modules, i.e., $\fg$-modules for 
	which the action of the Lie algebra of $B$ is integrable. 
	
	A problem which received significant study, beginning with Verma's thesis in 1966 \cite{verma66}, was the determination of the characters of the simple $(\fg, 
	B)$-modules. After their determination in several low rank cases by Jantzen \cite{jantz79} 
	and others, Kazhdan and Lusztig in 1979 formulated a general conjecture for 
	the simple $(\fg,B)$-modules with trivial central character \cite{kl79}. At the time of 
	their conjecture, they knew their formula was intimately related to the 
	Schubert subvarieties of $\sB$, or more precisely with the $!$-stalks of 
	their intersection cohomology D-modules. However, the connection with $\fg$-modules was not clear. 
	
	The work of Beilinson--Bernstein \cite{bb82}, and independently Brylinski--Kashiwara \cite{brylin81}, 
	provided the desired link. Namely, the localization theorem exchanged the 
	simple objects of $(\fg,B)\mod_0$ with the aforementioned intersection 
	cohomology sheaves, and the character formula followed.

	\begin{re} We have discussed above only the case of the 
	trivial central character. Beilinson--Bernstein in fact proved a similar 
	statement for any regular central character. For variants, including 
	generalized or singular central characters, see \cite{kashloc}, 
	\cite{beilinsonginzburg}, \cite{bmr}, \cite{krem}, \cite{ahc}. \end{re}

	\subsubsection{} We would like to next describe the behavior of the localization equivalence 
	in three basic examples. 
	
	\begin{ex} As we described in Section \ref{s:casa}, localization 
	interchanges the algebra of differential operators and the central quotient 
	of the enveloping algebra, i.e., 
		\[
		\Gamma(\sB, \sD) \simeq U_0(\fg). 
		\]
	\end{ex}
	
	\begin{ex}As we also discussed, localization interchanges the structure 
	sheaf and the trivial representation, i.e., 
		\[
		\Gamma(\sB, \OO) \simeq k. 
		\]
		As the reader may be aware, this is the simplest case of the 
		Borel--Weil--Bott theorem, which constructs all the irreducible 
		$G$-modules as global sections of line bundles on $\sB$. Informally, 
		the work of Beilinson--Bernstein  extends this to give a similar 
		construction of all irreducible $\fg$-modules. 
	\end{ex}

	\begin{ex} Finally, let us consider, for a point $x$ of $\sB$, the image of 
	$\delta_x$. Recall that the stabilizer of $x$ is a Borel subgroup $B$. 	
	With this, writing $\fb$ for the Lie algebra of $B$, we have that the delta 
	D-module is sent to the corresponding induced module
		\begin{equation}
		\Gamma(\Fl_G, \delta_x) \simeq U(\fg) \underset{U(\fb)} \otimes 
		\det(\fg/\fb).
		\end{equation}

		Let us sketch a proof of this, which the reader may enjoy thinking 
		through fully. To produce a map 
		\begin{equation} \label{e:hwv}
		U(\fg) \underset{U(\fb)} \otimes k_{\chi} \rightarrow \Gamma(\Fl_G, 
		\delta_x),
		\end{equation}
		we note that the identification of the generating line of the right hand side 
		with $\det(T_x \sB) = \det( \fg/\fb)$ is compatible with the 
		$\fb$-actions, cf. Remark \ref{r:hwdelta}. To see the obtained map is 
		an isomorphism, one notes that both sides are freely generated from 
		their determinant lines by the action of any subalgebra $\fn^-$ 
		transverse to $\fb$, i.e., the unipotent radical of an opposite Borel.

	\end{ex}

	\subsection{Equivariance and Harish-Chandra modules}
	
	\subsubsection{}Experts immediately understood that the localization theorem 
	carried far more 
	structure than simply an equivalence of categories. 
	
	For concrete 
	representation-theoretic applications, the main compatibility was as 
	follows. 
	Fix an algebraic subgroup $H \subset G$, which for simplicity we assume to 
	be 
	connected.
	
	\subsubsection{} Let us begin with the Lie theoretic side. Associated to 
	$H$ is 
	the full abelian subcategory of Harish-Chandra modules
	\[
	(\fg,H)\mod^\heartsuit \subset \fg\mod^\heartsuit. 
	\]
	Explicitly, a $\fg$-module is a Harish-Chandra module if the action of the 
	Lie algebra $\fh \subset \fg$ integrates to an action of $H$.  By our 
	assumption that $H$ is connected, this is a property and 
	not a further structure.

	Restricting to modules with trivial central character, we obtain a similar 
	full 
	category 
	\[
	(\fg, H)\mod_0^\heartsuit \subset \fg\mod_0^\heartsuit. 
	\]

	\subsubsection{}Here are some quick examples. 
	
	\begin{ex} 
		If $\fg$ is semisimple, the central quotient $U_0(\fg)$ is only a Harish-Chandra module for 
		the 
		trivial subgroup . 
	\end{ex}
	
	\begin{ex} The trivial module $k$ is a Harish-Chandra module for $G$ 
	itself and in particular for any subgroup. 
	\end{ex}

	\begin{ex} The Verma module $U(\fg) \underset{U(\fb)} \otimes 
	\det(\fg/\fb)$ is 
		a $(\fg,B)$-module. More generally, for any representation $W$ of $H$, 
		the 
		induced module 
		\[
		U(\fg) \underset{U(\fh)} \otimes W
		\]
		is a $(\fg,H)$-module, and any $(\fg, H)$-module is a quotient of such 
		a 
		module.  
	\end{ex}

	\subsubsection{} \label{s:eqdmod}Let us turn to the geometric side. 
	Associated to $H$ is 
	the 
	full abelian subcategory of $H$-equivariant D-modules 
	\[
	\Dmod(\sB)^{H, \heartsuit} \subset \Dmod(\sB)^\heartsuit.
	\]

	Let us review how this is defined in the generality of any algebraic 
	variety 
	$X$ equipped with an action of $H$. 
	
	\begin{re} For us, varieties are by definition quasicompact, quasiseparated, reduced schemes of finite type over $k$. Moreover, everything we state in this survey for varieties applies equally well to quasicompact, quasiseparated schemes of finite type over $k$, i.e., the reducedness plays no role.	The reader will lose little by ignoring these technicalities and sticking with e.g. quasi-projective varieties if they are more comfortable doing so. 
	\end{re}

	At first pass, for a D-module $M$ to 
	be 
	$H$-equivariant roughly means that it is constant along $H$-orbits within 
	$X$. To arrive at the actual definition, note that we should therefore 
	have, 
	for any 
	point $h$ of $H$ and $x$ of $X$, 
	the existence of an isomorphism of stalks
	\begin{equation}\label{e:stalks}
	i_x^! M \simeq i_{hx}^! M.
	\end{equation}
	Of course, since our spaces are not discrete, we need to put such 
	statements 
	into a family as we vary $h$ and $x$, compatibly with parallel transport. 
	To do 
	so, consider the 
	action and projection maps 
	\begin{equation}
	\alpha, \pi: H \times X \rightarrow X,\quad  \quad \alpha(h, x) = hx, 
	\quad 
	\pi(h, x) = x. 
	\end{equation}

	With this, an object $M$ of $\Dmod(X)^\heartsuit$ is $H$-equivariant if 
	there 
	exists an isomorphism \begin{equation} \alpha^! M \simeq \pi^! M. 
	\label{e:wholething}\end{equation}Since we only ask for the existence of an 
	isomorphism, it is clear this is a 
	property of a D-module and not further structure.

	It is instructive to compare this to the case of functions. Namely, fix a 
	finite group acting on a finite set, which we suggestively denote by 
	$H(\Fq)$ 
	and $X(\Fq)$, respectively. Note that a function $m$ in $\Fun(X(\Fq))$ is 
	$H(\Fq)$-invariant, i.e., constant along orbits, if and only if it satisfies 
	the 
	equality
	\[
	\alpha^! m = \pi^! m,
	\]
	where $\alpha, \pi: H(\Fq) \times X(\Fq) \rightarrow X(\Fq)$ are defined as 
	above, and the $!$-pullback of functions is as in Section \ref{ss:opsfinitesets}.

	\subsubsection{} Here is a pair of basic but important examples of equivariant D-modules. 
	
	\begin{ex} Suppose $X$ consists of a single $H$-orbit, i.e., after choosing 
	a base point, is of the form 
	 $H/K$ for some subgroup $K$ of $H$. 
	 
	 In this case, we claim there is a 
	 canonical equivalence between $H$-equivariant D-modules on $X$ and 
	 representations of the component group of $K$, i.e., 
	 \[
	        \Dmod(X)^{H, \heartsuit} \simeq \on{Rep}(\pi_0(K))^{\heartsuit}.
	 \]
	Indeed, it is straightforward to see that any equivariant D-module on $X$ 
	must be a local system, and the condition \eqref{e:wholething} is 
	equivalent to the pullback along 
	\[
	    H \rightarrow H/K
	\]
	being trivializable. On the other hand, the long exact sequence on homotopy 
	groups associated to the fibration $K \rightarrow H \rightarrow H/K$ gives 
	in low degrees an exact sequence 
	\[
	    \pi_1(H) \rightarrow \pi_1(H/K) \rightarrow \pi_0(K) \rightarrow 1,
	\]
	so the claim follows. \label{e:oneorbit}
	\end{ex}

	\begin{ex} Suppose $X$ consists of finitely many $H$-orbits. 
	In this case, the category of $H$-equivariant D-modules on $X$ is `glued' 
	from the individual orbits using recollement.
	
	The details of this are not so important for us. However, let us mention 
	the following basic consequence: the irreducible $H$-equivariant 
	D-modules on $X$ are in canonical bijection with the irreducible 
	$H$-equivariant local systems on the individual orbits. Concretely, each 
	irreducible object of $$\Dmod(X)^{H, \heartsuit}$$is supported on the 
	closure of a single orbit, and its restriction to the orbit is an 
	irreducible $H$-equivariant local system. In this way, the classification 
	of simples reduces, by the previous example, to an exercise in finite group 
	theory.

	\label{e:bigbrain}

	\end{ex}

	\subsubsection{} Having introduced the two sides, we may now state the 
	promised
	compatibility: localization exchanges Harish-Chandra modules and 
	equivariant 
	D-modules.
	
	\begin{theo}[Beilinson--Bernstein, \cite{beilinson-bernstein81}]\label{t:heqv} There is a commutative diagram 
	\[
	\xymatrix{\Dmod(\sB)^{H, \heartsuit} \ar[r]^\sim \ar[d] & (\fg, 
		H)\mod_0^\heartsuit \ar[d]  \\  \Dmod(\sB)^\heartsuit \ar[r]^\sim & 
		\fg\mod_0^\heartsuit,} 
	\]
	where the vertical arrows are the tautological inclusions and the bottom 
	horizontal arrow is localization. 
		\end{theo}
	
	\subsubsection{} In fact, Theorem \ref{t:heqv} holds, {mutatis 
	mutandis}, for any 
	$H$, 
	not necessarily connected and also for the appropriately defined derived 
	categories of equivariant D-modules and Harish-Chandra modules. The latter 
	assertion does not follow from the former, since these typically are not 
	the 
	derived categories 
	of their hearts, but rather something slightly more 
	sophisticated.

	\begin{re} The 
		reader is encouraged to view the aforementioned subtlety about derived categories as a statement of the form `to get the 
		correct 
		answers, one has to be a bit clever.' However, it does not take much to 
		see why 
		working naively with the derived category of the abelian category of equivariant sheaves can give the wrong answers. The orienting principle here 
		is that 
		$H$-equivariant D-modules on $X$ should be the same as D-modules on the 
		quotient stack $X/H$, cf. Remark \ref{r:introstacks} below, and should therefore carry the expected relation 
		with its 
		topology. In the simplest case of $X = \on{pt}$, we should 
		in particular want 
		\[
		\on{Hom}_{\Dmod(\on{pt})^H}(k, k) \simeq \on{C}^*(\on{pt}/H, k),
		\]
		i.e., the derived endomorphisms of the constant sheaf should be the 
		cohomology 
		of 
		the classifying space. But since on abelian categories
		\[
		\Dmod(\on{pt})^{H, \heartsuit} \simeq \on{Rep}(\pi_0(H)),
		\]
		we simply cannot get {\em any} higher Exts from the abelian category 
		alone. We 
		will 
		meet the correction, due to Bernstein--Lunts, in Section \ref{s:catrepgrp}.

	\end{re}
	
	\begin{re} \label{r:introstacks}For the reader unfamiliar with stacks, the following coarse description suffices for our purposes. An Artin stack $Y$ over an algebraically closed field $k$ is, roughly speaking, something like a variety, but with the further property that each $k$-point $y$ of $Y$  carries a group of `internal symmetries.'
	
		The most relevant class of Artin stacks for us arise as follows.  Given a variety $X$ acted on by an algebraic group $H$, one can form the quotient stack $X/H$. The points of $X/H$ are given by $H$-orbits on $X$, and for a $k$-point $x$ of  $X$, the internal symmetry group of the corresponding point $$H \cdot x \mathrm{\phantom{i}  of \phantom{i} } X/H$$is simply the stabilizer $H_x$ of $x$ in $H$. That is, this $k$-point of the stack is not the variety $\pt$, i.e.,   $\Spec k$, but rather its quotient $\pt / H_x$. In general, a $k$-point $y$ of an Artin stack is the classifying space $\pt / H_y$ of an algebraic group $H_y$, which is the internal symmetry group alluded to above.

		As a simple case, consider $X = \mathbb{A}^1$ with the  action of $H = \mathbb{Z}/2$ induced by the involution  $$z \mapsto -z.$$ If we pluck out the origin, this restricts to an action of $H$ on $\mathbb{G}_m$. Here, every $H$ orbit is free and, in particular, points of the quotient should have trivial internal symmetries, like an ordinary variety. In fact, the quotient $\mathbb{G}_m / H$ is simply the variety $\mathbb{G}_m$, where the quotient map 
		$
		   \mathbb{G}_m \rightarrow \mathbb{G}_m / H
		$
		identifies with the squaring map  $$\mathbb{G}_m \rightarrow \mathbb{G}_m, \quad w \mapsto w^2.$$ 
		However, the full quotient $\mathbb{A}^1 / H$ is not a variety, as the point corresponding to the origin has been `orbifolded' into a copy of $\pt / H$. Intuitively, for a general quotient stack, the internal symmetry groups remember the folding that takes place when forming the orbit space $X/H$ from $X$.

		 Finally, let us mention that for sheaves on stacks, the basic new feature compared to the case of varieties is that their fibres carry actions by the internal symmetry groups. Namely, for an algebraic group $H$, we have 
		 \[
		   \QCoh( \pt / H) \simeq \on{Rep}(H) \quad \text{and} \quad \Dmod( \pt / H) \simeq \on{C}_*(H)\mod, 
		 \]
		 where $\on{C}_*(H)$ denotes the dual of the de Rham cohomology of $H$, and is viewed as an dg-algebra under convolution. 
		 
		 Therefore, for a coherent sheaf on an Artin stack $Y$, its $*$-fibre at a $k$-point $y$ carries an action of $H_y$. Similarly, for a D-module on $Y$, its $!$-fibre at $y$ carries an action of $\on{C}_*(H_y)$, which one thinks of as a locally constant action of $H_y$. In particular, the cohomology groups of the $!$-fibre carry actions of the component group $\pi_0(H_y)$. 
	\end{re}

	\subsubsection{} Let us give a basic example of how Theorem \ref{t:heqv} is 
	used. 
	
It follows from \ref{t:heqv} that there exist only finitely 
	many irreducible objects in $(\fg, H)\mod_0$ if and only if there are 
	finitely many $H$-orbits on $\sB$. Moreover,  in this case, one can 
	explicitly enumerate the simple objects as in Example \ref{e:bigbrain}.

	As a concrete example,  write $T 
	\subset GL_2$ for a maximal torus, and let us classify the simple objects of
	\[
	(\fgl_2, T)\mod_0^\heartsuit. 
	\] 
	To do so, we note that $T$ has three orbits on $\mathbb{P}^1$ -- the north 
	pole $B_0$, the south pole $B_\infty$, and everything else. As $T$ acts 
	with connected 	stabilizers on all three orbits, it follows that one has 
	three irreducible objects. Concretely, these are the the simple Verma 
	module of highest weight $-2$, the simple Verma module of lowest weight 2,  
	and the trivial module, i.e., 
	\[
	      U(\fgl_2) \underset{U(\fb_0)} \otimes \det(\fgl_2/\fb_0), \quad 
	      U(\fgl_2) \underset{U(\fb_\infty)} \otimes \det(\fgl_2/\fb_\infty), 
	      \quad \text{and} \quad k. 
	\]
	Note also that, by similar reasoning, there are infinitely many simple irreducible 
	objects in $(\fgl_n, 
	T)\mod_0^\heartsuit$ for $n > 2$.

	While this discussion was rather soft, we emphasize that subtler questions, 
	such as the determination of simple characters, may be approached 
	geometrically and typically translate to questions regarding the topology of orbit 
	closures.  
	
	\subsubsection{} Finally, let us make a couple orienting remarks about why 
	Theorem \ref{t:heqv} holds. It is clarifying to work in the following 
	generality. \label{s:kandx}
	
	Suppose $K$ is an algebraic group acting on a smooth variety $X$. As before, if we 
	write $\fk$ for the Lie algebra of $K$, one 
	has a global sections functor 
	\[
	   \Gamma(X, -): \Dmod(X)^\heartsuit \rightarrow \fk\mod^\heartsuit.
	\]
	This admits a left adjoint, given by localization
	\[
	   \on{Loc}: \fk\mod^\heartsuit \rightarrow \Dmod(X)^\heartsuit, \quad 
	   \quad M \mapsto \sD_X \underset{U(\fk) } \otimes M.
	\]
	We claim that, for any algebraic subgroup $H \subset K$, this restricts to an 
	adjunction 
	\begin{equation} \label{e:locadj}
	  \on{Loc}:(\fk, H)\mod^\heartsuit \leftrightarrows 
	  \Dmod(X)^{H, \heartsuit}: \Gamma(X, -).
	\end{equation}
	In particular, specializing to $K = G$ and $X = \sB$, this recovers Theorem 
	\ref{t:heqv}.

	\subsubsection{}
	To see why \eqref{e:locadj} holds, we will need to first revisit 
	the definition of an equivariant 
	D-module. Let us follow the notation of Section \ref{s:eqdmod}. For a 
	general group $H$, not necessarily connected, as part of the data in 
	the 
	definition of an equivariant D-module one asks 
	for a distinguished isomorphism 
	\[
	\tau:  \alpha^! M \simeq \pi^! M,
	\]
	which is unital and satisfies a cocycle, i.e., associativity, condition on 
	$H 
	\times H \times 
	X$. At the 
	level of stalks, the cocycle condition asks that the two ways to identify
	\[
	i^!_{h_1 h_2 x} M \simeq i_x^! M
	\]
	by applying $\tau$ either once to the pair $(h_1h_2, x)$, or twice to the 
	pair 
	$(h_2, x)$ and then 
	$(h_1, h_2x)$, coincide. It is a nice fact that, if $H$ is connected, such a 
	datum is unique if it exists, and its existence is equivalent to the a 
	priori 
	weaker statement \eqref{e:wholething}.

		\begin{re}
	For a disconnected group, equivariance truly is a further structure one 
	equips 
	a 
	D-module with. For instance,  one can check that if $X = \on{pt}$, then a 
	datum of $H$-equivariance on an object $W$ of 
	\[
	\Dmod(\on{pt})^\heartsuit \simeq \Vect^\heartsuit 
	\]
	is the same as a representation of the component group $\pi_0(H)$ on $W$. 
	
 The reader may recognize such representations from Example 
	\ref{e:oneorbit}. This is not a coincidence, and indeed we have 
		\[ \Dmod(H/K)^{H, \heartsuit} \simeq \Dmod(H)^{H \times K, \heartsuit} 
		\simeq \Dmod(\on{pt})^{K, \heartsuit}.  \]
	In terms of D-modules on stacks, this follows from the identity $$H 
	\backslash H/K \simeq \on{pt}/K.$$	
		
	\end{re}

	\label{s:rigideq}
	
	\subsubsection{} \label{s:rigimon} The rigidified definition of equivariance in Section 
	\ref{s:rigideq} may remind the reader of the definition of an algebraic 
	$H$-module, i.e., an $\sO_H$-comodule. This is a vector space $W$ 
	equipped with a coaction map 
	\[
	   W \rightarrow \sO_H \otimes W,
	\]
	which satisfies unitality and cocycle conditions.

	\begin{re}In case the reader has not encountered this definition before, in the present setting it reduces to the following, perhaps more familiar description. Namely, $W$ carries an action of the abstract group $H(k)$, is a union of finite dimensional representations $M_\alpha$, and, with respect to any basis of such a subrepresentation $M_\alpha$, the matrix coefficients of the action of $H(k)$ are regular functions on $H$. 
	\end{re}

	Thinking of $W$ as a 
	quasicoherent sheaf on $\on{pt}$, this may be rephrased as the datum of an 
	$H$-equivariant quasicoherent sheaf on $\on{pt}$, i.e., descent along 
	$\on{pt} \rightarrow \on{pt/H}$.

	With this in mind, here is a reformulation of the data of 
	$H$-equivariance for a D-module on an $H$-variety $X$. First suppose for 
	simplicity that 
	$X$ is affine, and 
	for a D-module $M$ let us denote by the same letter the corresponding 
	$D_X$-module, i.e., its global sections.
	
	With this, a datum of $H$-equivariance 
	on $M$ is the same as 
	giving an algebraic action of $H$ on the underlying vector space of $M$, (i) 
	which 
	is compatible with the action of $H$ on $D_X$ (i.e., one upgrades from 
	a 
	$D_X$-module to a $H \ltimes D_X$-module), and (ii) such that the induced 
	infinitesimal 
	action of $\fh$ on $M$ agrees with the action by vector fields, i.e., with the 
	composition
	\begin{equation} \label{e:hch}
	\fh \rightarrow T_X \rightarrow D_X. 
	\end{equation}
	I.e., one specifies an integration of the action of $\fh$ from the 
	connection 
	on $M$ to an action of the group $H$.
	
	For general $X$, not necessarily affine, 
	one can given a similar definition. Explicitly, equivariance data for a 
	D-module is 
	(i) equivariance data for the underlying quasicoherent sheaf that is compatible 
	with 
	the canonical equivariance data on the quasicoherent sheaf of algebras 
	$\sD_X$, 
	which satisfies the analogue of condition (ii) above. 
	
	\subsubsection{} From this formulation, it is clear that global sections 
	restricts to a 
     functor 
     \[
         \Dmod(X)^{H, \heartsuit} \rightarrow (\fk, H)\mod^\heartsuit. 
     \]

	Let us also sketch why localization restricts to a functor 
	\[
	   (\fk, H)\mod^\heartsuit \rightarrow \Dmod(X)^{H, \heartsuit}.
	\]
	The starting point here is that $\sD_X$ is naturally a $H$-equivariant 
	quasicoherent sheaf of algebras -- we are simply saying functions and vector 
	fields may be moved along the $H$ action. However, $\sD_X$ is {\em not} an 
	$H$-equivariant D-module. The problem stems from condition (ii), namely the 
	coincidence of the two infinitesimal actions of $\fh$. 
	
	To see this, note 
	that the action of $\fh$ on $\sD_X$ induced by its left action on itself 
	takes the form $h \otimes P \mapsto h \cdot P$, where $h$ is an element of 
	$\fh$ 
	and $P$ is a local section of $\sD_X$. On the other hand, the induced 
	infinitesimal action of $\fh$ on $\sD_X$ obtained by $H$-equivariance must 
	be by 
	derivations of the algebra, so it cannot agree with the left action. Accordingly, it is instead the adjoint 
	action, i.e., takes the form 
	\[
	     h \otimes P \mapsto h\cdot P - P\cdot h.
	\]

      As the issue lies in the appearance of right multiplication, it should 
      not surprise the reader that tensoring on the right with an 
      $H$-integrable module should help. Indeed, given a $(\fk, H)$-module 
      $M$, it is straightforward to see that the underlying quasicoherent sheaf 
      of $$\on{Loc}(M) = \sD_X \underset{U(\fk)} \otimes M$$inherits a datum of 
      $H$-equivariance from the action of $H$ on its tensor factors. For 
      condition (ii), we may note that the derivative of this action takes the 
      form
      \[
          h \otimes (P \otimes m) \mapsto h \cdot P \otimes m - P \cdot h 
          \otimes m + P \otimes h \cdot m,
      \] 
where $h$ and $P$ are as before and $m$ is an element of $M$. But now the 
second and third factors cancel, since we tensored over $U(\fk)$, and this 
leaves only the desired first term.

	\subsection{Hidden Symmetries}
	
	\subsubsection{} In the previous subsection, we saw that localization 
	satisfied an important compatibility -- it interchanged equivariant 
	D-modules and Harish-Chandra modules. After some 
	unwinding, the assertion boiled down to a rather concrete and 
	elementary observation about the form of the functors.
	
	However, it is important to recognize that such a statement, especially at 
	the time of its introduction, did not fit easily in a known general 
	formalism in representation theory.

	We now would like to explain that this compatibility is in fact a shadow of 
	 a richer structure present on both sides of localization.
	
	 On the 
	geometric side, the further symmetry will be, if not obvious, not 
	wholly surprising. However, on the representation theoretic side, it is far 
	less obvious. It is reasonable to think of it as a hidden symmetry of the 
	category of Lie algebra representations. As we will explain, however, 
	pieces of it were witnessed and used to great effect on a case by case 
	basis, but the full structure did not emerge for several decades. 
	
	\subsubsection{} Our starting point will be the analogy between 
	\[
	   \Dmod(\sB) \quad \text{and} \quad \Fun(\sB(\Fq)). 
	\]
	Recall the latter carried an action of $G(\Fq)$, which was a linearization  
	of the action of $G(\Fq)$ on $\sB(\Fq)$ itself.

	This suggests the following idea -- we should view $\Dmod(\sB)$ as some sort 
	of representation of the group $G$, again linearizing the action of $G$ on 
	$\sB$. Crucially, the localization equivalence 
	\[
	    \Dmod(\sB) \simeq \fg\mod_0
	\]
	should be one of $G$-modules, and 
	in particular the category of Lie algebra representations should again be a 
	$G$-module. To make this precise, one needs a good notion of an action of 
	an algebraic group not on a vector space, but on a category.

	\subsubsection{} Let us run a plausibility check on this idea in the 
	somewhat more general setting of $K$ and $X$ as in 
	Section \ref{s:kandx}. 
	
	Recall the localization and global sections functors
	\[
	     \on{Loc}: \fk\mod \rightleftarrows \Dmod(X): \Gamma(X, -).
	\]
	Whatever an action of $K$ on a category $\sC$ should be, certainly for each 
	element $g$ of $K$ we should have an automorphism of the category
	$$
	    g: \sC \xrightarrow{\sim}\sC,
	$$
	and these should be suitably compatible under composition. Let us 
	now find such automorphisms of $\fk\mod$ and $\Dmod(X)$ in such a way that they are intertwined by localization.

	 We begin with the geometric side. Since $X$ is a $K$-variety, acting by 
	 $g$ gives an automorphism of the variety itself $$g: X 
	 \xrightarrow{\sim} X$$and in particular of its 
	 category of D-modules 
	 \[
	      g_*: \Dmod(X) \xrightarrow{\sim} \Dmod(X).
	 \]
	For example, for any point $x$ of $X$, this sends $\delta_x$ to 
	$\delta_{gx}$. Note that there are, for any $g$ and $g'$ in $K$, canonical 
	isomorphisms 
	\[
	     g_* \circ g'_* \simeq (gg')_*,
	\]
	i.e., one has a compatibility between multiplication in the group and 
	composition of endofunctors. 
	
	To get a hint of what the corresponding operators on $\fk\mod$ should be,  let us see what this does to global sections. Note that, for any D-module 
	 $M$ on $X$, the underlying vector spaces of 
	\[
	   \Gamma(X, M) \quad \text{and} \quad \Gamma(X, g_* M)
	\]
	are canonically identified, but the action of global differential operators 
	is twisted by $g$. 
	
	   In view of this, let us denote by 
	   \[
	       \on{Ad}_{g, *}: \fk\mod \xrightarrow{\sim} \fk\mod  
	   \]
	the automorphism which sends a $\fk$-module 
	$N$ with action $$\fk \otimes N \rightarrow N, \quad X \otimes n \mapsto X 
	\cdot n$$to the $\fk$-module with 
	the same underlying vector space and action $X \otimes n \mapsto 
	\on{Ad}_{g^{-1}}(X) \cdot n$. There again is a compatibility 
	\[
	   \on{Ad}_{g, *} \circ \on{Ad}_{g', *} \simeq \on{Ad}_{gg', *}.
	\]

	\subsubsection{} By the preceding discussion, we obtain, for every $g$ in 
	$K$, 
	canonically commuting diagrams 
	\[
	\xymatrix{\Dmod(X) \ar[d]^{g_*} \ar[rr]^{\Gamma(X, -)} && \fk\mod 
	\ar[d]^{\on{Ad}_{g, *}} \\ \Dmod(X) \ar[rr]^{\Gamma(X, -)} & &
	\fk\mod. }
	\]
	That is, $\Gamma(X, -)$ intertwines the actions of $K$ constructed so far 
	on either category. The reader may wish to check that the same holds for 
	$\on{Loc}$. 	

	\subsubsection{}Now comes an important point. The amount of structure we 
	have expected so far from a categorical representation of an algebraic group is clearly not 
	enough. For example, we have only used its closed points viewed as an abstract group.

	To proceed, we recall that the action of $G(\Fq)$ on $\Fun(\sB(\Fq))$ 
	extends to an action of the group algebra $\Fun(G(\Fq))$. Following the 
	functions-sheaves correspondence, we find that the derived category 
	$\Dmod(K)$ is monoidal under convolution. Explicitly, 
	given 
	two sheaves $\sF$ and $\sG$ in $\Dmod(K)$, their convolution is given by 
	\begin{equation} \label{e:convalgdmod}
	  \sF \star \sG := m_*\big( \hspace{.2mm} p_1^! \hspace{.2mm}\sF 
	  \hspace{.2mm}\overset{!}{\otimes}\hspace{.2mm} p_2^! \hspace{.2mm}\sG 
	  \hspace{.2mm} 
	  \big),
	\end{equation}
	where we denote by the letters $m, p_1, p_2: K \times K \rightarrow K$  
	the multiplication 
	map and projections onto the first and second factors, respectively, and we 
	use the evident functoriality for morphisms. 
	
	We 
	emphasize that formula \eqref{e:convalgdmod} matches the usual multiplication in the group 
	algebra under the functions-sheaves correspondence. Relatedly, the monoidal 
	unit is the delta function $\delta_e$ supported at the identity of the 
	group. Somewhat heuristically, one can think of a general D-module as a 
	`continuous linear combination' of delta D-modules, as in the case of the 
	usual group algebra.

	Similarly, the category $\Dmod(X)$ is naturally a module category for 
	$\Dmod(K)$. Explicitly, given a sheaf $\sF$ in $\Dmod(K)$ and $\sM$ in 
	$\Dmod(X)$, their convolution is given by 
	\begin{equation} \label{e:conv}
	     \sF \star \sM := a_*\big( \hspace{.2mm}p_1^! \hspace{.2mm}\sF 
	     \hspace{.2mm} 
	     \overset{!}{\otimes} \hspace{.2mm} p_2^! \hspace{.2mm}
	     \sM\hspace{.2mm}\big),\end{equation}
	where $a, p_1, p_2: K \times X \rightarrow X$ denote the action and 
	projection maps. 
	
	By construction, convolution with the delta function  supported at a point 
	$g$ of $K$ reproduces 
	the previously discussed symmetry, i.e., one has a canonical isomorphism of 
	functors 
	\begin{equation}
	\label{e:action}
	 \delta_g \star (-) \simeq g_*: \Dmod(X) \rightarrow \Dmod(X). 
	\end{equation}
	Informally, we have therefore extended the action from individual delta 
	functions to all their linear combinations, i.e., $\Dmod(K)$.

	Our next task is to find a corresponding action of $\Dmod(K)$ on Lie 
	algebra representations. Before doing so, we would like to 
 make some further comments on the D-module side.

	\begin{re} Note that while the discussion of the automorphisms $g_*$ and 
	$\on{Ad}_{g, *}$ worked equally well with abelian or derived categories, 
	for this monoidal structure to behave in the desired way, i.e., to even have a reasonable D-module pushforward, it is essential we 
	work with derived categories. So far, the corresponding triangulated 
	category is enough, but we will be forced to ask for more in the next section. 
	
	Note also that to be closed under convolution we immediately must leave the 
	world of finitely generated D-modules. The prototypical example is 
	the convolution of $\sD_K$ with itself, namely one can show that
	\[
	    \sD_K \star \sD_K \simeq \sD_K \underset{\underline{k}}{\otimes} \OO_K.
	\]
	\end{re}

	\begin{re} It is clarifying to rewrite the formulas for convolution and 
	action as
		\begin{align*}
		&\Dmod(K) \otimes \Dmod(K) \simeq \Dmod(K \times K) \xrightarrow{m_*} 
		\Dmod(K) \quad \text{and} \\ &\Dmod(K) \otimes \Dmod(X) \simeq 
		\Dmod(K \times X) \xrightarrow{a_*} \Dmod(X),
		\end{align*}
	exactly as in the function-theoretic case, cf. Equations \eqref{e:convfun} 
	and \eqref{e:actionfun}. However, to do so we need the tensor product of 
	derived categories, which can only be performed at the level of 
	dg-categories, as discussed in Section \ref{s:catrepgrp}.  	
	\end{re}

	\subsubsection{}\label{ss:weakinv} Finally, let us turn to the categorical action of $K$ on 
	$\fk\mod$. To begin, note that the construction of the present action is 
	less immediate 	than in the case of $\Dmod(X)$.  Indeed, there we only 
	needed D-modules to act	on D-modules, for which we could use geometric 
	functorialities. Here, we need to have D-modules act on Lie algebra 
	representations, which are seemingly more remote from one another.

	However, we may proceed by rewriting the category $\fk\mod$ as a category 
	of D-modules as follows. Recall that if $K$ acts on a variety $X$, we saw 
	that the sheaf of differential operators $\sD_X$ was not an equivariant 
	D-module. However, it is a prototypical example of a {\em weakly
	equivariant} D-modules -- the abelian category of such objects, which we 
	denote by 
\[
     \Dmod(X)^{K, w, \heartsuit}
\]
 is defined exactly as in Section \ref{s:rigimon}, but without the crucial condition 
 (ii) imposed.

If we take $X = K$, viewed as a $K$-variety under right multiplication, we 
claim that one has a canonical equivalence 
\begin{equation} \label{e:weakinvdmod}
    \Dmod(K)^{K, w, \heartsuit} \simeq \fk\mod^\heartsuit. 
\end{equation}
Let us sketch a proof for the convenience of the reader. First, one checks that 
the left hand side has projective generator $\sD_K$. This follows from the fact 
that it corepresents $K$-invariant global sections, i.e., $\Gamma(K, -)^K$, combined with the tautological isomorphism between $K$-equivariant quasicoherent sheaves on 
$K$ and $\on{Vect}$. Second, one concludes by recalling the identification of 
right invariant differential operators and the enveloping algebra, i.e., 
\[
     \Gamma(K, \sD_K)^K \simeq U(\fk). 
\]
Summing up, we have found a projective generator with endomorphisms $U(\fk)$, so we are done by Morita theory.

The equivalence of Equation \eqref{e:weakinvdmod} prolongs to an equivalence of 
derived categories 
\[
    \Dmod(K)^{K, w} \simeq \fk\mod,
\]
and the former carries a $\Dmod(K)$ action via left convolution, which 
completes the construction of the desired action.

\subsubsection{} With this, the desired upgrade of the localization and global 
section functors reads as follows.

\begin{theo}[Beilinson--Drinfeld, \cite{bdh}]\label{t:loc+equiv}  Let $X$ and $K$ be as in Section \ref{s:kandx}. There is a canonical adjunction of 
$\Dmod(K)$ equivariant functors 
	\[
	    \on{Loc}: \fk\mod \rightleftarrows \Dmod(X): \Gamma(X, -). 
	\]
\end{theo}

In particular, Beilinson--Bernstein localization, as stated in Theorem \ref{t:bbloc}, is an equivalence 
of $\Dmod(G)$ representations.\footnote{The careful reader will note that this localization theorem concerned not all of $\fg\mod$ but rather $\fg\mod_0$. The $\Dmod(G)$ action on $\fg\mod_0$ formally follows from the appropriate commutation of the $Z(\fg)$ and $\Dmod(G)$ actions on $\fg\mod$, or can be constructed by hand as in Section \ref{ss:weakinv}.} We defer an explanation of Theorem \ref{t:loc+equiv}, including a discussion of $\Dmod(K)$-equivariance, 
until Section \ref{s:catrepgrp}. 

\subsubsection{} Let us indicate the motivating application of this theory. 
To orient ourselves, consider a subgroup of a finite group, which we suggestively denote by $H(\Fq) \subset K(\Fq)$. Recall that, for any 
representation $\pi$ of $K(\Fq)$, its $H(\Fq)$ invariants $\pi^{H(\Fq)}$ carry 
an action of the 
convolution algebra $$\Fun(H(\Fq) \backslash K(\Fq) /H(\Fq)).$$

In a similar way, given an algebraic subgroup $H$ of $K$, for any categorical 
representation $\sC$ of $K$ one can form its $H$-invariant objects 
$\sC^H$, and 
these will carry an action of the convolution category 
\[
   \Dmod(H \backslash K  / H).
\]
In the current example, we obtain the following upgrade of Theorem \ref{t:heqv} 
and its derived variants.

\begin{theo}[Beilinson--Drinfeld \cite{bdh}]\label{t:lochchmod+eq}There is a canonical adjunction of $\Dmod(H 
\backslash K / 
H)$-equivariant functors 
	\[
	\on{Loc}: \fk\mod^H \leftrightarrows \Dmod(X)^H: \Gamma(X, -). 
	\]
\end{theo}

Once the formalism is in place, this is deduced from Theorem \ref{t:loc+equiv} 
by applying $H$-invariants. In particular, the general notion of categorical 
invariants, which we will discuss in Section \ref{s:catrepgrp}, recovers on the 
Lie theoretic side the category of Harish-Chandra modules
\[
   (\fk, H)\mod \simeq \fk\mod^H
\]
and on the geometric side the category of equivariant D-modules. 

\subsubsection{}\label{ss:hitchinhecke}
Theorem \ref{t:lochchmod+eq} was proven by Beilinson--Drinfeld in their fundamental work \cite{bdh} on quantization of the Hitchin system.\footnote{Strictly speaking, they needed 
and proved a certain infinite-dimensional analog.} In this motivating problem, 
the 
compatibility with the action of $\Dmod(H \backslash K / H)$ was crucial. 
Namely, the desired result was the existence of Hecke eigensheaves, certain D-modules with 
equivariance properties with respect to this action. 
They were constructed via localization from Lie algebra representations, 
which were shown directly to have the desired property.

\subsubsection{}Here is a slightly less sophisticated example. Consider the 
category of highest weight representations 
\[
    \fg\mod^B.
\]
By the previous formalism, it tautologically carries an action of the Hecke 
category 
\[
   \Dmod(B \backslash G / B). 
\]
We stress that, in its historical development, the connection between highest weight modules and D-modules on the flag variety was seen through the localization theorem and, in particular, was highly nonobvious. However, the formalism of categorical actions makes their relation manifest. Moreover, this quickly and conceptually reproduces several important endofunctors of 
this category.

 For example, for each element $w$ of the finite Weyl group, one 
attaches certain endofunctors 
\[
  j_{w, !} \quad \text{and} \quad j_{w, *}: \fg\mod^B \rightarrow \fg\mod^B,
\]
known as the Enright completion and Arkhipov twisting functors, respectively \cite{enright79}, \cite{bb82}, \cite{ark}. 
Geometrically, these are given by convolution with the $!$- and $*$-extensions 
of the constant D-module on 
\[
     BwB \subset G, 
\]
and many of their basic properties have simple geometric proofs.

As a related example, modulo nontrivial facts about the intersection cohomology 
of Schubert varieties, the Kazhdan--Lusztig conjecture follows from the assertion 
that 
\begin{equation} \label{e:brungus}
  \fg\mod_0^B
\end{equation}
is canonically equivalent to $\Dmod(B \backslash G / B)$. This was originally 
proven by establishing the localization theorem. However, using the action of 
the Hecke category on \eqref{e:brungus}, one can also directly check that 
convolution with the antidominant Verma module yields an equivalence
\[
  \Dmod(B \backslash G / B) \simeq \fg\mod_0^B, 
\]
from which the fully faithfulness of $\Gamma: \Dmod(\Fl_G) \rightarrow \fg\mod_0$ formally follows, cf. \cite{ahc}.

	\section{Categorical representations of groups}
	\label{s:catrepgrp}

	\subsection{Overview}

\subsubsection{} In Section \ref{s:bbloc} we followed the symmetries of 
Beilinson--Bernstein localization and arrived at the notion of a categorical 
representation of a group. In this section, we would like to give a more careful introduction to this formalism.

The basic ideas and results here are due to Beilinson--Drinfeld and Frenkel--Gaitsgory. The formalism for discussing the category of all categorical representations and the intertwining functors between them is due largely to Gaitsgory. 

 This formalism requires some finer points of modern homological algebra to get off the ground, so that is where we begin.

\subsection{From triangulated to dg-categories}

\subsubsection{} In Section \ref{s:bbloc}, we treated $\Dmod(K)$ as a monoidal 
triangulated category, and its modules were triangulated categories. However, 
to perform many basic manipulations with its modules, the 
triangulated structure is insufficient. 

Here is a basic problem of this form. Given modules $\sC_1$ and $\sC_2$ for 
$\Dmod(K)$, one would like to speak about the category of $\Dmod(K)$ equivariant 
triangulated functors 
\[
   \Hom_{\Dmod(K)\mod}(\sC_1, \sC_2).
\] 
Moreover, we would like these to behave like their function theoretic 
counterparts. For example, given subgroups $H_1(\Fq)$ and $H_2(\Fq)$ of a 
finite group $K(\Fq)$, one has a canonical isomorphism
\begin{align*}
  \Hom_{K(\Fq)\mod}\big(\Fun(K(\Fq)/H_1(\Fq)),\Fun(K(\Fq)/H_2(\Fq)) 
  \big) \simeq 
 \Fun\big( H_1(\Fq) \backslash K(\Fq) / H_2(\Fq)\big).
\end{align*}
One would like, for algebraic subgroups $H_1$ and $H_2$ of $K$, a similar 
isomorphism 
\begin{align*}
\Hom_{\Dmod(K)\mod}\big(\Dmod(K/H_1),\Dmod(K/H_2) 
\big)  
\simeq 
\Dmod( H_1 \backslash K / H_2).
\end{align*}

Unfortunately, the construction of suitable Hom spaces is an issue already if 
$K$ is the trivial group. In this 
case we are asking for a reasonable category of triangulated functors from 
$\sC_1$ to $\sC_2$. However, this is a well known issue with triangulated 
categories. For example, given triangulated functors 
equipped 
with a morphism $F_1 \rightarrow F_2$, it is not clear how to complete it to an 
exact triangle. Namely, one would like to take cones pointwise, i.e., for each object $c$ of $\sC_1$ set$$\on{Cone}(F_1 \rightarrow F_2)(c) \simeq \on{Cone}(F_1(c) \rightarrow F_2(c)),$$but turning this into a functor runs into trouble due to the non-functoriality of cones in $\sC_2$.

	As a more sophisticated example, one would like to be able to take limits, 
	colimits, and tensor products of $\Dmod(K)$ representations, which again are not 
	readily defined when working with triangulated categories.

	Happily, all these issues may be addressed by passing from triangulated to 
	dg-categories.

	\subsubsection{} Let us quickly explain the basic idea of dg-categories. 
	Recall that for most triangulated categories appearing in 
	representation theory, such as $\fk\mod$, the homomorphisms between two 
	objects $t_1$ and $t_2$ are computed as the zeroth cohomology of a complex 
	of vector spaces. Similarly, the homomorphisms between $t_1$ and the shift  
	$t_2[i]$, for any integer $i$, is the $i^{th}$ cohomology of the same 
	complex.
	
	Concretely, for the derived category of an abelian category, one takes 
	either a projective resolution $p_1$ of $t_1$ or an injective resolution 
	$i_2$ of 
	$t_2$\footnote{Of course, strictly speaking, since we are talking about 
	unbounded complexes for an abelian category possibly of infinite 
	cohomological 
	dimension, we mean K-projective and K-injective resolutions in the 
	sense of Bernstein and Spaltenstein \cite{spaltenstein}.} and uses the corresponding complex of 
	homomorphisms 
	\[
	  \underline{\Hom}(p_1, t_2) \xrightarrow{\sim} \underline{\Hom}(p_1, i_2) \xleftarrow{\sim} 
	  \underline{\Hom}(t_1, i_2).  
	\]
	We recall that the the above identifications are quasi-isomorphisms, i.e., induce isomorphisms on cohomology, but typically do not identify the entire chain complexes themselves.

	The basic idea of dg-categories is to remember not only the cohomology of 
	these complexes but also all the complexes themselves, up to simultaneous quasi-isomorphism. The reader familiar with 
	derived categories should not find this maneuver so surprising. 
	After all, derived categories arise from the usefulness of remembering not only individual derived functors, e.g., $\Ext^1$, but also the complex computing them. 
	Passing from triangulated to dg-categories simply repeats this idea not 
	only for objects, but also for morphisms. 
	 
  \subsubsection{} Plainly, a $k$-linear dg-category $\sC$ has a collection of 
  objects. 
  For 
  each pair of objects $c_1$ and $c_2$, one has a complex of $k$ vector spaces  
  $\Hom(c_1, c_2)$, and for each triple of objects one has a composition 
  morphism of complexes
  \[
  \Hom(c_1, c_2) \otimes \Hom(c_2, c_3) \rightarrow \Hom(c_1, c_3), 
    \]
	which is associative in the naive sense.

	Here are some basic operations within a dg-category. First, given an object 
	$c$ of $\sC$ and an integer $i$, the shift $c[i]$, if it exists, 
	corepresents the functor 
	\[
	\Hom(c[i], -) \simeq \Hom(c,-)[-i]. 
	\]
	Similarly, given a 0-cycle $f: c_1 \rightarrow c_2$, the cone of $f$, if it 
	exists, corepresents the functor 
	\[
	   \Hom( \on{Cone}(f), -) \simeq \on{Cone}( \Hom(c_2, -) \rightarrow 
	   \Hom(c_1, -))[-1],
	\]
	where the cone on the right is the usual one of complexes of vector spaces. 
	Finally, direct sums are defined in the usual way, i.e., when they exist 
	they corepresent the functor
\[
   \Hom( \underset{i}{\oplus} \hspace{.2mm} c_i, -) \simeq \underset{i}\Pi 
   \hspace{.2mm} \Hom( c_i, -).
\]
  A dg-category is said to be cocomplete if it is closed under the above 
  operations, i.e., shifts, cones, and arbitrary direct sums.\footnote{In what 
  follows, 
  we also tacitly assume that any cocomplete dg-category is also presentable. This is a technical assumption about its size satisfied in all relevant examples.} 

  \subsubsection{} Here are two basic examples of cocomplete dg-categories. 
  
  \begin{ex} The category $\Vect$ of $k$ vector spaces is a cocomplete 
  dg-category. Explicitly, objects of 
  $\Vect$ are chain complexes of $k$ vector spaces $(V_i, d_i)$, for $i \in 
  \mathbb{Z}$. 
  
  Given a pair of complexes 
  $(V_i, d_V)$ and $(W_i, d_W)$, we have 
  \[
  \Hom(V_1, V_2) := \underline{\Hom}(V_1, V_2),
  \]
  where $\underline{\Hom}(V_1, V_2)$ is the chain complex with $i$-chains 
  \[
      \underline{\Hom}(V, W)^i := \underset{j}{\Pi} \hspace{.2mm} \Hom(V_j, 
      W_{i+j}),
  \]
  and the usual differential
  \[
       d \phi := d_W \circ \phi -(-1)^i \phi \circ d_V, \quad \text{for } \phi 
       \in \underline{\Hom}(V,W)^i. 
  \]
  In particular, $i$-cycles are simply maps of chain complexes $V \rightarrow 
  W[i]$, and $i$-boundaries are maps which are null-homotopic. \label{e:dgvect}
  \end{ex}

\begin{ex}More generally, suppose $\sA$ is a $k$-linear Grothendieck abelian 
category, such as D-modules 
on an algebraic variety over $k$ or the left modules for a $k$-algebra. In this case, the 
unbounded derived category 
$
    \on{D}(\sA)
$
naturally enhances to a dg-category as follows. The objects are given by 
K-injective complexes of objects of $\sA$, and the homomorphisms are defined  
as in Example \ref{e:dgvect}.

\end{ex}

\subsubsection{}\label{s:homdgcats} Given two cocomplete dg-categories $\sC_1$ 
and $\sC_2$, one 
has a cocomplete dg-category 
\begin{equation} \label{e:bigbrainhoms}
\udHom(\sC_1, \sC_2),
\end{equation}
which one thinks of as the functors between them; the subscript `$\DGCat$' will be discussed further in Section \ref{s:dgcat} below. The construction of \eqref{e:bigbrainhoms} is not 
completely straightforward, as one needs it to be invariant under 
quasi-equivalences of dg-categories.\footnote{A naive dg-functor between 
dg-categories $F: \sC \rightarrow \sC'$ is specified by compatible maps between 
the sets of objects and chain complexes of morphisms. It is said to be a 
quasi-equivalence if (i) it induces quasi-isomorphisms on the complexes of 
homomorphisms and (ii) for any object $c'$ of $\sC'$, there exists an object $c$ of 
$\sC$ and a 0-cycle $F(c) \rightarrow  c'$ inducing quasi-isomorphisms 
\begin{equation} \label{e:bigbrainboms}
  \Hom_{\sC'}(-, F(c)) \simeq \Hom_{\sC'}(-, c'). 
\end{equation}
This may be reformulated as follows. After passing to the zeroth cohomology of all complexes of 
homomorphisms, a naive dg-functor $F$ induces a functor between the associated `homotopy' categories
\[
\on{Ho}(F): \on{Ho}(\sC) \simeq \on{Ho}(\sC'). 
\]
Under the assumption that $\sC$ and $\sC'$ are closed under shifts, $F$ is a quasi-equivalence if and only if $\on{Ho}(F)$ is an equivalence of ordinary categories. 

Plainly, in practice one does not care about the strict complexes of homomorphisms themselves, but instead the complexes up to quasi-isomorphism. This leads to both the notion of quasi-equivalence and the desire to obtain functor categories invariant under quasi-equivalence. 
}

 However, it may be 
constructed, as with 
maps in the derived category, by model category theoretic means as a naive 
category of dg-functors between suitable replacements of $\sC_1$ and $\sC_2$; see \cite{toenlectures} and references therein.  

Let us describe some basic properties of \eqref{e:bigbrainhoms}. First, the set of objects 
of \eqref{e:bigbrainhoms} are the quasi-functors in the sense of Keller \cite{kellerderiving} which are moreover continuous, i.e. commute with 
arbitrary direct sums. Second, the set of $i$-cycles in the morphisms between $F_1$ 
and $F_2$ are, up to (co)fibrant replacement, natural transformations $F_1 
\rightarrow F_2[i]$. 

\begin{re} We should mention that the adjective `continuous' here refers to commutation with colimits, whereas in other parts of the literature it refers to commutation with limits. \end{re}

Here are some important examples of such functor categories.

\begin{ex} For any cocomplete dg-category $\sC$, one has a canonical 
equivalence 
	\[
	     \udHom(\Vect, \sC) \simeq \sC,
	\]
	which exchanges an object $c$ of $\sC$ with the dg-functor of tensoring 
	with $c$, i.e., $V \mapsto V \otimes c$. 
	
	Explicitly, this dg-functor assigns a vector space $k^{\oplus I}$, for a set $I$, the object $c^{\oplus I}$, and given a map $k^{\oplus I} \rightarrow k^{\oplus J}$, the resulting map $c^{\oplus I} \rightarrow c^{\oplus J}$ is obtained from the map $k \rightarrow \on{Hom}_\sC(c, c)$, i.e. the $k$-linearity of $\sC$. The image of a bounded complex of vector spaces is obtained by taking cones, and for an unbounded complex by taking colimits. 
\end{ex}

\begin{ex}The present example, which is a version of the function theoretic discussion in Section \ref{sss:endos}, will be fundamental in what follows and is due to \cite{toenpaper} for varieties and  extended by \cite{benzvifrancisnadler} to a broad class of derived stacks.

For algebraic varieties $X$ and $Y$, one has a canonical equivalence 
\[
   \Hom_{\DGCat}(\QCoh(X), \QCoh(Y)) \simeq \QCoh(X \times Y),
\]
where one associates to sheaf $\mathscr{K}$ on the product the corresponding integral transform, i.e., 
\[
   \QCoh(X) \xrightarrow{\pi_X^*} \QCoh(X \times Y) \xrightarrow{\mathscr{K} \otimes - } \QCoh(X \times Y) \xrightarrow{\pi_{Y, *}} \QCoh(Y), 
\]
where $\pi_X: X \times Y \rightarrow X$ and $\pi_Y: X \times Y \rightarrow Y$ denote the projections. 

Similarly, one has a canonical equivalence 
\[
   \Hom_{\DGCat}(\Dmod(X), \Dmod(Y)) \simeq \Dmod(X \times Y),
\]
which again sends a D-module $\mathscr{K}$ on the product to the integral transform 
\[
   \Dmod(X) \xrightarrow{\pi_X^!} \Dmod(X \times Y) \xrightarrow{ \mathscr{K} \overset ! \otimes - } \Dmod(X \times Y) \xrightarrow{\pi_{Y, *}} \Dmod(Y). 
\]
\end{ex}

\subsubsection{}\label{s:dgcat} To talk about (co)limits of 
dg-categories, one uses the following structure. There exists an $(\infty, 
1)$-category, which we denote by $\DGCat$, whose objects are 
cocomplete dg-categories and 1-morphisms are continuous quasi-functors; see Chapter 1 of \cite{gaitsroz} and references therein. Here are 
some of its basic properties. 

First, it admits all (co)limits in the 
$\infty$-categorical sense. Briefly, an object of a limit of dg-categories 
$\varprojlim \sC_\alpha$ is a homotopy-coherent system of objects in each 
$\sC_\alpha$, and morphisms are the corresponding inverse limit of mapping 
complexes.\footnote{We highlight that this is again a basic point where 
triangulated 
categories would not suffice. That is, one really needs the complexes, and not 
simply their cohomology, due to the non-exactness of inverse limits.} 

Colimits of diagrams of dg-categories are then characterized by the identity 
\[
\udHom(\varinjlim \sC_\alpha, -) \simeq \varprojlim 
\udHom(\sC_\alpha, -)
\]
and, relatedly, limits of dg-categories satisfy the identity 
\[
 \udHom(-, \varprojlim \sC_\beta) \simeq \varprojlim 
 \udHom(-, \sC_\beta).
\]

It is not easy to describe the objects or morphisms in a general colimit of 
dg-categories (and indeed, one in general should only expect straightforward 
descriptions for limits, as they are characterized by a mapping in 
property). However, to orient the reader, let us describe an important example 
that admits an explicit description.

\begin{ex}\label{e:colimitcats}Recall that an object $c$ of a dg-category $\sC$ 
is said to be {\em 
compact} if $$\Hom(c, -)$$ commutes with colimits. These should be thought of 
as 
small objects, e.g., in $\Vect$ an object is compact if and only if it is 
equivalent to a bounded complex of finite dimensional vector spaces. Relatedly, 
$\sC$ is compactly generated if every object can be written as a colimit of 
compact objects. 

Given a filtered diagram of compactly generated categories $\sC_\alpha$, where 
each 1-morphism $\sC_\alpha \rightarrow \sC_{\alpha'}$ sends compact objects to 
compact objects, the colimit 
\[
  \varinjlim \sC_{\alpha}
\]
admits the following description. It is again compactly generated by insertions 
of compact objects from each $\sC_{\alpha}$, and homomorphisms between such 
objects are the filtered colimit of the homomorphisms at each step in the 
colimit. 
\end{ex}

\subsubsection{}A very useful feature of $\DGCat$ is that one can tensor together 
cocomplete dg-categories 
\[
   \sC \text{ and } \sD \rightsquigarrow \sC \otimes \sD. 
\]
More carefully, $\DGCat$ is a symmetric monoidal $\infty$-category with respect 
to the Lurie tensor product.

In many practical examples, the following hands-on description of the tensor 
product is useful. 

\begin{ex}Suppose that $\sC$ and $\sD$ are compactly generated, say by sets of 
objects $c_\alpha$ and $d_\beta$. Then the tensor product is compactly 
generated by the set of simple tensors $c_\alpha \otimes d_\beta$, with 
homomorphisms given by the K\"{u}nneth formula, i.e., 
\[
\Hom( c_\alpha \otimes d_\beta, c_{\alpha'} \otimes d_{\beta'}) \simeq 
\Hom(c_\alpha, c_{\alpha'}) \otimes \Hom(d_\beta, d_{\beta'}).
\]

\end{ex}

Here are some special cases of the above example, which hopefully convey to the 
reader the sense in which this tensor product aligns with their intuition from 
the functions-sheaves correspondence, cf. Section \ref{ss:opsfinitesets}.

\begin{ex} Given algebraic varieties $X$ and $Y$, one has a canonical 
equivalence of the categories of quasicoherent sheaves 
	\[
	   \on{QCoh}(X) \otimes \on{QCoh}(Y) \simeq \on{QCoh}(X \times Y).
	\]
Explicitly, given an object $\sM$ of $\on{QCoh}(X)$ and $\sN$ of 
$\on{QCoh}(Y)$, the simple tensor $\sM \otimes \sN$ of the left hand side is sent to the 
external product object $$\sM \boxtimes \sN	\simeq \pi_X^* \sM 
\overset{}\otimes \pi_Y^* \sN.$$

\end{ex}

\begin{ex} Let $X$ and $Y$ be as in the preceding example. One also has a 
canonical equivalence of the categories of D-modules 
	\[
	     \Dmod(X) \otimes \Dmod(Y) \simeq \Dmod(X \times Y).
	\]
   Again, given an object $\sM$ of $\Dmod(X)$ and an object $\sN$ of 
   $\Dmod(Y)$, the simple tensor $\sM \otimes \sN$ of the left hand side is sent to the 
   external product object 
   \[
      \sM \boxtimes \sN \simeq \pi_X^! \sN \overset{!}{\otimes} \pi_Y^! \sM. 
   \] 

\end{ex}

\begin{ex}Given $k$-algebras $A$ and $B$, one has a canonical equivalence 
	\[
	    A\mod \otimes \hspace{.2mm} B\mod \simeq (A \underset{k} \otimes 
	    B)\mod. 
	\]
	Note that this generalizes the affine cases of the previous two examples. 
\end{ex}

Let us mention that the monoidal structure is unital, with unit the category of 
vector spaces $\Vect$. In particular, for any dg-category $\sC$, one has a 
canonical equivalence 
\[
   \Vect \otimes \hspace{.7mm}  \sC \simeq \sC. 
\]
The reader may wish to check this directly in the case $\sC$ is compactly 
generated.

\subsubsection{} Given dg-categories $\sC$ and $\sD$, the category $\udHom(\sC, 
\sD)$ of Section \ref{s:homdgcats} is an internal hom object with respect to 
the monoidal structure on $\DGCat$.

A technical point is that this is not 
the same as the category of morphisms between $\sC$ and $\sD$, viewed as objects of 
$\DGCat$. Instead, the 
category of morphisms is the truncation of 
$\udHom(\sC, \sD)$ to an $\infty$-groupoid, i.e., one discards the non-invertible 
natural 
transformations. It is expected that one has a canonical enhancement of 
$\DGCat$ to a
symmetric monoidal $(\infty,2)$-category, where the homomorphisms are indeed 
$\udHom(-, -)$, but to our knowledge this is not yet available in the 
literature.

\subsubsection{} Using the monoidal structure, one can make sense of algebra 
objects of $\DGCat$, i.e., monoidal dg-categories. Informally, this 
consists of a dg-category $\sA$ with a multiplication 
\[
    \sA \otimes \sA \rightarrow \sA,
\]
which is associative up to coherent homotopy. We emphasize the coherent homotopy 
is further structure and is part of the full data of a monoidal dg-category.

 Here are some basic examples of monoidal dg-categories. 
 
 \begin{ex}\label{ex:qcohalg} For an algebraic variety $X$, $\QCoh(X)$ is a 
 monoidal dg-category 
 under $*$-tensor product, and $\Dmod(X)$ is a monoidal dg-category under 
 $!$-tensor product. 	
 
 Given a map $f: X_1 \rightarrow X_2$, the corresponding pullback functors 
 \[
 f^*: \QCoh(X_2) \rightarrow \QCoh(X_1) \quad \text{and} \quad f^!: \Dmod(X_2) 
 \rightarrow \Dmod(X_1)
 \]
 carry a canonical datum of monoidality. 
 \end{ex}

\begin{ex} Suppose a variety $Y$ is an algebra object in the category of 
varieties. That is, it is equipped with a map 
	\[
	    \mu: Y \times Y \rightarrow Y,
	\]
which is associative in the evident sense. Then $\QCoh(Y)$ and $\Dmod(Y)$ both 
acquire `convolution' monoidal structures, given by $*$-pushforward along 
$\mu$.

Given a map $f:Y_1 \rightarrow Y_2$ of algebras, the associated pushforward 
maps 
\[
 f_*: \QCoh(Y_1) \rightarrow \QCoh(Y_2) \quad \text{and} \quad f_*: \Dmod(Y_1) 
 \rightarrow \Dmod(Y_2)
\]
carry a canonical datum of monoidality. 

As our primary case of interest, we take $Y = K$ to be an algebraic 
group, with $\mu$ its multiplication, and thereby obtain the group algebra 
$\Dmod(K)$ we met previously (or rather, its canonical dg-enhancement).

Similarly, given a homomorphism $\phi: K_1 \rightarrow K_2$ of groups, one 
obtains a homomorphism of group algebras 
\[
   \phi_*: \Dmod(K_1) \rightarrow \Dmod(K_2).
\]

\end{ex}

\begin{re} One can replace in the two preceding paragraphs $\Dmod(K)$ with $\QCoh(K)$ to obtain another monoidal dg-category under convolution. Modules for $\Dmod(K)$ are sometimes called strong categorical representations and modules for $\QCoh(K)$ weak categorical representations. 

It may be clarifying to remark that, given any monoidal functor from varieties to a monoidal $\infty$-category $(\sS, \otimes)$, we obtain a similar convolution algebra linearizing $K$. In particular, one can speak of further higher categorical versions of $\on{Fun}(K(\Fq))$, e.g., by replacing sheaves with sheaves of  categories on $K$. \end{re}

\subsubsection{} Given a monoidal dg-category $\sA$, one has associated 
$\infty$-categories of its left and right modules, which we denote by 
\[
  \sA\mod \quad \text{and} \quad \rmod\sA,
\]
respectively. Somewhat informally, a left module consists of a cocomplete 
dg-category $\sM$ equipped with an action map 
\[
  \sA \otimes \sM \rightarrow \sM,
\]
which is `associative' up to coherent homotopy and similarly for a right 
module. Again, the coherent homotopy is part of the data of specifying a 
module.

Given two left modules $\sM_1$ and $\sM_2$, one has the 
dg-category of 
$\sA$-equivariant functors 
\[
   \Hom_{\sA\mod}(\sM_1, \sM_2)
\]
(and similarly for right modules over $\sA$). Explicitly this may be calculated via a bar complex as the limit of the 
semi-cosimplicial diagram of dg-categories  

\begin{equation} \label{e:hombar}
   \xymatrix{  \udHom(\sM_1, \sM_2)
   	\ar@<-.4ex>[r] \ar@<.4ex>[r] &  \hspace{.5mm} \udHom(\sA \otimes \sM_1, 
   	\sM_2) 
   	\ar@<-.8ex>[r] \ar@<.8ex>[r] \ar[r] & \hspace{.5mm} \udHom(\sA^{\otimes 2}  
   	\otimes \sM_1, \sM_2)   \hspace{.5mm} \cdots \hspace{.5mm},   
   }
   \end{equation}
where the arrows are induced in the standard pattern of bar resolutions by 
the action maps of $\sA$ on $\sM_1$ and $\sM_2$ and the multiplication map 
$\sA^{\otimes 2} \rightarrow \sA$.

  In addition, given a right module $\sM$ and a left module $\sN$, one can form 
  their tensor product 
  \[
     \sM  \underset{\sA}{\otimes} \sN. 
  \]
As for homomorphisms, this may be computed via a bar resolution as the colimit of 
the semisimplicial diagram of cocomplete dg-categories 
\begin{equation} \label{e:tensorbar}
\xymatrix{  \cdots \hspace{.5mm} \sM \otimes \sA^{\otimes 2} \otimes N  
\ar@<-.8ex>[r] \ar@<.8ex>[r] \ar[r]   &  \hspace{.5mm} \sM \otimes \sA \otimes 
\sN \ar@<-.4ex>[r] 
\ar@<.4ex>[r]  &  \hspace{.5mm} \sM \otimes  \sN}.
\end{equation}

Let us now give a few basic examples. 

\begin{ex} For any left $\sA$-module $\sM$, there are tautological equivalences 
	\[
	   \on{Hom}_{\sA\mod}(\sA, \sM) \simeq \sM \simeq \sA \underset \sA \otimes 
	   \sM
	\]
	and similarly for right modules. 
\end{ex}

\begin{ex} \label{e:fibprodvar}Given maps of varieties $X \rightarrow Z 
\leftarrow Y$, as in 
Example \ref{ex:qcohalg}, we obtain canonical actions of $\QCoh(Z)$ on 
$\QCoh(X)$ and $\QCoh(Y)$. With this, their tensor product is given by
\[
    \QCoh(X) 
    \underset{\QCoh(Z)}{\otimes} \QCoh(Y) \simeq \QCoh( X \underset{Z}{\times} 
    Y),
\]
where the above fibre product is derived, i.e., one obtains a 
derived scheme.\footnote{Plainly, this is glued from the derived rings $\OO_X 
\underset{\OO_Z}{\overset{L}{\otimes}} \OO_Y$ affine by affine, exactly as in 
the case of the usual fibre product. One thinks of the negative cohomology groups of the structure sheaf as adding `nilpotent fluff' to the variety, exactly like the perhaps more familiar the non-reducedness of a general scheme. All derived schemes we will meet in this survey arise in this way.} 
Concretely, the equivalence sends a simple tensor $\sM \otimes \sN$ to their 
external product, restricted to $X \underset{Z}{\times} Y$.  

Similarly, their homomorphisms are again given by 
\[
 \Hom_{\QCoh(Z)}(\QCoh(X), \QCoh(Y)) \simeq \QCoh( X \underset{Z}{\times} 
 Y).
\]
Explicitly, one associates to a quasicoherent sheaf $\sS$ on the fibre product 
the integral transform 
\[
  \QCoh(X) \xrightarrow{ \pi_X^*} \QCoh(X \underset Z \times Y) \xrightarrow{ - 
  \overset * \otimes \sS} \QCoh(X \underset Z \times Y) \xrightarrow{\pi_{Z, 
  *}} \QCoh(Y). 
\] 

The same assertions hold, mutatis mutandis, for D-modules. That is, we have 
equivalences 
\begin{align*}
  \Dmod(X) \underset{\Dmod(Z)}{\otimes} \Dmod(Y) &\simeq \Dmod(X \overset{\on{cl}}{{\underset Z {\times}}} Y) \text{ and}  \\ 
  \Hom_{\Dmod(Z)}(\Dmod(X), \Dmod(Y)) &\simeq \Dmod(X \overset{\on{cl}}{{\underset Z {\times}}} Y),
\end{align*}
again given by external tensor product and integral transform. Here the superscript `$\on{cl}$' denotes the classical, i.e., non-derived, fibre product. Note however that when the derived fibre product has only finitely many nonzero cohomologies in its structure sheaf, e.g., if $Z$ is smooth, due to the invariance of categories of D-modules under 
nil-thickenings it is 
unimportant whether one takes the derived or classical fiber product of 
varieties. 

As with other sheaf-theoretic statements we have met, these identities are 
analogues 
of the 
following identity for functions. Given a finite set $S$, consider $\Fun(S)$ as 
an algebra 
under 
pointwise multiplication. Given finite sets $X \rightarrow Z 
\leftarrow Y$, one then has canonical isomorphisms of vector spaces 
\[
    \Fun(X) \underset{\Fun(Z)}{\otimes} \Fun(Y) \simeq 
    \Hom_{\Fun(Z)}(\Fun(X), \Fun(Y)) \simeq \Fun(X \underset Z \times Y).
\]
\end{ex}

\subsection{Basics of categorical representation theory}

\subsubsection{} With these generalities in hand, let us describe some of the 
basic operations in categorical representation theory and thereby recover Theorems \ref{t:loc+equiv} and \ref{t:lochchmod+eq} encountered in Section \ref{s:bbloc}.

Recall that, for an algebraic group $K$, we have its associated $\infty$-category of 
categorical representations 
\[
   \Dmod(K)\mod.
\]
As elsewhere in representation theory, inversion on $K$ induces an isomorphism 
between $\Dmod(K)$ and its reverse monoidal dg-category, and in particular one 
may freely pass between left and right modules.

\subsubsection{} A basic family of representations may be constructed 
geometrically as follows. 

Given a variety $X$ with an action of $K$, one obtains a canonical action of 
$\Dmod(K)$ on $\Dmod(X)$. Explicitly, the underlying binary product is given by 
\[
   \Dmod(K) \otimes \Dmod(X) \simeq \Dmod(K \times X) \xrightarrow{a_*} 
   \Dmod(X),
\]
exactly as in the function-theoretic case of Section \ref{ex:actspace}. Given a 
map $f: X \rightarrow X'$ of $K$-varieties, the functors of $*$-pushforward and 
$!$-pullback along $f$ carry canonical data of $K$-equivariance. 

\subsubsection{} In particular, if we take $X = \on{pt}$, we obtain the trivial 
representation of $K$ on 
\[
    \on{Vect} \simeq \Dmod(\on{pt}).
\]
Concretely, if we write $\pi: K \rightarrow \on{pt}$ for the projection, the 
action map takes the form 
\[
    \Dmod(K) \otimes \Vect \rightarrow \Vect, \quad \quad \sM \otimes V \mapsto 
    \pi_*(\sM) \underset k \otimes  V.
\]

For any representation $\sC$ of $K$, we may accordingly form its invariants and 
coinvariants 
\[
 \sC^K := \Hom_{\Dmod(K)\mod}(\on{Vect}, \sC) \quad \text{and} \quad \sC_K := 
 \on{Vect} \underset {\Dmod(K)} \otimes \sC.
\]
These are canonically identified, as in the case of usual representations of 
finite groups, as follows. 

By the smoothness of $K$, we have a $\Dmod(K)$ equivariant adjunction 
\[
     \pi^![-2 \dim K]: \Dmod(K) \rightleftarrows \Vect: \pi_*,
\]
which upon applying $\Hom_{\Dmod(K)\mod}( -, \sC)$ induces an adjunction 
\[
    \on{Oblv}: \sC^K \rightleftarrows \sC: \on{Av}_*.
\]
As the notation suggests, one thinks of $\on{Oblv}$ as a forgetful functor and 
$\on{Av}_*$ as an averaging functor. For formal reasons, the averaging functor 
factors through a map out of the coinvariants
\[
 \on{Av}_*: \sC_K \rightarrow \sC^K.
\]
By a very useful theorem of Beraldo and Gaitsgory \cite{Bearaldo}, \cite{whitlocglob}, this is an equivalence.

Let us give some basic examples of (co)invariants. 

\begin{ex}Given $\Dmod(K)$ representations $\sC$ and $\sD$, one has a canonical 
equivalence 
	\[
       \Hom_{\Dmod(K)\mod}(\sC, \sD) \simeq \udHom(\sC, \sD)^K,
	\]
	where the invariants are taken with respect to the diagonal action of 
	$\Dmod(K)$ on the space of maps. 
\end{ex}

\begin{ex} \label{e:dmodquotstack} For a $K$-variety $X$, the $K$-equivariant 
objects in $\Dmod(X)$ 
canonically identify with the category of D-modules on the quotient stack 
$X/K$, i.e., 
	\[
	     \Dmod(X)^K \simeq \Dmod(X/K).
	\]
Plainly, D-modules satisfy smooth descent and the diagram \eqref{e:hombar} 
simply is the `$\check{\on{C}}$ech complex' for calculating D-modules on $X/K$ via the cover 
$X 
\rightarrow X/K$.  
\end{ex}

\begin{ex} Under the natural action of $\Dmod(K)$ on its category of Lie 
algebra representations $\fk\mod$, the invariants identify with the 
algebraic representations of the group itself, i.e., 
\[
      \fk\mod^K \simeq K\mod. 
\]
More generally, given a homomorphism $K \rightarrow K'$, the category 
\[
 \fk'\mod^K
\]
canonically identifies with the category of Harish-Chandra modules for the pair $(\fk', K)$. 
\end{ex}

\begin{ex} Here is another class of examples, which in particular gives an alternative perspective on the action on $\fk\mod$. Suppose $A$ is a $k$-algebra equipped with an action of $K$ by automorphisms. In this case, one obtains an action of $\QCoh(K)$ on $A\mod$.

Upon differentiation at the identity, the action of $K$ on $A$ yields  an action $\alpha$ of $\fk$ on $A$ by derivations. Suppose one is further given a trivialization of this action, i.e., a homomorphism of Lie algebras $$\tau: \fk \rightarrow A,$$such that $\alpha$ agrees with the adjoint action via $\tau$. In this case, one obtains an action of $\Dmod(K)$ on $A\mod$. \end{ex}

\begin{re}\label{r:invariantstwice}Let us record here a basic point of departure from the case of usual 
representations. Namely, the operation of taking invariants is not 
idempotent, i.e., the tautological map 
\[
   \on{Oblv}_{\sC^K}: (\sC^K)^K \rightarrow \sC^K
\]
is typically not an isomorphism. Instead, one has a canonical equivalence 
\[
   (\sC^K)^K \simeq (\sC^K) \otimes \Dmod(\on{pt/K}).
\]

This happens because, unless $K$ is 
contractible, the forgetful map $\on{Oblv}_{\sD}: \sD^K \rightarrow \sD$ for any (nonzero) representation $\sD$
will not be 
fully faithful, i.e., Exts in $\sD^K$ and $\sD$ do not agree. Indeed, in the universal case of $\sD = \Dmod(K)$, one has 
\[
\Dmod(K)^K \simeq \Dmod(K / K) \simeq \Vect,
\]
with generator the constant sheaf. In particular its self Exts in $\Dmod(K)$ are $\pi_*(\underline{k})$, i.e., the cohomology of $K$, while its self Exts in $\Vect$ are trivial.

In general, the 
composition 
\[
    \sD^K \xrightarrow{\on{Oblv}_{\sD}} \sD \xrightarrow{\on{Av}_*} \sD^K 
    \xrightarrow{\on{Oblv}_{\sD}} \sD
\]
canonically identifies with $\on{Oblv}_{\sD} \otimes \pi_*(\underline{k})$, i.e., 
tensors the underlying object of a $K$-equivariant object with the cohomology 
of $K$. The reader may wish to compare this, via the functions-sheaves correspondence, 
to the factor of the cardinality of a finite group that arises in the 
discussion of averaging for usual representations. 

Note, however, that if $K$ is unipotent then $\on{Oblv}$ is fully 
faithful, as $K$ is contractible and hence $\pi_*(\underline{k}) \simeq k$. \end{re}

\subsubsection{} Having discussed the functor of invariants, let us indicate how one obtains Theorems \ref{t:loc+equiv} and \ref{t:lochchmod+eq}, beginning with the latter.

Given a map of algebraic groups $H \rightarrow K$, one has 
induction and restriction functors 
\[
  \on{ind}_H^K: \Dmod(H)\mod \rightarrow \Dmod(K)\mod, 
\] 
\[
  \on{res}_K^H: \Dmod(K)\mod \rightarrow \Dmod(H)\mod.
\]
Plainly, $\on{res}_K^H$ restricts the action of $\Dmod(K)$ on a dg-category 
$\sC$ along the monoidal functor 
\[
         \Dmod(H) \rightarrow \Dmod(K),
\]
and $\on{ind}_H^K$ takes a dg-category $\sD$ with an action of $\Dmod(H)$ and 
tensors it up, i.e., 
\[
      \sD \mapsto \Dmod(K) \underset{ \Dmod(H) } \otimes \sD. 
\]
As with Frobenius reciprocity for finite groups, induction and restriction are 
canonically both left and right adjoint to one another. 

\label{ss:invariantscosets} In particular, one obtains a functor of $H$-invariants on 
$\Dmod(K)$ representations, i.e., 
\begin{align*}
     \big(\on{Res}_K^H(-)\big)^H &\simeq   \Hom_{\Dmod(H)\mod}( \on{Vect}, 
     \on{res}_K^H-). 
\intertext{By adjunction we may rewrite this as }
&\simeq  \Hom_{\Dmod(K)\mod}( \on{ind}_H^K \on{Vect}, 
-) \\ & \simeq \Hom_{\Dmod(K)\mod}( \Dmod(K) \underset{\Dmod(H)}{\otimes} \Vect 
, 
-) 
\intertext{By the identification of $H$ invariants and coinvariants and 
Example \ref{e:dmodquotstack}, we may recognize this as}
& \simeq \Hom_{\Dmod(K)\mod}( \Dmod(K/H)  
, 
-).
\end{align*}

The upshot is that the $H$-invariants of any $\Dmod(K)$ representation carry, by 
precomposition, a 
canonical action of 
\[
\Hom_{\Dmod(K)\mod}( \Dmod(K/H), \Dmod(K/H)) \simeq \Dmod(H \backslash K / H),
\]
i.e., the Hecke algebra action of Beilinson--Drinfeld we encountered at the end 
of Section 
\ref{s:bbloc}.

\subsubsection{} Let us finally sketch, for a smooth $K$-variety $X$, how one obtains the $\Dmod(K)$ equivariant adjunction 
\[
    \on{Loc}: \fk\mod \rightleftarrows \Dmod(X): \Gamma(X, -).
\]
To begin with, induction from $\OO_K$-modules to $\sD_K$-modules yields a monoidal functor
\[
      \QCoh(K) \rightarrow \Dmod(K).
\]
Similarly, one has a canonical action of $\QCoh(K)$ on $\Dmod(X)$ such that the forgetful functor $\Dmod(X) \rightarrow \QCoh(X)$ carries a canonical datum of $\QCoh(K)$-equivariance. In particular, the functor of global sections, i.e., the composition 
\[
    \Dmod(X) \rightarrow \QCoh(X) \xrightarrow{\pi_*} \QCoh(\pt) \simeq \Vect,
\]
carries a canonical datum of $\QCoh(K)$-equivariance. 

However, for any $\Dmod(K)$ representation $\sC$, one has by adjunction 
\begin{align*}
\Hom_{\QCoh(K)\mod}( \on{Res} \sC, \Vect) &\simeq \Hom_{\Dmod(K)\mod} (\sC, \Hom_{\QCoh(K)\mod}(\Dmod(K), \Vect)) \\ &\simeq \Hom_{\Dmod(K)\mod}(\sC, \fk\mod), 
\end{align*}
where for the last equivalence one uses the identification of $\fk\mod$ and weakly $K$-equivariant D-modules on $K$, cf. Section \ref{ss:weakinv}. This yields the desired $\Dmod(K)$ equivariance of the global sections functor. For localization, one can use the nice fact, due to Gaitsgory \cite{whitlocglob}, that any adjoint of a $\Dmod(K)\mod$-equivariant functor inherits a compatible datum of equivariance.

\subsection{Complements I - character sheaves}

\label{ss:charsheaves}
\subsubsection{} In this section, we began by sketching some preliminary 
generalities on dg-categories. We then discussed some of the basic definitions and 
constructions in categorical representation theory, sufficient to recover some 
results previewed in Section \ref{s:bbloc}.

In the remainder of this section, we would like to survey a few more basic 
results and constructions. However, the reader may wish to skip directly to 
Section \ref{s:catloop}.

\subsubsection{} Character sheaves are (certain) adjoint equivariant D-modules 
on groups. Lusztig initiated their study, particularly in the reductive case, in view of their close ties to irreducible characters of 
finite groups of Lie type \cite{luschar1, luschar23, luschar4, luschar5}. This is a particularly rich and nontrivial example of the 
functions-sheaves correspondence. 

Of course, the characters of a finite group are not simply adjoint equivariant 
functions on the group, but rather arise as traces from representations. One 
can ask whether character sheaves have such an interpretation. Indeed they 
do, and we turn to this next. 

\subsubsection{Duality} \label{ss:dualizable}Recall that to define the character of a 
representation of a finite group, one needs the representation to be finite 
dimensional so that one may safely take traces. The analogue of being `small enough' to take traces in the categorical setting is as follows.

Using the monoidal structure on $\DGCat$ one can speak 
of dualizable dg-categories. Under the analogy between $\DGCat$ and $\Vect$, 
these correspond to (bounded complexes of) finite dimensional vector spaces, 
and indeed a dg-category $\sC$ is dualizable if and only if the natural map 
\[
     \sC \rightarrow \udHom( \udHom(\sC, \Vect), \Vect)
\] 
is an equivalence. 

For a dualizable dg-category, one has a canonical 
equivalence $\sC^\vee \simeq \udHom(\sC, 
\Vect)$, i.e., its dual in the monoidal sense is simply its continuous 
dg-modules. For any $\sS$ the natural map 
\[
    \sS \otimes \sC^\vee \rightarrow \udHom(\sC, \sS)
\]
is an equivalence.\footnote{Replacing $\sC^\vee$ with $\Hom(\sC, \Vect)$, this 
is yet another characterization of dualizable dg-categories.} More generally, 
there are equivalences
\[
\Hom( \sS \otimes \sC, \sS') \simeq \Hom( \sS, \sC^\vee \otimes \sS')
\]
for any dg-categories $\sS$ and $\sS'$.

  Let $\sA$ be a monoidal dg-category. By functoriality, if $\sC$ is a left $\sA$-module, then $\sC^\vee$ is a right $\sA$-module. In particular, for a 
  dualizable $\Dmod(K)$ representation $\sC$, its dual acquires a canonical action of 
  $\Dmod(K)$, i.e., one may form the contragredient representation.

Before we discuss how to define the character sheaf of a dualizable $\Dmod(K)$ representation, we would like to discuss some basic examples of dualizable categories and modules. First of 
all, the following general fact ensures that most dg-categories one encounters in practice are 
dualizable. 

\begin{ex}If a category $\sC$ is compactly generated, then it is dualizable. 
Moreover, its dual is compactly generated, and one has a canonical equivalence 
of non-cocomplete dg-categories between the compact objects of $\sC$ and the 
opposite of the compact objects of $\sC^\vee$. 
\end{ex}

Here are some useful special cases of the previous example. 
\begin{ex} Given a $k$-algebra $A$, its category of left-modules $A\mod$ is 
dualizable with dual the category of right modules $\on{mod-}A$. Explicitly, 
the 
evaluation map 
	\[
	   \rmod A \otimes A\mod \rightarrow \Vect
	\]
	sends a right module $M$ and left module $N$ to their tensor product $M 
	\underset A \otimes N$. 
\end{ex}

\begin{ex} \label{ex:dmodselfdual}Given a variety $X$, its category of D-modules is made canonically 
self-dual by the operation of Verdier duality on compact D-modules. Explicitly, the 
resulting evaluation map is given by the composition 
	\[
	 \Dmod(X) \otimes \Dmod(X) \simeq \Dmod(X \times X) \xrightarrow{\Delta^!} 
	 \Dmod(X) \xrightarrow{\pi_*} \Vect,
	\]
	where, exactly as in Section \ref{s:funselfdual}, we denote by $\Delta: X 
	\rightarrow X \times X$ the diagonal map and $\pi: X \rightarrow \on{pt}$ 
	the projection to the 
	point. 
	
	If $X$ is a $K$-variety, this pairing is canonically 
	$\Dmod(K)$ equivariant, i.e., identifies $\Dmod(X)$ with its contragredient 
	representation. Similar statements hold, mutatis mutandis, for $\QCoh(X)$, viewed as a $\QCoh(K)$ representation. 
\end{ex}

\subsubsection{Character sheaves} We are now ready to explain how to take the character sheaf of a representation. An action of $\Dmod(K)$ on a category $\sC$ 
is the same as a monoidal functor 
\[
\Dmod(K) \xrightarrow{\rho} \udHom(\sC, \sC).
\]
If $\sC$ is dualizable, then its endomorphisms $\udHom(\sC, \sC)$ are again dualizable by the identification $\udHom(\sC, \sC) \simeq 
\sC^\vee \otimes \sC$. We may therefore dualize the action 
map to obtain a $\Dmod(K \times K)$-equivariant functor of matrix coefficients 
\[
   \udHom(\sC, \sC) \simeq \udHom(\sC, \sC)^\vee \xrightarrow{\rho^\vee} 
   \Dmod(K)^\vee \simeq \Dmod(K).
\]

In particular, by passing to diagonal $K$-invariants, we obtain a functor 
\begin{equation} \label{e:mcequivfunc}
   \Hom_{\Dmod(K)}(\sC, \sC) \rightarrow \Dmod(K / K ).
\end{equation}
 Applying this to the identity endomorphism $\id_{\sC}$, we obtain an adjoint 
 equivariant D-module $\chi_\sC$ on $K$, 
 namely the {\em character sheaf} of $\sC$.

\subsubsection{} Here is a more hands-on description of $\chi_\sC$. Given a 
dualizable category $\sC$ and an endofunctor $f: \sC \rightarrow \sC$, one can 
take its trace, which is a complex of vector spaces. Exactly as in linear 
algebra, this is computed via the composition 
\[
\Vect \xrightarrow{\id_{\sC}} \sC \otimes \sC^\vee \xrightarrow{f \otimes 
\on{id}_{\sC^\vee}} \sC \otimes \sC^\vee \xrightarrow{\on{ev}} \Vect.
\]
That is, the composite endofunctor of $\Vect$ is given by tensoring by a vector space, which is by 
definition $\on{tr}(f, \sC)$. We remark that, when $f$ is simply the identity 
endomorphism of $\sC$, the resulting complex is known as the Hochschild 
homology of $\sC$.

Having understood how to take traces in the present categorical context, we can now describe the character sheaf fairly explicitly. Namely, its stalks store the traces of the 
corresponding automorphisms of $\sC$, exactly as in usual representation 
theory. That is, for any closed point $g: \pt \rightarrow K$, one has a 
canonical isomorphism 
\[
     g^!(\chi_\sC) \simeq \on{tr}( \delta_g \star -, \sC). 
\]

\subsubsection{}Let us mention a few features and examples. 

\begin{ex} A $\Dmod(K)$ equivariant functor between dualizable modules 
$\sC \rightarrow \sD$ 
need not induce a map between their character sheaves. However, if it admits 
a continuous right adjoint, we do obtain a suitably functorial map $\chi_\sC 
\rightarrow 
\chi_{\sD}$. 
	
	Let us remark that, in practical situations, $\sC$ is compactly generated, 
	and the condition of being right adjointable is equivalent to sending 
	compact objects to compact objects, i.e., nothing `small' becoming `big'. 
	\end{ex}

\begin{ex} The induction and restriction functors on categorical 
representations preserve the property of dualizability.  On characters, they 
reproduce the usual functors of induction and restriction for character 
sheaves. 
\end{ex}

\begin{ex} If $X$ is a $K$-variety, its character sheaf may be described as 
follows. Consider the 
closed subscheme $\on{Stab}$ 
of $K \times X$, whose fibre over a point of $X$ is its stabilizer. That is, 
$\on{Stab}$ is the 
equalizer of the action and 
projection maps to $X$. Then under the tautological $G$-equivariant map 
\[
\on{Stab} \rightarrow G,
\]
the pushforward of the dualizing sheaf is $\chi_{\Dmod(X)}$. In particular, by 
proper base change, the 
trace of any element of $K$ on $\Dmod(X)$ is the Borel--Moore homology of its 
fixed points. The reader may wish to compare this with the analogous assertion for finite groups and 
sets, where the trace of an element of the group recovers the number of fixed points. 
\end{ex}

\begin{ex} If $K=G$ is a reductive group, and $\sC = \Dmod(\sB)$, then by the 
preceding example its 
character sheaf is the Grothendieck--Springer sheaf on $G/G$. More generally, 
the functor of matrix coefficients for equivariant endomorphisms, cf. Equation 
\eqref{e:mcequivfunc}, takes the form
\[  \Dmod(B \backslash G / B) \rightarrow \Dmod(G/G).\]
This recovers the Harish-Chandra, or horocycle, transform, i.e., push-pull along 
the correspondence $$G/G \leftarrow G/B \rightarrow B \backslash G / B,$$ where 
the left and middle terms are quotients with respect to the adjoint action. The 
latter was the point of departure for Lusztig's study of character sheaves. 
\end{ex}

Here are two related remarks. 

\begin{re}In the present discussion, we treated character sheaves as a 
direct construction, which takes as input dualizable representations and 
outputs adjoint equivariant D-modules. 

An alternative perspective is to interpret $\Dmod(K/K)$ as the horizontal trace 
of $\Dmod(K)$. Indeed, for any monoidal dg-category $\sA$, the latter is 
defined as $$\sA \underset{ \sA \otimes \sA^{\on{op}}}\otimes \sA$$or, 
equivalently, as 
the 2-categorical trace of $\sA\mod$, viewed as a dualizable $\DGCat$-linear 
category. If we were working one categorical level lower, i.e. with an algebra $A$ instead of a monoidal dg-category, the analogue of this discussion is the presentation of the Hochschild homology of $A\mod$ as
\[
  A \underset{A \otimes A^{\on{op}}} \otimes A. 
\]

In the present case, using the correspondence between categorical coinvariants 
and D-modules on quotient stacks, we 
have a tautological identification 
\[
   \Dmod(K) \underset{ \Dmod(K \times K)} \otimes \Dmod(K) \simeq \Dmod(K 
   \overset{K \times K} \times K) \simeq \Dmod(K / K).
\]
General formalism assigns to a dualizable $\sA$-module an 
object of its horizontal trace, which recovers the preceding discussion of 
character sheaves as a 
special case; see \cite{gkrvtoyshtuka} and \cite{bznadlerchartheory} for further discussion.
\end{re}

\begin{re}There are parallel theories to those discussed in this section, where 
one works with analytic or \'etale constructible sheaves on $K$ rather than 
D-modules, cf. Section \ref{ss:variants} below. When one takes \'etale sheaves on $K$ defined over a finite field 
$\Fq$, a version of the functions-sheaves correspondence 
identifies the category of modules for $K(\Fq)$ with the trace of Frobenius of 
$\on{Shv}(K)\mod$. This sends categorical representations to usual 
representations and, compatibly, character sheaves to usual 
characters; see \cite{Gaitsgoryinventshtukas}. 
\end{re}

\subsection{Complements II: generation}
\subsubsection{}Our final topic is several nice 
theorems about generation. Recall a basic and useful feature of life is Morita 
theory -- for a given ring $A$, while its category of modules $A\mod$ is 
compactly generated by $A$ itself, it admits many other compact generators $c$. For any choice of one, taking its endomorphisms yields an equivalence 
\[
  \rmod\End(c)  \simeq   A\mod, \quad \quad M \mapsto M \underset{\End(c)} 
  \otimes c.
\] 
Said somewhat differently, a given dg-category $\sC$ may have many compact 
generators and correspondingly many presentations as a category of modules. We 
would like to review here a couple useful statements of this form, which 
necessarily take place one categorical level higher, for 
$\Dmod(K)\mod$.

\subsubsection{} \label{ss:hchbimod}
 First of all, a nontrivial theorem of Beraldo \cite{beraldo}, building on work of Gaitsgory \cite{1affine}, states 
that $\Dmod(K)$ is generated by $\fk\mod$. Note that the endomorphisms of the 
latter are 
\begin{align*}
   \Hom_{\Dmod(K)\mod}(\fk\mod, \fk\mod) & \simeq \udHom(\fk\mod, \fk\mod)^K \\ 
   & \simeq (\fk\mod^\vee \otimes \hspace{.7mm} \fk\mod)^K \\ & \simeq (\fk\mod \otimes \hspace{.7mm}
   \fk\mod)^K \\ & \simeq (\fk \oplus \fk\mod)^K.
\end{align*}
Stated in words, the equivariant endomorphisms of $\fk\mod$ are the category of Harish-Chandra 
bimodules, i.e., $\fk$-bimodules for which 
the diagonal action is integrated to the group $K$. 

\begin{re} We emphasize that this identity is an attractive and useful feature 
of the 
formalism. Namely, Harish-Chandra 
bimodules have been studied and used in representation theory for many decades 
and have been understood to have many favorable properties; see \cite{bege1980} and references therein. The above 
articulates this in a precise sense -- they are exactly the endofunctors 
of $\fk\mod$ which commute with the categorical $K$-action and vice versa.
\end{re}

The theorem of Beraldo therefore affords an equivalence 
\[
    \Dmod(K)\mod \simeq (\fk \oplus \fk\mod^K)\mod.
\]
In particular, for any $\Dmod(K)$ representation $\sC$ we have a 
canonical equivalence 
\[
    \fk\mod \underset{ \fk\oplus \fk\mod^K} \otimes \sC^{K, w} \simeq \sC.
\]
That is, recalling the equivalence $\Dmod(K)^{K,w} \simeq \fk\mod$, a 
$\Dmod(K)$ representation can be reconstructed from its weak invariants. 

\subsubsection{}We would also like to mention, for a reductive group $G$, a 
nice analogue of highest weight theory in the present setting due to 
Ben-Zvi--Gunningham--Orem \cite{BZO}. Namely, they showed that, for any maximal unipotent subgroup $N$, an arbitrary $\Dmod(G)$ representation $\sC$ can be reconstructed from its 
$N$-invariants, or even the weak Cartan 
invariants of the latter 
\[
   \sC^N \quad \text{and} \quad (\sC^N)^{T,w}.
\]
Passing to their endomorphisms, which roughly are Hecke categories without 
prescribed monodromy along $T$-orbits, we obtain equivalences 
\[
  \Dmod(G)\mod \simeq \Dmod(N \backslash G / N)\mod \simeq \Dmod(N \backslash G 
  / N)^{T \times T, w}\mod
\]
and, in particular, identities
\[
   \Dmod(G/N) \underset{ \Dmod(N \backslash G / N)} \otimes  \sC^N \simeq  
   \Dmod(G/N)^{T,w} \underset{ \Dmod(N \backslash G / N)^{T \times T,w}}  
   \otimes 
   (\sC^{N})^{T,w} \simeq \sC. 
\]

\begin{re}The argument of \cite{BZO} yields the following useful variation. Given any 
subgroup $H$ of a 
group $K$ such that the quotient of $K$ by the normalizer of $H$ is proper, 
the tautological map 
	\[
	     \Dmod(K/H) \underset{\Dmod(H \backslash K / H)} \otimes \sC^K 
	     \rightarrow \sC,
	\]
which need not be essentially surjective, is always a fully faithful 
embedding. This may be applied to check equivalences 
of categories on `highest weight vectors,' which are often much easier to 
directly analyze, cf. \cite{ahc}, \cite{yangmoyprasad} for further discussion and examples. 	
\end{re}

\section{Categorical representations of loop groups}
\label{s:catloop}

	\subsection{Overview}
	
	\subsubsection{}  So far, we have discussed Beilinson--Bernstein 
	localization and categorical representation theory for reductive groups. 
	
	In the remainder of this survey, we would like to turn to their affine 
	analogues, i.e., 
	we pass from reductive groups to loop groups. The story here is 
	significantly richer and is the subject of the local geometric Langlands 
	correspondence.

	In this section, we will give some introductory comments on loop groups and 
	their categorical representations to prepare for a discussion of the 
	Langlands correspondence in Section \ref{s:affbbloc}.

	\subsection{Loop groups}

	\subsubsection{} As before, we continue to denote by $G$ a connected 
	reductive group over an algebraically closed field $k$ of characteristic 
	zero. 
	
	The main actor in this section is not $G$, but rather its loop group $LG$. At first 
	pass, just 
	as $G$ consists of certain matrices with entries in $k$, i.e., pick a 
	faithful representation,  $LG$ consists of the same form of matrices but 
	with entries in the field of Laurent series $k(\!(t)\!)$. Somewhat more 
	carefully, since we need to make sense of $LG$ as an algebro-geometric 
	object, it is a certain group ind-scheme whose closed points are as we just 
	described. 
	
	\begin{ex}\label{ex:ga} Here is an explicit construction in the case when 
	$G$ is the 
	(non-reductive) additive group $\mathbb{G}_a$, which contains most of the 
	representative ideas. In this case, we have 
	$\mathbb{G}_a(k) \simeq k$, i.e., it is a one dimensional affine space. As above, we 
	would like $L\mathbb{G}_a(k) \simeq k(\!(t)\!)$. The basic idea is that 
	$k(\!(t)\!)$ is again an infinite dimensional affine space, and, in fact, by 
	bounding the order of the pole of the Laurent series, an ascending union of 
	pro-finite dimensional affine spaces. Explicitly as a vector space
	\[
	  k(\!(t)\!) = \varinjlim_n \varprojlim_m t^{-n} \cdot k[[t]] / t^{-n+m} 
	  \cdot k[[t]],
	\]
	and we use the same formula to describe its structure as an ind-scheme. 
	That is, fixing $n$, we obtain a pro-finite dimensional affine space, which 
	explicitly is the spectrum of a polynomial ring in countably many variables, i.e., the functions corresponding to the Taylor coefficients of a series. 
	Allowing $n$ to vary, we obtain an ascending union of infinite dimensional 
	affine spaces along closed embeddings, i.e., an ind-scheme. 
	\end{ex}
	
	\begin{re} In general, $LG$ admits a succinct description via its functor 
	of points, which is defined as 
	follows. For any variety $X$, consider the functor $LX$ from affine 
	varieties to sets whose 
	$R$-points, for any $k$-algebra $R$, are
		\[
		     LX(R) = X(R(\!(t)\!) ).
		\]
	While this is not well-behaved for a general variety, if $X$ is affine, this 
	is representable by an ind-scheme of ind-infinite type. Taking $X = G$ 
	yields the loop group $LG$. 
	\end{re}

	Having given a formal definition, let us try to give a slightly more hands-on description in the spirit of Example \ref{ex:ga}. Recall the analogy 
	between $G(\Fq)$ and $G(k)$ and, more 
	generally, finite $G(\Fq)$-sets and finite dimensional $G(k)$-varieties as 
	discussed in the previous sections. A basic guiding principle for the 
	present situation is an 
	analogy between the $p$-adic group 
	$G(\Fq(\!(t)\!))$ and the loop group $LG$. 
	
\begin{re} While some authors use the term $p$-adic group to mean the points of $G$ in $\mathbb{Q}_p$ or a finite extension thereof, here by a mild abuse of notation we call the points in any non-Archimedean local field, such as $\Fq(\!(t)\!)$, a $p$-adic group. \end{re}

Recall that $p$-adic groups, with their natural topology, are ind-pro-finite 
sets. Explicitly, fix any compact open subgroup $K$, i.e., a subgroup such that $K \cap G(\Fq[[t]])$ is finite index in both $K$ and $G(\Fq[[t]])$. Then $K$ is 
naturally an inverse limit of finite groups, and the entire group decomposes as 
a disjoint union of its $K$ cosets. Concretely, for $G(\Fq[[t]])$, the 
relevant inverse limit comes from the topology on power series, i.e., 
\begin{equation}
\label{e:invlim}
      G(\Fq[[t]]) = \varprojlim_n G(\Fq[t]/(t^n)).
\end{equation}

In a similar way, one can make sense of compact open subgroups of $LG$. Namely, 
one has the arc group $L^+G$, which is a group scheme of infinite type with 
$k$-points $G(k[[t]])$. This is defined as an inverse limit of finite dimensional 
algebraic groups similarly to \eqref{e:invlim}. A sub-group scheme $K$ 
of $LG$ is said to be compact open if $K \cap L^+G$ is of finite 
codimension in both $K$ and $L^+G$. Then $K$ is naturally an inverse limit of algebraic 
groups, and the coset space $$LG/K$$is an ind-scheme of ind-finite type. We 
emphasize that, under the analogy between finite sets 
and finite dimensional varieties, this corresponds to being a countable discrete union 
of 
cosets  in the case of $p$-adic groups.

 \subsection{D-modules on loop groups and categorical actions}

	\subsubsection{} We would next like to discuss the convolution algebra of 
	D-modules on $LG$. Since this is an infinite-dimensional object, some care 
	is required. 
	
	Its definition may be arrived at as follows. Recall that we thought of 
	modules for $\Dmod(G)$ as akin to representations of $G(\Fq)$. In 
	particular, the convolution algebra of D-modules itself was analogous to 
	the group algebra $\Fun(G(\Fq))$. 
	
     We would similarly like modules for $\Dmod(LG)$ to 
	behave like representations of the $p$-adic group $G(\Fqt)$. As the latter 
	is a topological group, one correspondingly works with smooth 
	representations.\footnote{Explicitly,  a representation on a complex vector space $V$ is smooth if the action map $G(\Fqt) 
	\times V \rightarrow V$ is continuous with respect to the discrete 
	topology on $V$, i.e., the stabilizer of each vector is open. For the reader 
	unfamiliar with $p$-adic representation theory, the point here is that this condition mirrors a similar one for Lie 
	group representations, where one asks for the action to be suitably smooth in the sense of manifolds. In the $p$-adic setting, however, the 
	topologies on the group and on the vector space are, 
	informally speaking, `coprime.' Practically, smooth representations are what show up in 
	nature when studying automorphic forms.} For this reason, modules for $G(\Fqt)$ are not the same as modules for its group algebra, but instead for 
	its Hecke algebra 
	$$\on{Dist}(G(\Fqt)) $$
	of compactly supported distributions. Concretely, this can be written as 
	the 
	colimit 
	\begin{equation} \label{e:padicgroups}
	\on{Dist}(G(\Fqt))\simeq \varinjlim_K  \on{Dist}( K \backslash G(\Fqt) / K),
	\end{equation}
	where $K$ runs over compact open subgroups and each $\on{Dist}( K \backslash 
	G(\Fqt) / K)$ is simply finite linear combinations of delta distributions on 
	double cosets, i.e., distributions that integrate a function over a double 
	coset $KgK$ with respect to some Haar measure.

	With this in mind, the category of D-modules on the loop group is defined 
	as 
	\[
	\Dmod(LG) \simeq \varinjlim \Dmod(K \backslash LG / K),
	\]
	where $K$ runs over compact open subgroups and the colimit is taken in 
	$\DGCat$ along $*$-pullback functors. It remains then to describe, for 
	a fixed $K$, the category 
	$\Dmod(K\backslash LG/K)$. However, recalling that  $LG/K$ is an ind-scheme 
	of ind-finite type, we may moreover write it as an ascending union 
	\[
	   LG/K = \bigcup_i X_i,
	\]
	where each $X_i$ is a closed subscheme of finite type that is moreover $K$ 
	invariant. As $K$ acts on each $X_i$ via a finite dimensional quotient $Q_i$ 
	with prounipotent kernel, we have 
	\[
	   \Dmod(X_i)^K \simeq \Dmod(X_i)^{Q_i} \simeq \Dmod(Q_i \backslash X_i),
	\]
	i.e., we are dealing with D-modules on a usual finite dimensional stack. We then 
	set 
	\[
	  \Dmod(K \backslash LG/K) \simeq \varinjlim \Dmod( Q_i \backslash X_i). 
	\]

	The point here is that, roughly speaking, objects in 
	$\Dmod(LG)$ resemble compactly supported distributions on $G(\Fqt)$. More 
	precisely, by Example \ref{e:colimitcats}, $\Dmod(LG)$ is compactly 
	generated, with compact objects given by $*$-extensions of bounded 
	complexes of coherent D-modules from subvarieties of $LG/K$, for some $K$.

	While the exact details of this construction are not so important for us, here are a couple orienting remarks. 
	
	\begin{re}First, one can replace all possible $K$ in the above colimit \eqref{e:padicgroups} with 
	any cofinal subset. A concrete such choice would be the congruence 
	subgroups $K_n$, for $n \geqslant 1$, of $L^+G$, which we will simply 
	describe at the level of 
	$k$-points via the exact sequence 
		\[
		      1 \rightarrow K_n \rightarrow L^+G = G(k[[t]]) \rightarrow 
		      G(k[[t]]/(t^n)) \rightarrow 1. 
		\]

	Noting that $K_1$ is the prounipotent radical of $L^+G$, it follows that 
	any sufficiently small compact open subgroup $K$ is 
	prounipotent. In particular, for two such prounipotent subgroups $K' \subset K$, the 
	pullback 
	\[
	    \Dmod(K \backslash LG / K) \rightarrow \Dmod(K' \backslash LG / K')
	\]
	is a fully faithful embedding, cf. Remark \ref{r:invariantstwice}. 
	\end{re}

	\begin{re}While the above construction may seem somewhat ad hoc, we should 
	mention that it is a special case of a similar definition of D-modules for 
	a larger 
	class of infinite dimensional varieties known as placid ind-schemes \cite{rdm}.
	
	In fact, a basic subtlety in this infinite dimensional setting is that there are two dual versions of the category of D-modules. One comes with functoriality for $*$-pushforwards and the other with functoriality for $!$-pullbacks. The presentation above is for the category with $*$-pushforwards. This is the analogue of the distinction between functions and distributions on $G(\Fqt)$. Relatedly, the dual categories for $LG$ may be identified after trivializing the dimension torsor, i.e., fixing a compact open subgroup. 
	\end{re}

	\subsubsection{} As in the finite dimensional case, $\Dmod(LG)$ is a 
	monoidal dg-category under convolution. In particular, as before, we can 
	define its category of modules 
	\[
      \Dmod(LG)\mod.
	\]

	Here are some basic examples of representations. 
	
	\begin{ex} If $X$ is a placid ind-scheme acted on by $LG$, its category of 
	D-modules $\Dmod(X)$ carries a canonical action of $\Dmod(LG)$ by 
	convolution. 		
	\end{ex}

	\begin{ex}\label{ex:kmmods} Given a central extension of the loop group
		\[
		1 \rightarrow \mathbb{G}_m \rightarrow \widetilde{LG} \rightarrow LG 
		\rightarrow 1,
		\]
		one obtains upon passing to Lie algebras a central extension of the loop algebra, i.e., an affine Lie 
		algebra
		\[
		 0 \rightarrow k \cdot 1 \rightarrow \gk \rightarrow L \fg \rightarrow 
		 0. 
		\]
		There is a canonical action of $\Dmod(LG)$ on the derived category of 
		smooth modules
		$
		    \gk\mod,
		$
		induced by the fact that the adjoint action of $\widetilde{LG}$ on $\widehat{\fg}$ factors through $LG$. By definition, an object of the corresponding abelian category is a 
		representation of $\gk$ on which 1 acts by the identity and where each vector 
		is annihilated by a compact open subalgebra. The full dg-category, 
		following Frenkel--Gaitsgory \cite{fglocaff}, is a slight 
		modification of the naive unbounded derived category, such that objects 
		like Verma modules or, more generally, inductions of finite dimensional 
		representations of compact open subalgebras are compact generators.

		While 
		the details of this will not be important for us, let us simply remark 
		that the existence and utility of this modification stem from the fact that the abelian category 
		is of infinite cohomological dimension. Relatedly, the renormalized derived 
		category carries a $t$-structure and canonically identifies on bounded 
		below parts with 
		the usual derived category. That is, the renormalization leads only to 
		differences in  
		cohomological degree 
		$-\infty$. Similarly, the category 
		$\Dmod(LG)$ carries a $t$-structure, which is canonical up to an overall 
		shift,\footnote{This mirrors the uniqueness of the Haar measure on 
		$G(\Fqt)$ up to scaling.} with the property that $\Dmod(LG)$ agrees with the naive unbounded 
		derived 
		category of its heart, modulo issues in cohomological degree 
		$-\infty$.   
	\end{ex}

		\section{Local geometric Langlands and affine Beilinson--Bernstein localization}
	
	\label{s:affbbloc}
	
	\subsection{Overview}
	
	\subsubsection{}  In the previous sections, we discussed localization theory for reductive Lie algebras and explained its relation to the categorical representation theory of reductive groups. In Section \ref{s:catloop}, we met the loop groups and their categorical representations.

	In this final section, we will sketch some of 
	the basic conjectures and results in local 
	geometric Langlands, which may be understood as parametrizing the categorical representations of loop groups.  At the end, we will come full circle and meet an 
	analogue 
	of Beilinson--Bernstein localization for loop groups, as developed by 
	Beilinson--Drinfeld, Frenkel--Gaitsgory, and Raskin--Yang.

	\subsection{A naive formulation of the correspondence}

	\subsubsection{} For the $p$-adic group $G(\Fqt)$, the local Langlands 
	correspondence (conjecturally) parametrizes its irreducible representations 
	on $k$ vector spaces in terms of Galois theoretic data, to first 
	approximation homomorphisms 
	from the absolute 
	Galois group of $\Fqt$ into the Langlands dual $\check{G}(k)$.

	In this subsection we would like to review what is expected and what is 
	known about an analogous parametrization of $\Dmod(LG)\mod$.

	\subsubsection{} Recall that $k(\!(t)\!)$ consists of functions on a formal 
	punctured disk $D^\times$. Following the analogy between Galois 
	representations and local systems, cf. Appendix \ref{s:frmfun}, the counterparts to $\cG$-valued Galois 
	representations in the present setting are flat $\check{G}$ connections, i.e., $\cG$ de Rham local 
	systems 
	on 
	$D^\times$.
	
	 Concretely, these are given by connections modulo gauge 
	equivalence, i.e., 
	\[
	    \LS \simeq \{ d + \cfg(\!(t)\!)dt \} / \cG(k(\!(t)\!)). 
	\]
	The same formula defines $\LS$ further as an object of algebraic geometry, 
	i.e., one takes the quotient prestack of the action of the $L\cG$ on the 
	ind-scheme of connection forms $d + \cfg(\!(t)\!)dt$, which is defined 
	similarly to Example \ref{ex:ga}; see \cite{raskinspectraldecomp} for further discussion.

	\subsubsection{} \label{s:infdiscuss} Let us try to arrive at the statement of the local 
	geometric Langlands 
	correspondence in a somewhat heuristic fashion. To first approximation, we 
	can hope that irreducible $\Dmod(LG)$ representations are in 
	bijection with the points of $\LS$. That is, to each $\cG$-connection 
	$\sigma$ one should be able to canonically attach an irreducible categorical 
	representation $\sC_{\sigma}$ of $LG$. 
	
	Because a general object $\sC$ of $\Dmod(LG)\mod$ should decompose as a direct 
	integral of irreducible representations, it should define a corresponding 
	family of categories over $\LS$, whose fibres over each point measure the 
	corresponding multiplicity of $\sC_{\sigma}$ in the direct integral. 
	
	Since irreducibility is a feature of objects of abelian categories, 
	this discussion is difficult to formulate precisely. However, there is a 
	good notion of a family of categories over a space, namely a 
	quasicoherent sheaf of categories, which we review below. 
 
    Modulo corrections, which we will meet duly in Equation \eqref{e:loclang} below, this gives the 
    statement of the local geometric Langlands correspondence. Namely, the conjecture is 
    that there should be an equivalence of $(\infty,2)$-categories between all 
    categorical representations of the loop group and sheaves of categories on 
    the moduli space of local systems, i.e., 
    \begin{equation}
    \label{e:loclangnaive}       \mathbb{L}: \Dmod(LG)\mod \simeq 
    \on{2-QCoh}(\LS). 
    \end{equation}

	Here is a basic orienting remark. The asserted equivalence of 2-categories is not saying that for a given categorical representation $\sC$, the underlying 
	dg-categories of $\sC$ and $\mathbb{L}(\sC)$ are equivalent. Rather, it is saying that, for any pair of objects $\sC$ and $\sD$, there is an
	equivalence 
	\[
         \Hom_{\Dmod(LG)\mod}(\sC, \sD) \simeq 
         \Hom_{\TwoCoh(\LS)}(\mathbb{L}(\sC), \mathbb{L}(\sD)),
	\]
	and these are compatible with composition as one varies $\sC$ and $\sD$. 
	
	\begin{re}  As a toy model, a Morita equivalence for usual algebras 
    \[
        A\mod \simeq B\mod
    \]
	matches not the underlying vector spaces of 
	modules, but 
	only their homomorphisms. Explicitly, one can consider, for a finite 
	dimensional $k$ vector space $V$, the Morita equivalence 
	\[
	   V \otimes -  :\Vect \simeq  \End(V)\mod.
	\] 
	This exchanges $k$ and $V$, which certainly are not isomorphic, but on 
	endomorphisms we do have
	\[
	     \Hom_{\Vect}(k, k) \simeq k \simeq \Hom_{\End(V)}(V,V).
	\]
		\end{re}

	Before turning to corrections, we would like to first discuss  
	sheaves of categories and  some of the basic compatibilites $\mathbb{L}$ 
	is 
	expected to satisfy. 
	
	\subsubsection{}   To discuss sheaves of categories, let us start with the 
	case of an affine variety $X = \on{Spec }O_X$. In this case, a 
	quasi-coherent 
	sheaf of categories on $X$ is simply an $O_X$-linear, as opposed to merely
	$k$-linear, dg-category. Since $O_X$-linearity is the same as being able to 
	tensor by complexes of $O_X$-modules, we have equivalently that 
	\[
	  \TwoCoh(X) \simeq \QCoh(X)\mod.
	\]
	In particular $\TwoCoh(\on{pt})$ is simply $\DGCat$, and, for any closed 
	point $ \on{pt} \rightarrow X$, one has a pullback 
	\[
	  \Vect \underset{ \QCoh(X)} \otimes  -: \TwoCoh(X) \rightarrow \DGCat, 
	\]
	i.e., one can speak of the fibres of a sheaf of categories, which are usual 
	dg-categories. In this way, one thinks of a sheaf of categories as a 
	family of categories over $X$.

	\begin{re}\label{r:pushpull2coh} More generally, given any map $f: X 
	\rightarrow Y$ pullback 
	defines a monoidal functor $\QCoh(Y) \rightarrow \QCoh(X)$, and, in 
	particular, adjoint induction and restriction functors 
		\[
		    f^*: \TwoCoh(Y) \leftrightarrows \TwoCoh(X): f_*. 
		\]
	\end{re}

	Here is a basic example of a sheaf of categories. 
	
	\begin{ex}Given a map $\tilde{X} \rightarrow X$ of varieties,   
	$\QCoh(\tilde{X})$ is a $\QCoh(X)$-algebra and, in particular, a sheaf of categories over $X$. For any 
	point $x$ of $X$, the fibre of the sheaf of categories is, by Example 
	\ref{e:fibprodvar}, simply 
	the quasicoherent sheaves on the actual fibre, i.e., 
	\[
	  i_x^* \QCoh(\tilde{X}) \simeq \QCoh( \tilde{X} \underset X \times x).
	\]
	\end{ex}

	\subsubsection{} For a general $X$, not necessarily affine, the definition 
	proceeds by gluing, 
	exactly as in the definition of quasicoherent sheaves. While the details are not so important for us, we should mention the 
	following beautiful phenomenon established by Gaitsgory in \cite{1affine}. Namely, for general $X$, one by definition produces sheaves of categories by 
	patching them affine by affine. On the other hand, one also has the global 
	monoidal category of quasicoherent sheaves and hence its category of left 
	modules $\QCoh(X)\mod$. One can therefore ask for the relation between the two. As for usual 
	quasicoherent sheaves, one 
	has a 
	global sections functor 
	\begin{equation} \label{e:globsectshvcat}
	 \Gamma: \TwoCoh(X) \rightarrow \QCoh(X)\mod.
	\end{equation}
    For usual sheaves, the analogous map 
    \[
    \Gamma: \QCoh(X) \rightarrow \Gamma(X, \OO_X)\mod
    \]
	is rarely an equivalence if $X$ is not affine.\footnote{For the reader 
	wondering about the use of `rarely', the point is that when one considers 
	the full 
	derived algebra of global functions, $\Gamma$ can still be an isomorphism 
	for non-affine varieties, e.g., for open subvarieties of affine varieties.}
	By contrast, for most finite dimensional $X$, e.g., any algebraic stack of 
	finite type, the functor \eqref{e:globsectshvcat} is an {equivalence}, a 
	phenomenon called 1-affineness.

	A basic conjecture in the subject, due to Gaitsgory, is that $\LS$ is also 
	1-affine; see \cite{raskinspectraldecomp}. This is a subtle assertion, as 1-affineness does not hold in
	general for ind-schemes or classifying stacks of group schemes of infinite 
	type. So, in some sense, one needs the two types of issues to `cancel 
	out' for $\LS$. This was partially addressed in an 
	important work of Raskin \cite{raskinspectraldecomp}, which 
	showed that 
	\eqref{e:globsectshvcat} is fully faithful for $X = \LS$ .

	\subsection{Properties of the correspondence I}

	\subsubsection{} Having dispensed with a discussion of definitions, 
	let us 
	get to the heart of the matter and discuss what $\mathbb{L}$ is expected to 
	look like. More precisely, let us list a few basic objects on each side 
	that are expected to be exchanged with one another. Combined with a 
	compatibility with parabolic induction, this is expected to largely pin 
	down the equivalence.

	The structure of our discussion will be as follows. We have listed the categories exchanged in Figure \ref{f:table}. For each row in the figure, we will spend a little 
	time describing the objects that appear on either side and motivate why one 
	should expect them to be 
	exchanged. Finally, we will see what matching intertwining operators 
	amounts to in these cases and, in particular, we will meet the promised 
	correction of \eqref{e:loclangnaive}.

	\begin{center}
		\begin{figure}

			\begin{tabular}{|c|c|} 
				
				\hline &
				

\\
			 \hspace{1cm}Categorical representations of $LG$  \hspace{1cm} 
				& 
				  \hspace{1cm}Sheaves of categories on $\LS$  \hspace{1cm} 
				\\  & \\  
				\hline &
				\\ 
	$\hspace{0cm} \Dmod(\on{Gr_G})$ 	& 
				$ \on{QCoh}( \on{pt}/\cG)$  
				\\ & \\ [2mm] 
				$\hspace{0cm}\Dmod(\Fl^{\on{aff}}_G)$ 	& $\hspace{0cm} \QCoh( 
				\widetilde{{\sN}}/\cG)$  \\ & \\ [2mm] 
				$\hspace{0cm}\Dmod(LG/LN, \psi) $	& $\hspace{0cm} 
				\QCoh(\LS)$  
				\\ & \\ [2mm] 
				$\hspace{0cm}\gkc\mod $	& $\hspace{0cm} \QCoh(\Op)$  \\  &
				\\ 
			\hline	
			\end{tabular}
			\caption{Some objects (conjecturally) exchanged under local 
			geometric Langlands.}
			
			\label{f:table}
		\end{figure}
	\end{center}

\subsubsection{}For the first two examples from Figure \ref{f:table}, it is 
clarifying to recall some 
basic features of $p$-adic representation theory. By smoothness, any 
representation $\pi$ of a $p$-adic group contains nontrivial invariants $\pi^K$ 
for all sufficiently small compact open subgroups $K$. For such a $K$,	if $\pi$ is 
irreducible,  the fixed vectors $\pi^K$ will be an irreducible 
representation of the Hecke algebra $\on{Dist}(K \backslash G(\Fqt) / K)$, 
which naturally acts on it by convolution. In this way, one obtains a bijection 
between irreducible representations of 
$G(\Fqt)$ containing $K$-invariant vectors and irreducible representations of 
this Hecke algebra. 

As one shrinks the subgroup $K$, the corresponding Hecke algebra grows, and 
this gives a 
basic measure of the complexity of a representation. Namely, for a given 
irreducible $\pi$, one can ask for a maximal $K$ for which one has nontrivial 
invariants, and, in some sense, the larger $K$ is, the less complicated $\pi$ is.

On the other side of local Langlands, the (Weil--Deligne variant of the) absolute 
Galois group of $\Fqt$ carries a canonical filtration by descending normal 
subgroups called the ramification filtration. For a given $\cG(k)$-valued Galois 
representation, one can measure its complexity by asking through what quotient 
of this filtration the action factors. A basic and striking organizing principle is that 
these two measures of complexity are exchanged under local Langlands.\footnote{More precisely, the depth of a representation, which is a certain nonnegative rational number defined by Moy--Prasad \cite{moyprasad}, should coincide with the slope of the corresponding Galois representation, as measured by the upper numbering ramification filtration.}

\subsection{Properties of the correspondence II: unramified representations} 
	
	\subsubsection{} From the preceding discussion, 
	it follows that the most easily understood representations of $p$-adic 
	groups should be those that have fixed vectors with respect to the maximal 
	compact open subgroup $G(\Fq[[t]])$.\footnote{We should note that, for 
	groups 
	other 
	than $GL_n$, there 
	are other non-conjugate maximal compact subgroups.} Similarly, the simplest 
	Galois representations should be those on which the action factors through 
	the final quotient of the ramification filtration, i.e., they extend to 
	$\check{G}$-local systems on the non-punctured disk. Explicitly, these 
	are precisely the representations which factor through the Galois group of 
	the residue field 
\begin{equation} \label{e:galunramified}
 1 \rightarrow I \rightarrow \on{Gal}( \Fqt^s / \Fqt) \rightarrow 
 \on{Gal}(\overline{\mathbb{F}}_q / 
 \Fq) \simeq \widehat{\mathbb{Z}} \rightarrow 1,
\end{equation}
where the superscript `$s$' denotes the separable closure. 
	
	Since one actually works with the Weil group, which by definition replaces 
	$\widehat{\mathbb{Z}}$ with $\mathbb{Z}$, unramified Galois 
	representations, up to isomorphism, are 
	specified by elements of $\check{G}$, up to conjugacy. That is, we follow the image 
	of $1 \in \mathbb{Z}$ under a map $\mathbb{Z} \rightarrow \cG$. In this simplest case, local 
	Langlands reduces to the Satake isomorphism between the spherical Hecke 
	algebra and the character ring of $\cG$, i.e., 
	\[    \on{Dist}(G(\Fq[[t]]) \backslash G(\Fqt) / G(\Fq[[t]])) \simeq 
	\on{Fun}(\cG(k) / \cG(k)),
	\]
	where the right hand side denotes the adjoint invariant polynomial 
	functions on $\cG$, viewed as a variety over $\on{Spec} k$ \cite{satake63}. By passing to 
	their characters, we obtain that 
	irreducible unramified 
	representations of $G(\Fqt)$ are 
	parametrized by (semisimple) conjugacy classes in $\cG(k)$, as desired.
	
	\begin{re} For 
	the reader 
	encountering these ideas for the first time, we emphasize that even the 
	commutativity of the appearing Hecke algebra is not a priori obvious. This is illustrative of how passing to bi-invariant distributions for a large 
	compact subgroup dramatically simplifies life.\end{re}
	
	\subsubsection{} 
	To prepare for the geometric version, we should mention that each 
	irreducible unramified 
	representation embeds, uniquely up to scaling, into the universal 
	unramified representation 
	\begin{equation} \label{e:fungr}
	    \Fun( G(\Fqt) /G(\Fq[[t]])).
	\end{equation}
	Moreover, to extract 
	each from the universal case, one takes the submodule of \eqref{e:fungr} 
	which transforms by the corresponding character under right convolution by 
	the spherical Hecke algebra.

	Finally, let us discuss the analogues  of the above in  local geometric 
	Langlands. On the automorphic side, we consider the category of 
	D-modules on the affine Grassmannian  
	\[
	   \Dmod(\Gr_G) = \Dmod(LG/L^+G)
	\]
	 as a substitute in a way that is hopefully familiar by now to the 
	 above universal 
	unramified representation. On the spectral side, any de Rham local system 
	on 
	the punctured disk that extends to the full disk is trivializable, i.e., we 
	obtain the stacky point
	\[
	    \on{pt}/\cG \rightarrow \LS
	\]
	and correspondingly view $\QCoh( \on{pt}/\cG)$ as a sheaf of categories 
	over $\LS$.

	We should also mention that one has analogues of not only the universal 
	case, but also of an individual irreducible unramified representation. We defer 
	a more detailed discussion until Section \ref{ss:locaff}. 
	
	\begin{re}We should note that, while the unramified case is the most basic, 
	it is completely fundamental in both arithmetic and geometric Langlands. 
	Namely, in the global Langlands correspondence, which we do not discuss in 
	detail here, any irreducible smooth representation of the adelic group 
	$G(\mathbb{A}_F)$ factorizes as a restricted tensor product of irreducible 
	representations of the $p$-adic groups $G(F_x)$ at each place, and 
 all but finitely many are unramified. To even make sense of such an infinite tensor product, one makes essential use of the fact that unramified representations have canonical generating lines, namely their spherical invariants. In the global geometric theory, almost 
 all that is known so far is for the analogue of everywhere unramified (or 
 possibly tamely ramified) automorphic representations. 
	\end{re}

	\subsection{Properties of the correspondence III: unramified principal 
	series} 
	\label{s:unramps}
	\subsubsection{} 
	Given the unramified 
	representations of $p$-adic groups, one can produce further representations 
	as follows. As in other parts of representation theory, one has 
	parabolic induction functors that allow one to build representations of $G$ 
	starting with representations of smaller groups, such as the Cartan $T$. 
	Feeding unramified representations into this construction yields an 
	important family of modules, the unramified principal series. Their simple subquotients are precisely the simple $G(\Fqt)$-modules with nontrivial 
	invariants for the compact open Iwahori subgroup $I$, i.e., the preimage of 
	the Borel of $B(\Fq)$ of $G(\Fq)$ under the projection $$G(\Fq[[t]]) 
	\rightarrow G(\Fq).\footnote{This theorem is often attributed to Borel and Matsumoto, and we have also heard it attributed to Bernstein or Casselman; see Proposition 2.4 of \cite{cass80spherical} and references therein.}$$That is, applying parabolic induction to unramified 
	representations yields the next simplest family of representations of 
	$p$-adic groups. 
	
	In spectral terms, such irreducible representations of $G(\Fqt)$ correspond to Galois representations 
	that are tamely ramified with unipotent 
	monodromy. To explain this, recall from the sequence \eqref{e:galunramified} the inertia subgroup 
	\[
	   I \simeq \on{Gal}( \overline{\mathbb{F}}_q(\!(t)\!)^s / \overline{\mathbb{F}}_q(\!(t)\!)),
	\]
	where $\overline{\mathbb{F}}_q$ denotes an algebraic closure of $\Fq$. The next quotient in the ramification filtration splits this as 
	\[
	    1 \rightarrow I_{\on{wild}} \rightarrow I \rightarrow I_{\on{tame}} \rightarrow 1. 
	\]
	The basic idea for $I_{\on{tame}}$, the tame inertia, is that it behaves like its characteristic zero counterpart, i.e., covers of a formal punctured disk over $k$. These covers are classified exactly like the covering spaces of a circle in topology, i.e., one adjoins an $n^{th}$ root of $t$ for a positive integer $n$. In positive characteristic, these are separable if and only if $n$ is coprime to the characteristic of $q$, and such covers exhaust all Galois extensions of degree coprime to $q$. Explicitly, one may choose an isomorphism 
	\[
	   I_{\on{tame}} \simeq \underset{\ell} \Pi \hspace{.7mm} \mathbb{Z}_\ell, 
	\]
	where $\ell$ runs over primes not equal to the characteristic of $\Fq$. A Galois representation is said to be tamely ramified if it factors through the quotient by $I_{\on{wild}}$, and it is further said to be unipotently monodromic if any element of $I_{\on{tame}}$, or equivalently a topological generator, is sent to a unipotent element of $\cG$.

	\subsubsection{} 
	To orient ourselves for the geometric discussion, let us review in slightly more detail how the relevant parabolic induction of $p$-adic representations goes. Given the Borel 
	$B$, whose Levi quotient is the abstract Cartan  $T$, one obtains a correspondence on rational points
	\[
	    T(\Fqt) \leftarrow B(\Fqt) \rightarrow G(\Fqt),
	\]
	where the leftward arrow is a surjection with kernel the rational points of 
	the unipotent radical $N(\Fqt)$, and the rightward arrow is a closed 
	embedding. This gives rise to a functor of parabolic induction, 
	\begin{equation} \label{e:parindpadic}
	    T(\Fqt)\mod \xrightarrow{\on{Res}} B(\Fqt)\mod \xrightarrow{\on{Ind}} 
	    G(\Fqt)\mod.
	\end{equation}
	Explicitly, one inflates a smooth representation $\pi$ of $T(\Fqt)$ to a 
	representation of $B(\Fqt)$ by having the unipotent radical act trivially, 
	and then one induces up to $G(\Fqt)$. This induction, as usual, is concretely given 
	by smooth sections of the corresponding vector bundle over the flag 
	manifold, i.e., $$G(\Fqt) \overset{B(\Fqt)} \times \pi \rightarrow 
	G(\Fqt)/B(\Fqt).$$

	In particular, it follows that any parabolic induction of an 
	unramified character of $T$ embeds into the parabolic induction of the 
	universal unramified representation, i.e. the universal principal series module
	\[
	    \Fun\big( G(\Fqt) / T(\Fq[[t]]) \cdot N(\Fqt)\big).
	\]
	%
	\begin{re} For a reader who has not explicitly encountered parabolic 
	induction before, it may be orienting to realize that the classification of 
	irreducible algebraic representations of reductive groups by their highest 
	weight is an example of parabolic induction. 
	
	Indeed, for the algebraic 
	groups $G,B,T$ and $N$ as above, one has adjunctions between their 
	categories of algebraic representations over any field
	\[
	    \on{Rep}(G) \overset{\on{Res}}{\underset{\on{Ind}}{\rightleftarrows}} 
	    \on{Rep}(B) 
	    \overset{\on{Coinv}_N}{\underset{\on{Res}}{\rightleftarrows}} 
	    \on{Rep}(T),
	\]
	where $\on{Res}$ denotes restriction, $\on{Coinv}_N$ denotes the functor of 
	$N$ coinvariants, and $\on{Ind}$ is the functor of coinduction.  
	Concatenating these, one obtains the functor of parabolic restriction, 
	which takes a $G$-module and returns its $N$-coinvariants, viewed as a 
	$T$-module, as well as its right adjoint functor of parabolic coinduction, 
	which takes a $T$-module, extends it to a $B$-module with a trivial 
	action of $N$, and coinduces this to $G$.
	
	 With this, the theorem of 
	the highest weight, which can be phrased as each irreducible $G$-module 
	having a line of $N$-coinvariants, in particular implies that every 
	irreducible $G$-module appears as a submodule of a parabolically coinduced 
	module. Moreover, in characteristic zero, by complete reducibility it 
	follows that every irreducible coincides with a parabolically coinduced module. Via the 
	explicit model of the coinduction of a character as sections of a line bundle, this recovers  
	the 
	Borel--Weil theorem about sections of line bundles on $G/B$.  \end{re}

	\subsubsection{} 	Let us turn to the geometric case. In this case, the 
	automorphic side will 
	similarly be D-modules on a semi-infinite flag manifold, i.e., 
	\[
	  \Dmod( LG / L^+ T \cdot LN),
	\] 
	which, as before, may be constructed from $\Dmod(\Gr_T)$ by a suitable 
	functor of 	parabolic induction similar to \eqref{e:parindpadic}. 
	
	\begin{re} Here is a technical comment regarding this category. It is as 
	yet unknown how to treat the space $LG/L^+T \cdot LN$ as an object of 
	algebraic 
	geometry. However, its category of D-modules may be defined by enforcing 
	smooth descent, i.e., by taking the $\Dmod(L^+T \cdot LN)$ invariants  of 
	$\Dmod(LG)$. It is known, by an argument due to Raskin  \cite{raskin2ndadjoint}, that this canonically identifies 
	with the coinvariants and also with the (co)invariants for the 
	Iwahori subgroup, i.e., D-modules on the
	affine flag variety. This parallels the result of Borel and Matsumoto in the $p$-adic setting, and in fact is proven similarly.
	\end{re}

	 \subsubsection{} On the spectral side, one again applies a parabolic 
	 induction functor, 
	 i.e., pull-push along the correspondence 
	 \[
	    \on{LocSys}_{\check{T}} \leftarrow \on{LocSys}_{\check{B}} \rightarrow 
	    \LS
	 \]
	  to $\QCoh(\on{pt}/\check{T})$. Here we mean pullback and pushforward of 
	  sheaves of categories, as in Remark \ref{r:pushpull2coh}.\footnote{More 
	  carefully, the promised correction of Equation \eqref{e:loclangnaive} 
	  will involve passing from 
	  quasicoherent sheaves of categories to certain ind-coherent sheaves of 
	  categories, but this correction will not affect this calculation.}

	   Let us walk through how to compute this. 
	  For the first step, i.e., the pullback, note that we have a Cartesian 
	  diagram 
	  \[
	  \xymatrix{ \on{LocSys}_{\check{N}} \ar[r] \ar[d] & 
	  \on{LocSys}_{\check{B}} \ar[d] \\ \on{pt} \ar[r] & 
	  \on{LocSys}_{\check{T}}.}
	  \]
	  This simply says that a $\check{B}$ local system whose 
	  associated $\check{T}$ local system is trivialized is the same as 
	  an $\check{N}$ local 
	  system. Even more colloquially, if the transition matrices have 
	  ones along the 
	  diagonal, it's not a $\check{B}$-bundle but an $\check{N}$-bundle! By 
	  the unipotence of $\check{N}$, the moduli of $\check{N}$-local systems 
	  on the punctured disk is simply $\check{\fn} / \check{N}$, via the 
	  map 
	  \begin{equation} \label{e:rsmap}
	    \check{\fn} / \check{N} \xrightarrow{\sim}   
	    \on{LocSys}_{\check{N}}, 
	    \quad  X \mapsto d + 
	    X \cdot 
	    \frac{dt}{t}.
	  \end{equation}
	 From this discussion, it follows that the pullback is given by 
	 \begin{equation}
	 \QCoh(\on{LocSys}_{\check{B}}) \underset{ \QCoh(\on{LocSys}_{\check{T}})} 
	 \otimes \QCoh(\on{pt}/\check{T}) \simeq \QCoh( \on{LocSys}_{P}  \underset{ 
	 \on{LocSys}_{\check{T}}} \times \on{pt}/\check{T}  ) \simeq 
	 \QCoh(\check{\fn} / \check{B} ). 	 \end{equation} 

	 Finally, the pushforward again yields $\QCoh(\check{\fn} / \check{B} )$ 
	 viewed as a sheaf of categories on $\LS$ via the projection 
	 \begin{equation} \label{e:projection}
	     \check{\fn} / \check{B} \rightarrow \on{LocSys}_{\check{B}} 
	     \rightarrow \LS. 
	 \end{equation}
	It is clarifying to slightly rewrite this purely in terms of $\cG$. To do 
	so, we may tautologically rewrite our stack in terms of the cotangent 
	bundle of the
	 flag variety, i.e., 
	\[
	      \check{\fn} / \check{B} \simeq (\cG \overset{ \check{B}} 
	      \times \check{\fn}) / \cG =: 
	      \widetilde{\sN} / \cG.
	\]
	Explicitly, $\widetilde{\sN}$ is the moduli space of pairs $({\check{B}}', X)$, where ${\check{B}}'$ is a Borel subgroup and $X$ is an element of its unipotent radical. With this, the projection \eqref{e:projection} may be written via the Springer 
	resolution, i.e., as a composition 
	\[
	    \widetilde{\sN}/ \cG \rightarrow \cfg / \cG \rightarrow \LS.
	\]
	Here, the first map sends a pair $(\check{B}', X)$ to $X$, and  the second 
	map is defined similarly to \eqref{e:rsmap} (but is very 
	far from an isomorphism).

	\begin{re} As we will review in Section \ref{s:2realz} below, the relationship 
	between principal series representations of $G$ and coherent sheaves on 
	the Springer resolution for $\cG$ showed up first in the $p$-adic theory. 
	\end{re}

	\subsubsection{} The two preceding examples, namely the first two rows 
	of Figure \ref{f:table}, concern objects of low ramification. These cases 
	of the local geometric Langlands correspondence, while the most basic,  
	are already serious results. 
	
	Indeed, on the one hand, as we mentioned 
	before, they are sufficient local inputs for the simplest cases of the 
	(conjectural) global geometric 	Langlands correspondence,  where one 
	studies D-modules on the moduli space of $G$-bundles on a global curve, possibly with 
	parabolic reductions at finitely many points.

	On the other hand, they have important applications within usual 
	representation theory. A notable example is their implications for the  
	representation theory of critical level affine Lie algebras, which we will 
	discuss in more detail in Section \ref{ss:locaff} below. Another key application, after passing to sheaves with modular coefficients,\footnote{Strictly speaking, these fit more naturally into the local \'{e}tale or Betti versions of local geometric Langlands, cf. Section \ref{ss:variants} below, as they have natural analogues with varying coefficients.} is to the modular representation theory of 
	reductive groups; see \cite{finkelbergmirkovicconjecture}, \cite{mirkovic-vilonen}, \cite{acharrider}, \cite{richewilliamsontilting}, \cite{bezgaitsmiriwahotiwhittaker}, \cite{bezricheriderregcentralizer}, \cite{bezricheprinblock} and references therein, as well as the surveys \cite{willimasonsurvey1}, \cite{williamsonsurvey2}.

Given the importance of these developments, it is rather remarkable that they 
only skim the surface of the full conjecture, i.e., equate only very small 
subcategories of each side. Indeed, in spectral 
	terms, the previous cases concern the full subcategory of $\TwoCoh(\LS)$ 
	supported on 
	the formal neighborhood of local systems with regular singularities and 
	unipotent monodromy, which is the smallest in an infinite increasing sequence of loci whose union exhausts $\LS$.

The next two examples of the correspondence which we discuss, by contrast,
concern categorical loop group representations and sheaves of categories on 
$\LS$ with arbitrary ramification. As 
such, these are rather spectacular conjectures, among the deepest currently 
available in this part of
representation theory, and their resolutions likely require new ideas beyond 
the 
standard 
toolkit of geometric representation theory.

\subsection{Properties of the correspondence IV: Whittaker models}
	
\subsubsection{} In many situations in representation theory, 
e.g., algebraic representations of reductive groups or $(\fg, B)$-modules, 
irreducible representations admit a canonical generating line, namely the 
highest weight vectors. 

By contrast, for the $p$-adic group $G(\Fqt)$, most irreducible modules do not 
admit a canonical generating line, roughly due to the noncommutativity of $\on{Dist}(K 
\backslash G / K)$ for all sufficiently small $K$. 

\begin{re} For unramified representations, there is a canonical generating 
line, namely the $G(\Fqt)$-invariants. This basic miracle, i.e., the 
commutativity of the spherical Hecke algebra, is crucial in the subject, e.g., 
for the tensor product factorization of adelic representations and the 
reduction of the study of unramified automorphic representations to the 
study of Hecke eigenfunctions and eigensheaves. \end{re}

 However, if one asks instead for a canonical line not inside  $\pi$, but rather 
 as a quotient of $\pi$, the situation is better. Namely, there is a canonical 
 functor from $G(\Fqt)\mod$ to vector spaces, which sends each irreducible to 
 either a line or zero. This is the functor of Whittaker coinvariants and may 
 be understood as follows. 
 
 To approximate the definition, we may first consider a seemingly direct
 analogue of highest weight theory, namely 
 the Jacquet functor of $N(\Fqt)$-coinvariants
 \[
       \on{Jac}: G(\Fqt)\mod \rightarrow \on{Vect}, \quad \pi \mapsto 
       \pi_{N(\Fqt)}. 
 \]
While this is very useful, it sends many interesting irreducible representations to 
zero\footnote{This phenomenon is more or less the existence of cuspidal 
representations, i.e., of simple $G(\Fqt)$-modules that may not be found 
inside parabolically induced modules. We emphasize that this is a basic fact of life 
in the $p$-adic setting and in the simpler setting of finite groups of Lie 
type, which is absent in the algebraic representation theory of reductive groups.} 
and others to greater than one dimensional vector spaces. 

To correct for this, one twists the construction by a generic character of 
$N(\Fqt)$. Explicitly, let us index the simple roots by $i \in I$, and consider 
the 
composition 
\begin{equation} \label{e:whitchar}
\psi:  N(\Fqt) \rightarrow (N/[N,N])(\Fqt) \simeq \underset{i} \Pi 
\hspace{.2mm} 
\Fqt \xrightarrow{\on{sum}} \Fqt 
\xrightarrow{\on{res}} \Fq,
\end{equation}
where sum is the addition of Laurent series and res sends a Laurent series to 
its residue, i.e., the coefficient of $t^{-1}$. By composing this with a 
suitably generic character $\mathbb{F}_q \rightarrow k^\times$, one obtains a 
character of $N(\Fqt)$, which we again denote by $\psi$, and the Whittaker 
coinvariants are given by 
\[
   \on{Whit}: G(\Fqt)\mod \rightarrow \Vect, \quad \pi \mapsto \pi_{N(\Fqt), 
   \psi}. 
\]
This sends simple modules to either a line or zero, with the crucial property 
that many cuspidal representations, e.g., all for $G = GL_n$, are sent to 
lines. 

A basic property here is that $\on{Whit}$ admits a right adjoint, given by 
sending $k$ to the space of Whittaker functions 
\[
     \Fun( G(\Fqt) / N(\Fqt), \psi),  
\]
i.e., the coinduction of $\psi$ to $G(\Fqt)$.

\subsubsection{}
The above constructions admit the following geometric analog. First of all, the 
homomorphism of topological groups \eqref{e:whitchar} lifts to a map of group 
ind-schemes 
\[
      \psi: LN \rightarrow L(N/[N, N]) \simeq \underset i \Pi \hspace{.2mm} 
      L\mathbb{G}_a \rightarrow L\mathbb{G}_a \rightarrow \mathbb{G}_a.
\]
From here, we would like to pull back from $\mathbb{G}_a$ the analogue of a 
generic character. Recalling the analogy between the trivial representation of 
$G(\Fq)$ and the tautological categorical action of $G$ on $\Vect$, we 
therefore need a generic action of $\mathbb{G}_a$ on $\on{Vect}$. An action of an 
algebraic group $H$ on $\Vect$ is the same data as a character sheaf on $H$, i.e., 
a D-module $\chi$ equipped with a suitably associative isomorphism 
\[
     \mu^! \chi \simeq \chi \boxtimes \chi,
\]
where $\mu: H \times H \rightarrow H$ denotes the multiplication map. We emphasize that in contrast to the broader usage of character sheaf in Section \ref{ss:charsheaves}, which was the analogue of the character of a finite dimensional representation of $G(\Fq)$, the present usage of character sheaf is the analogue of a character, i.e. one dimensional representation, of $G(\Fq)$.

\begin{re} Slightly 
more informally, if for any closed point $h$ of $H$ we write $\chi_h$ for the corresponding 
$!$-stalk of $\chi$, the above is giving isomorphisms of stalks 
\[
\chi_{h_1 h_2} \simeq \chi_{h_1} \otimes \chi_{h_2}, \quad \text{for } h_1, 
h_2 \in H.
\]
That is, a character sheaf gives a homomorphism into lines, i.e., the 
invertible elements of $\Vect^\heartsuit$, just as a character gives a homomorphism into 
nonzero scalars, i.e., the invertible elements of $k$. 
\end{re} 	

With this, the analogue of a generic character of $\Fq$ is the exponential 
D-module $\on{exp}(z)$ on $\mathbb{G}_a$. Explicitly, if we write $z$ for the 
standard coordinate on $\mathbb{G}_a \simeq \mathbb{A}^1$, this D-module stores 
the differential equation satisfied by the exponential function $e^z$, i.e., 
\[
\on{exp}(z) \simeq D_{\mathbb{G}_a} / ( D_{\mathbb{G}_a} \cdot (\partial_z - 
1)), 
\]
and its multiplicativity essentially follows from the identity $e^{z_1 + z_2} = 
e^{z_1} \cdot e^{z_2}$. With this, we write $\psi$ for the corresponding 
character sheaf on $LN$, i.e., $\psi^! \on{exp}(z)$, and we may again form the 
functor of Whittaker coinvariants 
\[
    \on{Whit}: \Dmod(LG)\mod \rightarrow \Vect, \quad \sC \mapsto \sC_{LN, 
    \psi}
\]
and the category of Whittaker D-modules 
\begin{equation} \label{e:catwhitsheaves}
    \Dmod(LG / LN, \psi). 
\end{equation}
Before turning to its spectral counterpart, we would like to make a few 
orienting comments. 

\begin{re}Let us try to provide three ways of thinking about the analogy 
between 
the exponential D-module on $\mathbb{G}_a(k)$ and a generic character of 
$\mathbb{G}_q(\Fq)$. 

The first way is rather down to earth. For simplicity, let us suppose that 
$\mathbb{F}_q$ coincides with its prime subfield $\mathbb{F}_p$. In this case, 
if we write $e^{2\pi i / p}$ for a nontrivial $p$th root of unity in 
$k^\times$, a concrete choice of generic character is 
\[
   \mathbb{F}_p \rightarrow k^\times, \quad a \mapsto e^{ a\cdot 2\pi i / p }, 
\]
and any other differs by multiplication by a nonzero element of $\mathbb{F}_p$. 
That is, a generic character looks like an exponential. The same holds, {
mutatis mutandis}, for a nontrivial extension $\Fq$, where one first uses the trace 
map $\Fq \rightarrow \Fp$. 

The second way is also rather practical. Namely, for a $k$ vector space $V$, its 
character sheaves are in canonical bijection with closed points of $V^*$, where 
to a covector $\xi$ one attaches the character sheaf $\xi^! \on{exp}(z)$. In 
particular, for $V = \mathbb{G}_a$, there is only one nontrivial choice of 
character sheaf up to dilation by $\mathbb{G}_m$, just as there is only one nontrivial character 
for $\Fq$ up to dilation by $\Fq^\times$.

The third way is more general and goes via the functions-sheaves correspondence. To 
begin, a beautiful feature of positive characteristic is an interesting 
self-isogeny of any commutative group scheme $H$ over $\Fq$, the Lang isogeny, 
which takes 
the form 
\[
   1 \rightarrow H(\Fq) \rightarrow H \rightarrow H \rightarrow 1,
\]
i.e., it cuts out as its kernel the finite group of rational points. For $H = 
\mathbb{G}_a$, this is concretely the short exact sequence 
\[
   1 \rightarrow \mathbb{G}_a(\Fq) \rightarrow \mathbb{G}_a \xrightarrow{x 
   \mapsto x^q - x} \mathbb{G}_a \rightarrow 1. 
\]
Pushing forward the constant sheaf along this isogeny and taking its direct 
summands yields a family of character sheaves on $\mathbb{G}_a$ indexed by 
irreducible representations of $\mathbb{G}_a(\Fq)$. Upon taking trace of 
Frobenius, these character sheaves yield exactly the irreducible characters of 
$\mathbb{G}_a(\Fq)$. That is, characters of $\Fq$ have canonical lifts to 
character sheaves on $\mathbb{G}_a$. 

For these character sheaves, viewed simply as \'etale local systems on 
$\mathbb{A}^1$, one can ask about their behavior at $\infty$, i.e., after compactification  
to $\mathbb{P}^1$. One finds that they are wildly ramified at this remaining 
point. The correspondence with $\on{exp}(z)$ in characteristic zero then fits into a broader analogy between wildly ramified Galois 
representations and de Rham local systems with irregular singularities. 
\end{re}

\begin{re} Let us comment on the definition of the category of 
Whittaker D-modules $\Dmod(LG / LN, 
\psi)$. This is defined by taking the $(LN, \psi)$-coinvariants of $\Dmod(LG)$, similarly to the categories of D-modules on semi-infinite flag manifolds. It is also identified with the $(LN, 
\psi)$-invariants by an 
important theorem of Raskin \cite{whit}. While an identification of Whittaker invariants and coinvariants does not literally hold in the $p$-adic setting, it is analogous to the compact approximation of the Whittaker model by Rodier \cite{rodier}.   
\end{re}

\subsubsection{}Let us now describe the expected spectral counterpart to the 
Whittaker 
	model. This 
	is simply $$\QCoh(\LS),$$ i.e., all quasicoherent sheaves on the moduli of 
	local systems. 
	
	The basic intuition here again can be taken from the 
	$p$-adic setting. Namely, each irreducible representation $\pi$ of 
	$G(\Fqt)$ has, at 
	most, a line of Whittaker coinvariants. It will be convenient to work instead with
	\[
	   \Hom( \on{ind}_{N(\Fqt)}^{G(\Fqt)} k_{-\psi}, \pi) =: \ell,
	\]
	where $\on{ind}$ is the functor of compactly supported induction. This identifies with the Whittaker coinvariants of the contragredient representation of $\pi$ and, in particular, is again a line or zero.

	On the other hand, under 
	local Langlands one associates to $\pi$ a Galois representation $\sigma$ and imagines $\pi$ as `sitting over' $\sigma$ in the space $\on{Gal}_{\cG}$ 
	of Galois 
	representations. 
	
	In the modern formulations of arithmetic local Langlands as an 
	equivalence of categories, cf. Section \ref{ss:variants} below, it is expected that any $\pi$ admitting nonzero Whittaker coinvariants is exchanged with a one dimensional skyscraper sheaf $\mathbb{L}(\pi)$ supported over $\sigma$. Moreover, this line canonically associated to $\pi$ is expected to simply be $\ell$. So, combining the identifications 
	\[
	    \Hom( \on{ind}_{N(\Fqt)}^{G(\Fqt)} k_\psi, \pi) \simeq \ell \simeq 
	    \Hom( \OO_{\on{Gal}_{\cG}},  \mathbb{L} \pi ), 
	\]
	we see that $\mathbb{L}$ ought to exchange the (compactly) induced Whittaker module and the structure sheaf $\OO_{\on{Gal}_{\cG}}$. 
	
	In the geometric setting, we have already replaced the Whittaker module by 
	the 
	corresponding space of D-modules $\Dmod(LG/LN, \psi)$.\footnote{Note that 
	the difference between 
	induction and 
	compact induction disappears, due to Raskin's identification of Whittaker 
	invariants and coinvariants \cite{whit}.} Similarly, we replace the structure sheaf 
	$\OO_{\on{Gal}_{\cG}}$ by $\QCoh(\LS)$, which again corepresents global 
	sections, but now of sheaves of categories.

    As we will shortly discuss when we get to intertwining operators, this 
    correspondence leads to many striking and concrete conjectures in geometry and representation theory. However, for now let us 
    conclude our 
    discussion of this row of Figure \ref{f:table} with a couple remarks. 
    
    \begin{re} Unlike in the $p$-adic setting, where not every tempered irreducible representation 
    admits a Whittaker model, the expectation in local geometric Langlands is 
    that every tempered $\Dmod(LG)$ representation does. 

    The basic idea is the following. For a point $x$ of a stack $X$ with 
    isotropy group $I_x$, the abelian category of quasicoherent sheaves on $X$ 
    scheme theoretically supported at $x$ identifies with $\on{Rep}(I_x)$, and 
    in particular the 
    skyscraper will not generate unless $I_x$ is unipotent. Under local 
    arithmetic Langlands, this matches Shahidi's generic packet conjecture, which states that, within a tempered L-packet,\footnote{Beyond the case of $GL_n$, the map from 
    irreducible representations of $G(\Fqt)$ to Galois representations is not a 
    bijection, but instead finite-to-one, and the fibres are called L-packets.} exactly one irreducible representation will admit a 
    Whittaker model \cite{shahidi90}. 
    
     By contrast, the skyscraper module $\Vect$ for $\QCoh(X)$ will generate 
     the subcategory of $\TwoCoh(X)$ supported at $x$, and in this sense every 
     irreducible tempered module is generic. Despite this, it 
     will not be the case that $\Dmod(LG/LN, \psi)$ is a generator for 
     $\Dmod(LG)$, for issues having to do with the promised 
     correction to the naive formulation \eqref{e:loclangnaive}. 
    \end{re}

	\subsection{Properties of the correspondence V: Kac--Moody representations}

	\subsubsection{} We have arrived at the last row of Figure \ref{f:table}, 
	and here something rather remarkable occurs. Namely, we will meet a 
	conjectural property of the correspondence, proposed by Frenkel--Gaitsgory \cite{frenkgaitslocallangds}, 
	which concerns Kac--Moody representations. As such, it has no $p$-adic 
	counterpart and is an aspect of local Langlands which is special to 
	the geometric context. 
	
	\begin{re} We should mention that the importance of Kac--Moody 
	representations in geometric Langlands first appeared in the global 
	story. Namely, as we alluded to Section \ref{ss:hitchinhecke}, a method of construction 
	for automorphic D-modules, due to Beilinson--Drinfeld, is via localization of 
	Kac--Moody representations. This may be understood as a global analogue of the work of Frenkel--Gaitsgory and was the first instance where the relationship to 
	conformal field theory provided tools unavailable in the arithmetic 
	setting. \end{re}
	
	\subsubsection{}  We would first like to give the reader a basic feel for 
	Kac--Moody algebras. Their study for a general connected reductive group 
	essentially reduces to the case when $\fg$ is simple. So, to avoid mostly 
	notational distractions, we shall assume that $G$ is a simply connected, almost simple 
	group.

	First, we begin with the loop algebra 
	\[
	   \fg(\!(t)\!) \simeq \fg \underset  k \otimes k(\!(t)\!).
	\]
	This is simply the extension of scalars of $\fg$ to the field of Laurent 
	series. Explicitly, the bracket takes the form 
	\[
	   [ X \otimes f, Y \otimes g] = [X,Y] \otimes fg, \quad X,Y \in \fg, \quad 
	   f, g \in k(\!(t)\!). 
	\]
	It is clarifying to think of it as the space of maps from 
	the punctured disk $\sD^\times$ to $\fg$, where the Lie bracket is 
	performed pointwise. That is, it is the infinitesimal symmetries 
	of the trivial $G$-bundle on the punctured disk. 
	
	Let us regard the loop algebra as an infinite dimensional Lie algebra over $k$. A basic fact of life is that many $k$ vector spaces that naively seem like 
	$\fg(\!(t)\!)$-modules in fact carry an action of a central extension, i.e., of an affine Lie algebra
		\[
	0 \rightarrow k \cdot \mathbf{1} \rightarrow \gk \rightarrow 
	\fg(\!(t)\!) \rightarrow 0.
	\]

	\begin{re} This may be understood as an instance of the general phenomenon 
	that quantization of field 
	theories introduces anomalies for their local symmetries. \end{re}
	
	Fortunately, all the possible central extensions of the loop algebra admit 
	a simple 
	classification -- they are in canonical bijection with 
	the line of invariant bilinear forms on the finite dimensional algebra 
	$\fg$. 
	Explicitly, to such a form $\kappa \in (\fg^* \otimes \fg^*)^G$, one 
	associates the central extension 
	\[
	 0 \rightarrow k \cdot \mathbf{1} \rightarrow \fgk \rightarrow \fg(\!(t)\!) 
	 \rightarrow 0,
	\]
	where the bracket of $\mathbf{1}$ with anything is zero, by centrality, and 
	the other brackets are modified by 
	\[
	[X \otimes f, Y \otimes g] = [X,Y] \otimes fg - \kappa(X,Y) \on{Res} (f 
	\cdot dg). 
	\]

\subsubsection{} For a level $\kappa$, let us denote the category of smooth 
$\fgk$-modules on 
which the central element $\mathbf{1}$ acts via the identity by 
$
  \fgk\mod, 
$
cf. Example \ref{ex:kmmods}. If $\kappa$ is integral, meaning the extension of the 
loop algebra arises from a central extension of $LG$, then, 
as discussed in Example \ref{ex:kmmods}, this category carries an action of 
$\Dmod(LG)$, and we therefore shall focus 
on this case. 

Recalling that $G$ was assumed to be simple and simply connected, the integral levels form a 
lattice$$ \mathbb{Z} \hookrightarrow (\fg^* 
\otimes \fg^*)^G,$$generated by the so-called basic form, which explicitly gives the short coroots squared length two. 

\subsubsection{}  A basic observation in the study of representations of 
affine Lie algebras is that as one varies the level $\kappa$, the 
corresponding representation theory has two `phases.'

Namely, for all sufficiently negative definite $\kappa$, the representation 
theory behaves similarly. For example, all Verma modules have finite length, 
and the character formulas for their simple quotients across these levels take 
essentially the same form \cite{kashiwaratanisakinegativeintegrallevel}. 

Likewise, for all positive definite $\kappa$, the 
representation theory behaves similarly as one varies the level. Here, the Verma modules now have infinite length, and the character formulas for 
their 
simple quotients do not substantially vary with the level \cite{kashiwara90}. 

We may therefore anticipate that there is a special level $\kappa_c$, the critical level, 
where the representation theory undergoes a phase transition. Naively, one 
might expect this central point to be the zero form. Instead, due to 
an 
affine analogue of the $\rho$-shift in the representation theory of $\fg$, it is 
minus half of the Killing form, i.e., 
\[
\kappa_c = - \frac{1}{2} \cdot \kappa_{\on{Killing}}.
\]

\subsubsection{} As one might expect, the behavior of the representation theory at 
the phase transition, i.e., at critical level, displays many features not 
present at other levels. 

A fundamental such property was discovered by Feigin and Frenkel \cite{feiginfrenkel91}. Namely, 
while the center of the enveloping algebra\footnote{As a technical remark, we 
mean the center of the appropriate completed enveloping algebra, or 
equivalently the 
Bernstein center of the abelian category of smooth $\fgk$-modules.} at 
noncritical levels is trivial, 
i.e., only scalars, they showed that the center at critical level was nontrivial.  
Crucially, they further identified it with the algebra of functions on the moduli space  $\Op$ of $\cG$-opers on the punctured disk. While we will momentarily recall its precise definition, at first pass 
the main point is that representation theory at the critical level has local 
Langlands duality, particularly the geometry of $\cG$-connections, hardwired into it.

\subsubsection{} To first approximation, a $\cG$-oper consists of a $\cG$ bundle with a 
connection and a $\cB$ reduction, which are suitably transverse. 

Let us spell this out explicitly. Recall 
that to speak of Langlands dual groups, one fixes pinnings of $G$ and $\cG$ 
and, in particular, Chevalley generators $f_i, \check{\alpha}_i, e_i, i \in I$, 
of 
$\check{\fg}$. If we write $f$ for the element of $\cfg$ such that 
\[
  (f, 2 \check{\rho}, \underset i \Sigma \hspace{.2mm} e_i )
\]
form an $\fsl_2$ triple, then the space of $\cG$-opers is the moduli space of 
connections 
\begin{equation} \label{e:defopers}
   \Op \simeq \{ d  + f dt + \check{\fb}(\!(t)\!)dt \} / L\check{N}.
\end{equation}
 By thinking of the trivial bundle as equipped with its 
tautological $\cB$-reduction, the term $f dt$ in the above formula implements 
the aforementioned transversality of the connection and the reduction. For 
alternative formulations of the definition of an oper, e.g., with the advantage of visibly being independent of the coordinate $t$, see for example 
\cite{BD}.

The formula \eqref{e:defopers} may be read as taking a certain affine subspace of all connections 
$d + \cfg(\!(t)\!)dt$ on the trivial $\cG$ bundle and quotienting by the 
adjoint 
action of $L\check{N}$.
Indeed, \eqref{e:defopers} is the Hamiltonian reduction of the affine space $d 
+ \cfg(\!(t)\!)dt$ with respect 
to a generic character of $L\check{N}$ and, in particular, 
carries a canonical Poisson structure. Gelfand--Dikii and Drinfeld--Sokolov \cite{drinfeldsokolov} studied its Poisson geometry in the context of integrable hierarchies 
prior to its appearance in 
geometric Langlands; for $\fg = \fsl_2$ this recovers the celebrated KdV hierarchy.

Here are two orienting remarks for a reader who has not encountered this moduli 
space before. 

\begin{re} Although $\Op$ is presented as space of connections modulo the 
action of a loop group, similarly to $\LS$, its geometry is astronomically 
simpler. In fact, it is simply an infinite dimensional affine space, and, in 
particular, its points do not have automorphisms.  
\end{re}

\begin{re}  The previous remark is typically proven by showing every $L\check{N}$ 
orbit is simply transitive, i.e., without nontrivial stabilizers, and exhibiting 
an explicit slice intersecting each orbit exactly once. 

It can be clarifying to note that this whole story, including the above argument, is an affine analogue of a 
work of Kostant \cite{kostantslice}. Namely, he showed that one has a canonical isomorphism 
\[
   \cfg /\!/ \cG \simeq (f + \check{\fb}) / \check{N},
\]
i.e., that a similar finite dimensional Hamiltonian reduction of $\cfg$ cuts out 
a copy of the Poisson center. Moreover, he showed the natural quantizations of 
both sides are matched, namely that the center $Z(\cfg)$ may be cut out from  
$U(\cfg)$ via a quantum Hamiltonian reduction with respect to $\cN$. 
	\end{re}

\subsubsection{} Via the projection map $\Op \rightarrow \LS$, where one 
passes from an oper 
to its underlying connection, we may view $\QCoh(\Op)$ as a sheaf of categories 
over $\LS$.

With this, it is expected that $\gkc\mod$ is exchanged with 
$\QCoh(\Op)$ under local geometric Langlands, so that the tautological monoidal 
functor
\[
\QCoh(\Op) \rightarrow \Hom_{\TwoCoh(\LS)}(\QCoh(\Op), \QCoh(\Op)) 
\]
is exchanged, via the Feigin--Frenkel isomorphism $\OO_{\Op} \simeq Z(\gkc)$, 
with the central action 
\[
  Z(\fgk)\mod \rightarrow \Hom_{\Dmod(LG)\mod}(\gkc\mod, \gkc\mod).
  \]

 In the next subsections, we will 
discuss several interesting conjectures this correspondence yields about representation theory 
at critical level. 

However, we would like to already mention here one particularly striking consequence. 
It is known that any $\cG$-connection $\sigma$ admits an oper structure $\chi$, 
i.e., the projection $\Op \rightarrow \LS$ is surjective on $k$-points \cite{frenkelzhuoper}. By the 
above, it formally follows that the `irreducible' categorical representation 
$\sC_{\sigma}$ of the loop group, which corresponds to the skyscraper $\Vect$ 
over $\sigma$, is given by 
\begin{equation} \label{e:kmcenc}
  \gk\mod \underset{\QCoh(\Op)} \otimes \Vect,
\end{equation}
i.e., the category of Kac--Moody representations with central character $\chi$.

On the one hand, this gives a uniform conjectural construction of all the 
$\sC_{\sigma}$ via Kac--Moody representations, i.e., the spectral to automorphic 
direction of the local geometric Langlands correspondence. On the other hand, this 
predicts that for any two opers $\chi, \chi'$ lying over $\sigma$, the 
corresponding 
categories of Kac--Moody modules are equivalent, which is by itself a deep 
conjecture in representation theory.

\begin{re} In fact, due to the correction to the naive formulation 
\eqref{e:loclangnaive} of local geometric Langlands, which we 
have yet to explain, the category $\eqref{e:kmcenc}$ is really a full subcategory of 
$\sC_{\sigma}$. In practice, this means that 
one obtains $\sC_{\sigma}$ by 
renormalizing \eqref{e:kmcenc},  cf. Section \ref{ss:locaff} for 
further discussion. 
\end{re}

\subsubsection{} This completes our introductory discussion of the objects appearing in Figure \ref{f:table}. However, before going further, we should emphasize that there are many more expected compatibilities coming from the arithmetic story and connections with physics, cf. Section \ref{ss:variants} below. 

We cannot help but mention one particularly spectacular example of this. Braverman--Finkelberg have proposed for $G = GL_n$ a candidate for the Langlands dual to the convolution algebra $\Dmod(LGL_n)$ itself \cite{bravermanfinkelbergconjecturecoloumb}.  This is based on physical considerations, ultimately going back to a striking formula for computing S-dual boundary conditions for four dimensional $\mathscr{N}=4$ Yang--Mills theory discovered by Gaiotto-Witten \cite{gaiottowitten}. Explicitly, they conjecture that $\Dmod(LGL_n)$ is exchanged under local Langlands duality with coherent sheaves on the moduli space parametrizing degenerate flags of local systems, i.e., local systems $\sL_i$ of rank $i$, $1 \leqslant i \leqslant n$, along with not necessarily injective maps of connections $$\sL_i \rightarrow \sL_{i+1}, \quad 1 \leqslant i < n.$$An arithmetic analogue of this conjecture has since been formulated by Zhu \cite{zhufamilies}.

\subsection{Intertwining operators}

Next, we would like to explain how considering the intertwining operators 
between the categories appearing in Figure \ref{f:table} recovers several 
fundamental theorems and produces many striking conjectures. 

At this point, the promised technical correction to the naive formulation 
\eqref{e:loclangnaive} of the local geometric Langlands correspondence 
naturally asserts itself in the following manner.

\subsubsection{} Recall the simplest representations in the $p$-adic 
setting are the unramified ones. Pleasantly, we will already meet the necessary correction to the 
naive formulation \eqref{e:loclangnaive} of the geometric conjecture in the unramified case. 

To see this, recall that local Langlands ought to exchange the 
D-modules on the 
affine 
Grassmannian $\Dmod(\Gr_G)$ with quasicoherent sheaves on the trivial bundle 
$\QCoh(\on{pt}/\cG).$ In particular, their endomorphisms ought to be 
canonically identified, so let us see what they are. 

On the automorphic side, recalling that $\Gr_G$ is the coset space $LG/L^+G$, 
as in Section \ref{s:catrepgrp} we have 
\begin{equation} \label{e:dersat1}
 \Hom_{\Dmod(LG)\mod}( \Dmod(\Gr_G), \Dmod(\Gr_G)) \simeq \Dmod(L^+G \backslash 
 LG / L^+G), 
\end{equation}
also known as the derived Satake category of bispherical sheaves. As a 
potential warning for readers encountering this for the first time, this 
category is not the derived category of the corresponding abelian category, but 
instead richer. 

On the spectral side, using that the map $\cfg / \cG \rightarrow \LS$ is 
formally \'etale at the trivial bundle, we may compute, as in Example 
\ref{e:fibprodvar}, that 
\begin{equation} \label{e:dersat2}
\Hom_{\TwoCoh(\LS)}(\QCoh(\on{pt}/\cG), \QCoh(\on{pt}/\cG)) \simeq \QCoh( 
\on{pt}/\cG \underset{\LS} \times \on{pt}/\cG) \simeq \QCoh( \on{pt} 
\underset{\cfg} \times 
\on{pt} / \cG).
\end{equation}
A couple points of explanation are in order. First, for the last equivalence in \eqref{e:dersat2}, one may use the following general observation. Suppose one has  varieties $X, Y,$ and $Z$ acted upon by an algebraic group $K$. Then, given a $K$ equivariant correspondence $$X \rightarrow Z \leftarrow Y,$$one has a canonical isomorphism of quotient (derived) stacks
\[
      X/H \underset{Z/H} \times Y/H \simeq (X \underset Z \times Y)/H.
\]
We will use this fact repeatedly going forwards without further comment. 

Second, we also remind the reader that the appearing fibre products in \eqref{e:dersat2} are derived schemes. Concretely, 
this just means we also remember the Tor, so that the self-intersection $$\on{pt} \underset{ \cfg} \times 
\on{pt}$$is the spectrum of the exterior algebra $\Sym \cfg^*[1]$, by a standard 
calculation using Koszul resolutions.

We can therefore naively hope for an equivalence of monoidal categories between 
\eqref{e:dersat1} and \eqref{e:dersat2}. Instead, the derived Satake 
equivalence, which is an important theorem of Bezrukavnikov--Finkelberg \cite{bezfinksat}, 
identifies \eqref{e:dersat1} with an enlargement of \eqref{e:dersat2}. More 
precisely, they proved an equivalence
\begin{equation} \label{e:dersat} \Dmod(L^+G \backslash LG / L^+G) \simeq  
\ICn( \on{pt} / \cG  \underset{ 
\cfg / \cG} \times \on{pt} / \cG), \end{equation}
where the right hand side is the monoidal category of ind-coherent sheaves with 
{nilpotent singular support}.

\subsubsection{} Let us briefly indicate what is meant by an ind-coherent sheaf, and what it means to impose a singular support condition. These are
admittedly somewhat technical corrections to the idea of a quasicoherent sheaf, and may be safely skipped by a reader meeting 
this circle of ideas for the first time. We begin with ind-coherent sheaves, 
and then turn to the singular support 
condition. 

In the simplest possible terms, the theory of ind-coherent sheaves 
stems from the difference between finitely generated modules and finitely 
generated projective modules, i.e., vector bundles. In the abelian category, say 
for an affine scheme $X$, there are certainly fewer vector bundles than finitely 
generated modules, unless $X$ is a 
disjoint union of reduced points.

However,  in the derived category, the distinction is subtler. Indeed, if $X$ 
is smooth, then any  finitely generated module admits a finite resolution by 
vector bundles. For this reason, one cannot distinguish between finite 
complexes of finitely generated modules and finite complexes of vector bundles. 
If $X$ is not 
smooth, they again differ in the derived category, as a general finitely generated module now needs an infinite projective resolution. Explicitly, one could take the skyscraper sheaf at any singular point of $X$. In particular, the distinction between these types of complexes gives a measure of 
the singularities of $X$.

\subsubsection{}Let us treat the preceding paragraphs more carefully, following Krause \cite{krause} and Gaitsgory \cite{gaitsgoryindcoh}. To begin 
with, recall that, for an algebraic variety $X$, its category of 
quasicoherent sheaves $\QCoh(X)$ is compactly generated by the perfect 
complexes, i.e., objects locally equivalent to finite complexes of vector 
bundles. If $X$ 
is affine, this reduces to the tautology that $O_X\mod$ is generated by the ring of 
regular functions $O_X$. For non-affine $X$, the claimed generation is 
straightforward for quasi-projective varieties, but is in general a somewhat delicate 
theorem of Thomason \cite{thomasontrobaugh}.

As the global analogue of finite complexes of finitely generated modules, 
consider the bounded derived 
category of coherent sheaves$$\on{Coh}(X)^b \subset \QCoh(X).$$ By 
definition, 
$\on{IndCoh}(X)$ is its ind-completion, i.e., $\on{IndCoh}(X)$ is compactly generated by 
$\on{Coh}(X)^b$. Plainly, objects of $\on{IndCoh}(X)$ are formal 
colimits of bounded complexes with coherent cohomology, and complexes of  
homomorphisms are computed as the homotopy limit of complexes
\[
\Hom_{\on{IndCoh}(X)}( \varinjlim_i \sE_i, \varinjlim_j \sF_j) \simeq 
\varprojlim_i 
\varinjlim_j \Hom_{\on{QCoh}(X)}( \sE_i, \sF_j).
\]
For any coherent sheaf $\sE$, this formula forces $\Hom_{\on{IndCoh}(X)}(\sE, -)$ to commute with colimits, which fails in $\QCoh(X)$ unless $\sE$ is perfect. 

Since perfect complexes are bounded coherent, one obtains upon ind-completing their inclusion functor an adjunction 
\[
\on{QCoh}(X) \rightleftarrows \on{IndCoh}(X), 
\]
where the left adjoint is fully faithful and an equivalence if $X$ is smooth. 
In particular, one thinks of $\on{IndCoh}(X)$ as an enlargement of 
$\on{QCoh}(X)$, and one has a truncation of any ind-coherent sheaf to a 
quasi-coherent sheaf. Morever, $\on{IndCoh}(X)$ carries a natural 
$t$-structure and agrees with $\on{QCoh}(X)$ on bounded below objects. 
Informally, one obtains $\on{IndCoh}(X)$ by adjoining objects concentrated in 
cohomological degree $-\infty$. 

\begin{re} It may be orienting to recall that we already met such renormalizations of derived categories of infinite cohomological dimension when discussing Kac--Moody representations and D-modules on the loop group, cf. Example \ref{ex:kmmods}. 
\end{re}

\begin{re} The idea of systematically distinguishing between perfect and bounded coherent complexes on singular varieties also arose earlier in algebraic geometry, e.g., the quotient $\on{IndCoh}(X)/\QCoh(X)$ is the ind-completion of Orlov's singularity category \cite{orlovsingularity}. 
\end{re}

\subsubsection{} Let us next explain the idea of singular support, as developed by  Gulliksen \cite{gulliksen}, Eisenbud \cite{eisenbud}, Avramov, and Buchweitz \cite{avramovbuchweitz} for complete intersections, and extended by Arinkin--Gaitsgory to derived complete intersections \cite{arinkingaitssingsupp}; see also the important \cite{bensoniyengarkrause}. By definition, the difference between ind-coherent and quasicoherent sheaves stems 
from singularities of $X$, and, in particular, the quotient $$\IC(X)/\QCoh(X)$$ 
will be supported along the singular locus of $X$. Given a coherent sheaf 
$\sF$, one might like an invariant of it that records, to first approximation, 
where along the singularities of $X$ it fails to be perfect. Singular 
support provides a refined answer to this question.

Here is a slightly more precise formulation. Let us suppose $X$ is a 
quasi-smooth derived scheme, i.e., one whose cotangent complex $T^*(X)$ is 
supported in 
cohomological degrees -$1$ and $0$. Then, given an object $\sF$ as above, its 
singular support is a closed conical subset $\Lambda$ of $H^{-1} T^*(X)$. The singular support measures, affine by affine, the local nilpotence 
of elements of $\OO_X$, recovering the usual support, as well as certain natural degree two self-maps of $\sF$. Moreover, given any such subset $\Lambda$, one obtains the full subcategory 
\[
   \IC_\Lambda(X) \subset \IC(X)
\]
of ind-coherent sheaves with singular support contained in $\Lambda$. For 
example, if $\Lambda$ is the zero section, this recovers $\QCoh(X)$. 

Although we will not provide a general definition of how to compute singular 
support, we will discuss some special cases momentarily. First, however, we include a few orienting remarks. 

\begin{re}\label{r:singsupp}A derived scheme is 
quasi-smooth if and 
only if it locally may be written as a fibre product $X \times_Z Y$, where 
$X,Y$, and $Z$ are smooth varieties. In this local picture, the cotangent 
complex takes the following explicit form. Namely, for a point 
$(x,y)$ 
of the fibre product with common image $z$, one has a distinguished triangle of tangent complexes
\begin{equation} \label{e:dts}
    T_{(x,y)} (X \times_Z Y) \rightarrow T_x(X) \oplus T_y(Y) \rightarrow 
    T_z(Z) \xrightarrow{+1}.
\end{equation}
 For the reader unfamiliar with this distinguished triangle, it may be helpful to first understand the case when $X$ and $Y$ map transversally to $Z$. Here, the fibre product is again smooth, and the above distinguished triangle becomes a hopefully intuitive short exact sequence. For the general case, the point is that the derived structure on the fibre product modifies its tangent complex so that the most naive modification of the transverse situation holds.

Upon dualizing \eqref{e:dts}, one sees that $H^{-1}$ of the cotangent bundle is 
concretely the 
codirections $\xi \in T^*_z(Z)$ normal to the images of $T_x(X)$ and $T_y(Y)$ 
or, even more informally, the remaining directions to walk in $Z$.

An important special case of the preceding is a fibre product of the form \begin{equation} 
\label{e:vanlocus} \{0\} 
\times_{\mathbb{A}^n} Y,\end{equation}i.e., the derived vanishing locus of $n$ 
functions $f_1, \ldots, f_n$ on $Y$.  In 
this way, quasi-smoothness is an example of the hidden smoothness philosophy in 
derived algebraic geometry. Namely, while the underlying classical scheme of 
the fibre product may be arbitrarily singular, it underlies a derived scheme 
that behaves like a complete intersection, i.e., is only mildly singular.
 
\end{re}

Let us next explain what singular support amounts to in a few examples.

\begin{ex}\label{ex:singsupp1}For a coherent sheaf $\sF$, the intersection of its singular support 
with the zero section is simply the usual support. 
\end{ex}

\begin{ex}\label{ex:singsupp2} Suppose $X$ may be exhibited as a hypersurface in a smooth variety 
$U$. In this case, $$H^{-1} T^*(X)$$consists of $X$, 
thought of as the zero section, along with a line $\ell_x$ placed at 
each singular point $x$ of $X$. Suppose a coherent sheaf $\sF$ contains such 
an $x$ in 
its usual support. Then its singular support contains $\ell_x$ if and only if 
$\sF$ 
is not isomorphic to a perfect complex in any Zariski neighborhood of $x$.
\end{ex}

\begin{ex}Suppose $X$ is given, as in the spectral side of derived Satake, as 
the self-intersection 
of the origin in a vector space$$X = 0 \underset V \times 0.$$
In this case, one has that $T^*_{(0,0)} X \simeq V^*[1]$. Here, 
singular supports may be understood as usual supports on the other side of 
Koszul duality. Namely, the skyscraper 
sheaf $k$ generates $\IC(X)$, which yields upon taking its endomorphisms an 
equivalence 
\[
    \IC(X) \simeq \Sym V[-2] \mod. 
\]
Via this identity, for any coherent sheaf $\sF$, its singular support agrees 
with the usual support of the corresponding 
$\Sym 
V[-2]$-module, ignoring grading shifts. One can replace $0$ and $V$ here by any smooth point of any variety. \end{ex}

\subsubsection{} Having discussed the basics of ind-coherent sheaves and 
singular support, let us return to the derived Satake equivalence, i.e.,

\begin{equation} \label{e:der} \Dmod(L^+G \backslash LG / L^+G) \simeq  
\ICn( \on{pt} / \cG  \underset{ 
	\cfg / \cG} \times \on{pt} / \cG). \end{equation}
In this equivalence, the subscript `nilp' refers to the subcategory of 
ind-coherent sheaves 
with singular support lying in the nilpotent cone $$\sN \subset \cfg^* \simeq 
H^{-1}T^*( \on{pt}/\cG \underset {\cfg / \cG}   \times \on{pt}/\cG).$$

\subsubsection{} As  we will see shortly, the pattern encountered in the 
unramified case, where one has to enlarge the spectral side by allowing 
ind-coherent sheaves with nilpotent singular support, occurs in several other 
cases of matching intertwining operators. This necessitates tweaking the 
spectral side of local geometric Langlands, which is accomplished by passing from quasicoherent 
sheaves of categories to ind-coherent sheaves of categories with nilpotent 
singular support
\[
\TwoCoh(\LS) \rightsquigarrow \TwoIC(\LS).
\]
We will touch on what this means momentarily, but for now let us say that with 
this correction, we obtain the current formulation of the local geometric 
Langlands conjecture, namely an equivalence
\begin{equation}
\label{e:loclang}\Dmod(LG)\mod \simeq \TwoIC(\LS). 
\end{equation}

\subsubsection{} Let us briefly sketch the theory of ind-coherent sheaves of 
categories with singular support, which is due to Arinkin \cite{arinkinindcohcats}. 

When one increases the categorical level from ind-coherent sheaves to 
ind-coherent sheaves of categories, to set up a theory of singular support one no longer asks that the base variety $X$ be quasi-smooth, i.e., has
its cotangent complex in degrees -1 and 0, but instead that it be smooth, i.e., 
has its cotangent complex in degree 0. We will see why momentarily. 

So, let $X$ be a smooth algebraic variety. Using 1-affineness, one can show 
that, for any proper map $f: Z \rightarrow X$ from a smooth variety, $\QCoh(Z)$ 
compactly generates the full subcategory of $\TwoCoh(X)$ consisting of objects supported over 
$\on{im}(f)$. In particular, this subcategory identifies with modules for the 
convolution algebra 
\[
  \QCoh(Z \underset X \times Z).
\]

The $(\infty,2)$-category of ind-coherent sheaves of categories on $X$ has the 
following basic properties. First, it carries a monoidal action of 
$\TwoCoh(X)$, so that one can speak of the supports of objects.  Moreover, it is
generated by objects of the form $\QCoh(Z)$, for $Z$ as above. However, the 
endomorphisms of such an object are now increased to
\[
   \IC( Z \underset{X}{\times} Z), 
\] 
which, roughly speaking, grows larger relative to its quasicoherent analog as $Z \rightarrow X$ becomes further from a 
smooth map. More precisely, as in Remark \ref{r:singsupp}, the singular support 
of objects in this convolution category lie along codirections in $X$ 
normal to the tangent bundle of $Z \times Z$.  As a consequence, the less 
smooth $f$ is over a given point $x$, the more ind-coherent sheaves of 
categories supported over $x$ lie in the subcategory of $\on{2-IndCoh}(X)$ 
generated by $\QCoh(Z)$.

In particular, $\QCoh(Z)$ will no longer generate a full subcategory of 
$\on{2-IndCoh}(X)$ corresponding to objects supported over $\on{im}(f)$, but 
rather 
the following microlocal refinement.\footnote{Explicitly, to see that a refinement is necessary, 
note that, in its absence, by taking $Z = X$ we would equate $\on{2-IndCoh}(X)$ and $\TwoCoh(X)$.}  Since we 
assumed $f$ is proper, we may form the conormal variety $N^\vee_{Z}$, i.e., the 
closed isotropic subvariety of $T^*X$ obtained by projecting the kernel of $f^* 
T^*X \rightarrow T^*Z$ along $$Z \underset X \times T^*X \rightarrow T^*X.$$ 
One can associate to any ind-coherent sheaf of categories its support in 
$T^*X$, and $\QCoh(Z)$ generates the full subcategory of objects with 
microsupport lying in $N^\vee_Z$.

More generally, given any closed, conical isotropic subvariety $\Lambda$ of 
$T^*X$, one associates the full subcategory
$$
\on{2-IndCoh}_\Lambda(X) \subset \on{2-IndCoh}(X).
$$
This is again generated by the (truncations) of $\QCoh(Z)$, with endomorphisms 
$$\on{IndCoh}_\Lambda(Z \underset X \times Z),$$where now $\Lambda$ is used to 
impose a usual singular support restriction on usual ind-coherent sheaves as in 
Remark \ref{r:singsupp}. In particular, taking $\Lambda$ to be the zero 
section, we recover $\TwoCoh(X)$.

\begin{re} An alternative approach to the formalism of ind-coherent sheaves and ind-coherent sheaves of categories, including a theory of singular supports beyond the quasi-smooth and smooth cases, respectively,   is under development by di Fiore--Stefanich. \end{re}

\subsubsection{} The previous discussion does not literally apply to $\LS$. This is because, due 
roughly to its infinite dimensionality, it is not a smooth Artin stack. However, it is 
expected that it carries a similar theory of ind-coherent sheaves of 
categories.

Admitting the existence of such a formalism, it only remains from \eqref{e:loclang} to explain what the 
singular 
support condition `nilp' refers to. However, this should  not be surprising 
given what we met in the unramified case. Namely, the cotangent complex to $\LS$ 
at a 
connection $(\check{\sP}, \nabla)$ is given by the de Rham complex computing the flat sections on the punctured disk of the coadjoint bundle 
\begin{equation} \label{e:datboi}
   \check{\sP} \overset{\cG} \times \cfg^*. 
\end{equation}
In particular, the cotangent space, i.e., its zeroth cohomology, canonically identifies with the flat sections of 
\eqref{e:datboi}. With this, `nilp' cuts out within each cotangent space the flat 
sections which are nilpotent.

\subsubsection{} Having dealt with the necessary corrections, let us at last 
have the promised fun and inspect the categories of intertwining operators 
between the 
objects in Figure \ref{f:table}. In this way, we will meet many beautiful theorems and conjectures, which will hopefully give the reader more of a hands-on feel for this subject. 

On either side, all the appearing categories are self-dual. For the categories of D-modules and quasicoherent sheaves, up to manageable issues of an infinite dimensional nature, this follows from Example \ref{ex:dmodselfdual}. We will discuss the case of Kac--Moody representations in more detail in Section \ref{ss:affhc} below. Admitting this self-duality, via the 
canonical equivalences
\[
\Hom_{\Dmod(LG)\mod}(\sC, \sD) \simeq \Hom_{\Dmod(LG)\mod}(\sD^\vee, \sC^\vee),
\]
and similarly for categories over $\LS$, we need only write down the categories 
of intertwining operators `above the diagonal.'

We record these intertwining operators in Figures \ref{f:tableaut} and \ref{f:tablespec}. The calculations of these categories on the automorphic and spectral sides 
essentially follow the pattern of those performed in Section 
\ref{s:catrepgrp}, particularly Example \ref{e:fibprodvar} and the discussion of Section \ref{ss:invariantscosets}. Let us now compare the entries in the two tables and 
comment on the arising equivalences, both known and unknown. 

\footnotesize
	
\begin{center}
	\begin{figure}
		
		\setlength{\tabcolsep}{.5em}
		\begin{tabular}{|c |c| c| c| c|}
			
			\hline & &&& \\ & $\Dmod(\Gr_G)$ & $\Dmod(\Fl_G)$ & $\Dmod(LG/LN, \psi)$& 
			$\fgc\mod$  \\ & &&& \\ \hline   & &&& \\$\Dmod(\Gr_G)$  \hspace{0mm}& 
			\hspace{0mm} 
			$\Dmod(L^+G 
			\backslash LG / L^+G)$ \hspace{0mm}& \hspace{0mm} $\Dmod(L^+G 
			\backslash LG / I)$ \hspace{0mm} &  \hspace{0mm}$\Dmod(L^+G 
			\backslash LG / LN, \psi)$ \hspace{0mm}&  
			\hspace{0mm}$\fgc\mod^{L^+G}$  \hspace{0mm}\\  & &&&

			\\ \hline   & &&& \\$\Dmod(\Fl_G)$  \hspace{0mm} &  \hspace{0mm} 
			\hspace{0mm} &  \hspace{0mm}$\Dmod(I \backslash LG / I)$ \hspace{0mm} 
			&  \hspace{0mm} 
			$\Dmod(I \backslash LG / LN, \psi)$ \hspace{0mm} &  \hspace{0mm}  
			$\fgc\mod^I$  \hspace{0mm}\\  
			& &&&

			\\ \hline   & &&& \\ \hspace{0mm}$\Dmod(LG/LN, \psi)$  \hspace{0mm}&& 
			\hspace{0mm} 
			&  
			\hspace{0mm}$\Dmod(LN, \psi \backslash 
			LG / LN, \psi)$ \hspace{0mm} &  \hspace{0mm}  $\fgc\mod^{LN, \psi}$  
			\hspace{0mm}\\  & &&&

			\\ \hline   & &&& \\$\fgc\mod$ &&&&  \hspace{0mm}$\fgc\oplus 
			\fgc\mod^{LG}$  \hspace{0mm}\\  & &&& 
			\\ 
			\hline 		\end{tabular}
			    		\caption{Some intertwining operators on the automorphic side of local 
		geometric Langlands.}
		

		\label{f:tableaut}

	\end{figure}
\end{center}

\begin{center}
	\begin{figure}
		

		\setlength{\tabcolsep}{.5em}
		\begin{tabular}{|c |c| c| c| c|}
			
			\hline & &&& \\ & $\QCoh(\on{pt}/\cG)$ & $\QCoh(\widetilde{\sN} / \cG)$ 
			& $\QCoh(\LS)$& 
			$\QCoh(\Op)$  \\ & &&& \\ \hline   & &&& \\$\QCoh(\on{pt} / \cG)$  
			\hspace{0mm}& 
			\hspace{0mm} 
			$\ICn( \pt / \cG  \underset {\cfg/\cG} \times  \pt / \cG)$ 
			\hspace{0mm}& 
			\hspace{0mm} 
			$\ICn(\pt / \cG \underset { \cfg / \cG} \times \widetilde{\sN} / \cG 
			)$ 
			 \hspace{0mm} &  \hspace{0mm}$\QCoh( \pt/ \cG)$ \hspace{0mm}&  
			\hspace{0mm}$\QCoh( \Op^{\on{unramified}})$  \hspace{0mm}\\  & &&&

			\\ \hline   & &&& \\$\QCoh(\widetilde{\sN} / \cG)$  \hspace{0mm} &  
			\hspace{0mm} 
			\hspace{0mm} &  \hspace{0mm}$\ICn( \widetilde{\sN}/ \cG \underset{ 
			\cfg / \cG } \times \widetilde{{\sN}} / \cG)$ 
			\hspace{0mm} 
			&  \hspace{0mm} 
			$\QCoh( \widetilde{\sN}/ \cG)$ \hspace{0mm} &  \hspace{0mm}  
			$\QCoh(\widetilde{\sN}/\cG \underset \LS \times \Op)$  
			\hspace{0mm}\\  
			& &&&

			\\ \hline   & &&& \\ \hspace{0mm}$\QCoh(\LS)$  \hspace{0mm}&& 
			\hspace{0mm} 
			&  
			\hspace{0mm} $\QCoh(\LS)$ \hspace{0mm} &  \hspace{0mm}  $\QCoh(\Op)$
			\hspace{0mm}\\  & &&&

			\\ \hline   & &&& \\$\QCoh(\Op)$ &&&&  \hspace{0mm}$\QCoh( \Op 
			\underset \LS \times \Op)$ \hspace{0mm}\\  & &&& 
			\\ 
			\hline 		\end{tabular}
		\caption{Some intertwining operators on the spectral side of local 
			geometric Langlands.}
		
		\label{f:tablespec}
		

	\end{figure}
\end{center}

\normalsize

\subsubsection{}  Let us begin with the unramified case, i.e., the first row of 
the tables. The first equivalence 
\begin{equation}
\Dmod(L^+G \backslash LG / L^+G) \simeq \ICn( \pt / \cG \underset {\cfg/\cG} 
\times \pt / \cG)
\end{equation}
is one we have already met, namely the  derived Satake equivalence of 
Bezrukavnikov--Finkelberg \cite{bezfinksat}. 

The derived Satake equivalence is $t$-exact and on abelian 
categories recovers the geometric Satake equivalence
\[
 \Dmod(L^+G \backslash LG / L^+G)^\heartsuit \simeq \QCoh(\pt / \cG) = 
 \on{Rep}(\cG),  
\]
which is an earlier fundamental theorem due to Lusztig \cite{lusztigqchar}, Drinfeld, Ginzburg \cite{ginzburgsatake}, Mirkovi\'c--Vilonen \cite{mirkovic-vilonen}. Both equivalences, but slightly more transparently the abelian one, categorify the analogous assertion in the $p$-adic setting 
\[
   \on{Dist}( G(\Fq[[t]]) \backslash G(\Fq(\!(t)\!)) / G(\Fq[[t]])) \simeq K_0( 
   \on{Rep}(\cG)) \underset{\mathbb{Z}}{\otimes} k,
\]
which, along with its mixed characteristic counterpart, is an isomorphism due 
to Satake \cite{satake63}.

Let us make two comments about this.

\begin{re} \label{r:paramcount} First, for readers meeting this for the first time, to get some 
purchase on this equivalence it is 
instructive to match parameters, i.e., bases for the Grothendieck groups of 
(almost) compact objects. 
	
 The basic point here is that there are 
	only countably many $L^+G \times L^+G$ orbits on $LG$, and these are 
	naturally indexed by irreducible representations of $\cG$. More carefully, 
	any 
	cocharacter
	\[
	   \check{\mu}: \mathbb{G}_m \rightarrow T	\]
	yields, upon restriction to the formal disk about the origin, a point of 
	$LT$, which we denote by $t^{\check{\mu}}$. If we let 
	$\check{\Lambda}^+$ denote the dominant cocharacters, $LG$ is then stratified by 
	the orbits
	\[
   L^+G \cdot t^{\check{\mu}} \cdot L^+G, \quad \text{for } \check{\mu} \in 
   \check{\Lambda}^+. 
\] 
But these are the highest weights of $\cG$-modules, and the equivalence 
exchanges 
the corresponding intersection cohomology complexes on $\Gr_G$ and simple 
$\cG$-modules. 	\end{re}

\begin{re} Second, here is an orienting comment. The underlying space on the 
automorphic side of derived Satake 
parametrizes triples $(\sP_1, \sP_2, \tau)$, where $\sP_1$ and $\sP_2$ are 
$G$-bundles on the disk and $\tau$ is an identification of them on the 
punctured disk.
	
	Similarly, the underlying space on the spectral side parametrizes triples 
	$(\sE_1, \sE_2, \tau)$, where $\sE_1$ and $\sE_2$ are $\cG$-connections on 
	the formal disk and $\tau$ is an isomorphism of their restrictions 
	to the punctured disk. Note that since every $\cG$-connection on the formal disk is trivializable, the underlying classical stack is $\pt / \cG$. However, this moduli problem is asking us to intersect $\pt / \cG$ with itself in $\LS$, which affords the derived structure crucial to making the equivalence hold. 
	 	
	 Thus, somewhat informally speaking, the equivalence exchanges 
	 constructible objects on a coherent moduli space for $G$ with coherent 
	 objects on an analogous constructible moduli space for $\cG$. 	
	 	
\end{re}

\subsubsection{} \label{s:abg}The next equivalence of the table, namely
\[
   \Dmod(L^+G \backslash LG / I) \simeq \ICn(\pt/\cG \underset {\cfg / \cG} 
   \times \widetilde{\sN} / \cG)
\]
is also a very nice theorem due to Arkhipov--Bezrukavnikov--Ginzburg \cite{abg}. As in 
Remark 
\ref{r:paramcount}, we may count parameters as follows. On the automorphic 
side, the relevant stratification of $LG$ is now by the strata
\[
     L^+ G \cdot t^{\check{\mu}} \cdot I, \quad \text{for } \check{\mu} \in 
     \check{\Lambda}.
\]
That is, to shrink from having $L^+G$ to $I$ on the right, we must grow by a factor of the 
finite Weyl group, i.e., from dominant cocharacters to all cocharacters. 

On the spectral side, recalling that $\widetilde{\sN}$ is the cotangent bundle 
to the (finite) flag variety, we are considering $\cG$-equivariant coherent 
sheaves on a derived thickening of the zero section, i.e., of the flag variety 
itself. Here, up to ignoring the thickening, the same lattice $\check{\Lambda}$ occurs in 
\[
   \QCoh( \Fl_{\cG}  / \cG  ) \simeq \QCoh(  \cB \backslash \cG / \cG) \simeq 
   \QCoh( \cB \backslash \pt) \simeq \on{Rep}(\cB), 
\]
as  equivariant line bundles on  $\Fl_{\cG}$ or, equivalently, the 
characters of $\on{Rep}(\cB)$.

\subsubsection{}We next pair the spherical and Whittaker categories to obtain 
\[
   \Dmod(L^+G \backslash LG / LN, \psi) \simeq \QCoh( \pt / \cG). 
\]
This is also a known and important theorem due to Frenkel--Gaitsgory--Vilonen \cite{fgv}. It may be understood as a geometric refinement of the 
Casselman--Shalika formula from $p$-adic representation theory \cite{cassshalika}, which 
explicitly 
determines the Whittaker covector of an unramified representation. 

The parameter count here is similar. Namely, $LG$ is stratified by the double 
cosets
\[
   L^+G \cdot t^{\check{\mu}} \cdot LN, \quad \text{for } \check{\mu} \in 
   \check{\Lambda}, 
\]
but an orbit supports Whittaker sheaves if and only if the coweight 
$\check{\mu}$ is antidominant.

\begin{re}Let us also briefly comment on a basic difference with the 
bi-spherical situation, i.e., derived Satake. There, the spectral side is a 
derived enhancement of $\on{Rep}(\cG)$. Roughly speaking, this occurs since 
the stalks of the intersection cohomology sheaves on $$L^+G \backslash LG / 
L^+G$$are quite interesting along orbit closures and correspond to $q$-weight 
multiplicities for $\cG$-representations \cite{lusztigqchar}. 

In the present case, one has no such derived enhancement. This occurs because 
the intersection cohomology objects on $$L^+G \backslash LG / LN, \psi$$are 
clean, 
i.e., have vanishing $!$-stalks and $*$-stalks along the boundary of the open 
orbit in their 
support.

\end{re}

\subsubsection{} \label{ss:sphericalkmmod}Finally, let us pair spherical vectors and Kac--Moody 
representations to obtain 
\[
   \fgc\mod^{L^+G} \simeq \QCoh(\Op^{\on{unramified}}), 
\]
where by definition the latter moduli space parametrizes oper structures on the trivial bundle, i.e. 
\[
\Op^{\on{unramified}} \simeq \pt / \cG \underset \LS \times \Op. 
\]

This is also a known theorem, due to Frenkel--Gaitsgory \cite{fg-sph}, which contains and 
strengthens many prior results about $L^+G$-equivariant Kac--Moody 
representations at critical level due to Feigin--Frenkel, Beilinson--Drinfeld, 
and others. 

In this case, the Grothendieck groups are no longer countable. Instead, one has 
continuous families of simple modules parametrized by pro-finite 
dimensional affine spaces. Nonetheless, we can give the following rough 
count. 

On the automorphic side, one has a collection of compact generators, the Weyl 
modules, which are by definition the parabolic inductions of the irreducible 
finite dimensional $G$-modules$$ \mathbb{V}_\lambda := \on{pind}_{\fg}^{\fgc}( 
L_\lambda), \quad \text{for } \lambda \in \Lambda^+.$$
Parabolic induction in the present setting of Lie algebras is given by the correspondence $$\fg \leftarrow \fg[[t]] \rightarrow \fgc,$$i.e. is the composition of restriction and induction functors
\[
\on{pind}_{\fg}^{\fgc}: \fg\mod \xrightarrow{\on{Res}} \fg[[t]]\mod \xrightarrow{\on{ind}} \fgc\mod. 
\]

It is a nontrivial fact that each $\mathbb{V}_\lambda$ has endomorphisms given by a quotient of the center $Z(\fgc)$, and is moreover projective in $\gkc\mod^{L^+G}_\lambda$, where $\gkc\mod_\lambda$ denotes the category of $\gkc$-modules on which $Z(\fgc)$ acts through the same quotient.

\begin{ex}\label{e:regop} If $\lambda = 0$, then we obtain the vacuum 
module\footnote{The 
terminology `vacuum' here stems from its origins in conformal field theory. 
Namely, this is the underlying vector space of a vertex operator algebra, which 
occurs as local operators in chiral conformal field theories with 
$G$-symmetry.} at 
critical level 
	\[
	  \mathbb{V}_0 \simeq \on{pind}_{\fg}^{\fgc}(k). 
	\]
Its endomorphisms canonically identify, under the Feigin--Frenkel isomorphism, 
with the algebra of functions on opers on the non-punctured formal disk
\[
    \Op^{\on{reg}} \simeq \{ d + f dt + \check{\fb}[[t]]dt \} / \cN[[t]].
\]
\end{ex}

On the spectral side, as we vary $\lambda$, these endomorphism algebras correspond to disjoint closed 
subschemes inside $\Op$. The moduli space $\Op^{\on{unramified}}$ is a formal thickening of their union, and the theorem further matches these nilpotent normal directions with the self-extensions of Weyl modules.

\subsubsection{}\label{s:2realz} We are ready to move on to the next row of Figures \ref{f:tableaut} and \ref{f:tablespec}, i.e., 
the case of tame ramification and unipotent monodromy. We first meet the 
endomorphisms, i.e., a monoidal equivalence
\begin{equation} \label{e:2realz}
 \Dmod(I \backslash LG / I) \simeq \ICn( \widetilde{\sN} / \cG \underset {\cfg 
 / \cG} \times \widetilde{\sN} / \cG). 
\end{equation}
This attractive statement is remarkable for many reasons. For one, there is not 
such a coherent realization of the finite Hecke category, and, in this sense, the equivalence 
really is something particular to the affine setting and, specifically, its 
connection with arithmetic.\footnote{Admittedly, 
the construction via Soergel bimodules can be understood as writing the 
monodromic Hecke 
category as the homotopy category of certain coherent sheaves on $$\check{T} 
\underset 
{\check{T}/W_f}{\times} \check{T},$$cf. \cite{bezrichetop}. However, this identification is somewhat subtler, as it is 
only fully faithful on tilting objects non-derivedly, i.e., the corresponding 
Soergel bimodules have nontrivial Exts in the world of bimodules. In any case, 
another pleasant feature 
of \eqref{e:2realz} is that the Soergel bimodule presentation of the finite 
Hecke category is indeed embedded into it.} As with other statements of low ramification we 
have met so far, it is a celebrated theorem due to 
Bezrukavnikov \cite{bez2realz}.

As far as 
we know, this equivalence was first formulated as a conjecture by Ginzburg \cite{chrissginzburg}. A major 
motivation comes from the $p$-adic 
setting, where an analogous isomorphism between $\on{Dist}(I \backslash G(\Fqt) 
/ 
I)$, i.e., the affine Hecke algebra, and an appropriate specialization of the 
equivariant $K$-theory of the Steinberg variety $$\widetilde{\sN} \underset 
\cfg 
\times \widetilde{\sN}$$is due 
to Kazhdan--Lusztig \cite{kazhdanlusztig87}. They, and independently Ginzburg, used this isomorphism to classify the irreducible representations of affine Hecke algebra and determine their characters. Equivalently, they determined the characters of the simple subquotients of the unramified principal series of $G(\Fqt)$. 

Let us count parameters. On the automorphic side, the relevant stratification 
of $LG$ has constituents 
\[
    I\cdot w \cdot t^{\check{\mu}} \cdot I, \quad \text{for } w \in W_f \text{ 
    and } 
    \check{\mu} \in \check{\Lambda}.  
\]
This is an affine analogue of the Bruhat decomposition of $G$ into the double cosets $B \cdot w \cdot B$, for 
$w \in W_f$, and indeed the indexing set here is the (extended) affine Weyl 
group $W_f \ltimes \check{\Lambda}.$ Note that, as in Section \ref{s:abg}, we once again shrink from $L^+G$ 
to $I$, and therefore our stratification grows by a factor of $W_f$. 

On the spectral side, these parameters arise as follows. To begin, note that 
 the Steinberg variety\footnote{Although this is the standard 
terminology, we emphasize that for the monoidal equivalence we really are 
working instead with the derived enhancement.} explicitly parametrizes triples $$(X, \cB', 
\cB''),$$where $X$ is a nilpotent element of $\cfg$ 
and $\cB'$, $\cB''$ are two Borel subgroups containing it. The Steinberg is the union of 
$W_f$ many locally closed Lagrangians in $\widetilde{\sN} \times \widetilde{\sN}$, where each 
Lagrangian is 
specified by the relative position of $\cB_1$ and $\cB_2$, i.e., a point of $\cB 
\backslash \cG / \cB$. 

\begin{re}  Recall that we are really interested in the quotient $$\widetilde{\sN} / 
\cG \underset {\cfg/\cG} 
\times \widetilde{\sN} / \cG.$$It is clarifying to recognize this as the 
cotangent (derived stack) of 
the 
Hecke stack
\begin{equation} \label{e:strestopesto} \cB \backslash \cG / \cB \simeq \Fl_{\cG} \times \Fl_{\cG} / \cG
\end{equation}
in its presentation via Hamiltonian reduction from the right hand side of \eqref{e:strestopesto}. 
Relatedly, the $W_f$ components in question are simply the conormal bundles to 
the 
corresponding Schubert cells. 
\end{re}

As in Section \ref{s:abg}, each component admits 
$\check{\Lambda}$-many equivariant line bundles, given by pullback along either 
projection to $\widetilde{\sN}/\cG$,\footnote{The two maps from 
$\check{\Lambda}$ to its Picard group differ by the action of the corresponding 
element of $w$.} so, summing over components, we get the desired count of $W_f 
\ltimes \check{\Lambda}.$

\subsubsection{}Let us now pair the unramified principal series and Whittaker 
model to obtain an equivalence
\begin{equation} \label{e:arkbez}
    \Dmod( I \backslash LG / LN, \psi) \simeq \QCoh( \widetilde{\sN} / \cG).
\end{equation}
This assertion is again known, and is a very nice theorem of  Arkhipov--Bezukravnikov \cite{ab}. Their argument is widely regarded as one of the gems in this subject. The result plays a basic role in Bezrukavnikov's proof of \eqref{e:2realz}, in parallel to the role of an analogous identity in K-theory for the 
identification of the affine Hecke algebra with the equivariant K-theory of 
Steinberg by Kazhdan--Lusztig.  

The parameter count here goes as follows. 
On the automorphic side, it is slightly more convenient to pass to the 
semi-infinite variant, i.e., 
\begin{equation} \label{e:semiinfre}  \Dmod( LN \cdot L^+T \backslash LG / LN, 
\psi).\end{equation}
In this presentation, under the stratification of $LG$ by the cosets
\[ (LN \cdot L^+T) \cdot w \cdot t^{\check{\mu}} \cdot LN,   \quad \text{for } 
w 
\in W_f \text{ and }  \check{\mu} \in \check{\Lambda}, \]
an orbit supports Whittaker sheaves if and only if $w$ equals the longest 
element $w_\circ$ of $W_f$. In particular, the set of relevant orbits is 
canonically parametrized by $\check{\Lambda}$. 

On the spectral side, we have
\[
    \widetilde{\sN} / \cG = (\cG \overset \cB \times \check{\fn}) / \cG \simeq 
    \check{\fn} / \cB,  
\]
and $\check{\Lambda}$ again corresponds to the equivariant line bundles 
pulled back from $\pt / \cB$. 

Here are two comments about this equivalence. 

\begin{re} In the $p$-adic setting, the relevant space of 
distributions\begin{equation} \label{e:antisph}\on{Dist}( I \backslash 
G(\Fqt) / N(\Fqt), \psi)\end{equation}
is an induced representation of the affine Hecke algebra. More precisely, if we 
write $H$ for the (extended) affine Hecke algebra, $H_f$ for the finite Hecke 
algebra, and $k_{\on{sgn}}$ for its sign representation, then \eqref{e:antisph} 
is the so-called antispherical representation
\[ H \underset {H_f} \otimes k_{\on{sgn}}.\]
One has a similar formula in the geometric setting, namely 
\[ \Dmod(I \backslash LG / LN, \psi) \simeq \Dmod(I \backslash LG / I) 
\underset {\Dmod(B\backslash G / B)} \otimes \Dmod(B \backslash G/ N, \psi),\]
cf. \cite{ahc}.  
\end{re}

\begin{re} We rewrote the automorphic side of 
\eqref{e:arkbez} in semi-infinite terms in \eqref{e:semiinfre}. This was to minorly simplify the combinatorics of relevant orbits. However, there are less frivolous reasons to bear this alternative expression in mind. Namely, recall that 
parabolic induction for categorical representations really does produce the 
semi-infinite version of the affine flag variety, cf. Section \ref{s:unramps}. 
Relatedly, in the proof of 
Arkhipov--Bezukravnikov, it 
is objects of semi-infinite origin, e.g., the Wakimoto sheaves, which play a 
central role. \end{re}

\subsubsection{} \label{ss:iwahorikmmod}By pairing the unramified principal series with Kac--Moody 
representations, we arrive at an equivalence
\[ 
 \fgc\mod^I \simeq \QCoh(\widetilde{\sN}/\cG \underset \LS \times \Op). 
\]
This assertion is an extension of the previously discussed results of
Frenkel--Gaitsgory in the spherical case, cf. Section \ref{ss:sphericalkmmod}. They first conjectured it in \cite{fg06loc}, and  obtained a version of it, after a further specialization of the central character, several years later as a consequence of their work on critical level localization on the affine flag variety \cite{fglocaff}.\footnote{Strictly speaking, in {\em loc. cit.} one finds an Iwahori monodromic variant.}

\subsubsection{} Let us pass to the third row of Figures \ref{f:tableaut} and \ref{f:tablespec} and study 
intertwiners out of the Whittaker model. We first meet its endomorphisms and 
an equivalence of monoidal categories
\begin{equation} \label{e:whitend} \Dmod(LN, \psi \backslash LG / LN, \psi)  
\simeq \QCoh(\LS).\end{equation}
This is a striking and essentially wide open conjecture. Even  symmetric 
monoidality for the left hand side is far from obvious. Morally, it is a 
manifestation of the uniqueness of Whittaker covectors in the current geometric 
setting under the usual correspondence between multiplicity freeness and 
commutative endomorphism algebras.

We emphasize that this equivalence would pin down an enormous swathe of 
Conjecture \eqref{e:loclang}, namely all the sheaves of categories on the 
spectral side with singular support along the zero section. Such categorical 
representations are called tempered.

\begin{re}Although \eqref{e:whitend} is very much open, here is a special case 
which is known. Namely, by work of Raskin \cite{whit}, the depth filtration on categorical 
representations of the loop group gives rise to a filtration on $$\Dmod(LN, \psi 
\backslash LG / LN, \psi),$$with each quotient a monoidal colocalization. The 
first quotient, corresponding to representations with vectors fixed by 
the prounipotent radical of an Iwahori subgroup, is the Iwahori--Whittaker 
category 
\[
     \Dmod(\mathring{I}, \psi \backslash LG / \mathring{I}, \psi). 
\]
This in turn admits a further quotient, which corresponds to passing from all 
tamely ramified representations to those with unipotent monodromy, namely the 
Steinberg--Whittaker category 
\[ 
   \Dmod(\mathring{I}, \psi, \on{s} \backslash LG / \mathring{I}, \psi, 
   \on{s}). 
\]	
Plainly, this is the full subcategory generated by averages of objects from $\Dmod(I^- \backslash LG / I^-)$, where $I^-$ is an Iwahori subgroup of $L^+G$ defined with respect to a Borel $B^-$ in general position to $B$, cf. \cite{ahc}. If we write $\sN^{\wedge}$ for the formal completion of $\cfg$ along the 
nilpotent cone $\sN$, we have a monoidal equivalence
\[ 
    \Dmod(\mathring{I}, \psi, \on{s} \backslash LG / \mathring{I}, \psi, 
    \on{s}) \simeq \QCoh( \sN^\wedge / \cG),
\]	
which is the promised piece of \eqref{e:whitend}.

 Bezrukavnikov proved the equivalence of affine Hecke categories \eqref{e:2realz} 
essentially by combining this with the previously mentioned result of Arkhipov--Bezrukavnikov. More 
precisely, let us write $\widetilde{\sN}^\wedge$ for 
the formal completion of the Grothendieck--Springer variety
$$\widetilde{\cG} := \cG \overset \cB \times \check{\fb}$$ along the Springer variety. That is, the Grothendieck--Springer variety parametrizes pairs $(X, \check{\fb}')$, where $X$ is an element of $\cfg$ and $\check{\fb}'$ is a Borel subalgebra containing it, and we complete this along the locus where $X$ is nilpotent. Then the desired equivalence \eqref{e:2realz}\footnote{Even more carefully, we mean its 
monodromic variant, from which the strict statement may be deduced formally.} is 
obtained via renormalizing the equivalence
\begin{align*} \QCoh( \widetilde{\sN}^\wedge / \cG  \underset {\cfg / \cG} 
\times 
\widetilde{\sN}^\wedge / \cG) 
& \simeq \QCoh(\widetilde{\sN}^\wedge / \cG)  \underset {\QCoh(\sN^\wedge / \cG)} 
\otimes  \QCoh(\widetilde{\sN}^\wedge / \cG ) \\ & \simeq \Dmod(\mathring{I} 
\backslash LG / \mathring{I}, \psi, \on{s}) \underset {\Dmod(\mathring{I}, 
\psi, \on{s} \backslash LG / \mathring{I}, \psi, \on{s})} \otimes 
\Dmod(\mathring{I}, \psi, \on{s} \backslash LG / \mathring{I}).  \end{align*} 

\end{re}

\subsubsection{} We next pair the Whittaker model with Kac--Moody 
representations to obtain an equivalence
\begin{equation} \label{e:affskr} \fgc\mod^{LN, \psi} \simeq \QCoh(\Op).  
\end{equation}
This is again a strong statement concerning objects of arbitrary ramification. 
Unlike other such statements we are encountering, it is in fact a theorem of 
Raskin \cite{whit}.

Namely, the affine Skryabin theorem of Raskin\footnote{We note in passing that Skryabin is an anagram of `by Raskin.'}  gives, for any level $\kappa$, a 
canonical equivalence
\[  
  \fgk\mod^{LN, \psi} \simeq \sW\mod, 
\]
where $\sW$ is the affine W-algebra. On the other hand, at the critical level Feigin--Frenkel had identified $\mathscr{W}_{\kappa_c}$ with the center $Z(\fgc)$ of the enveloping algebra \cite{feiginfrenkel91}, and therefore also with $\OO_{\Op}$. This in combination with the affine Skyrabin theorem yields \eqref{e:affskr}. 

\begin{re} The reader may wish to consult the survey of Arakawa in the present 
volume for more details about the W-algebra; see also \cite{araintro}, \cite{fbz}. Here let us only mention in 
passing that the affine Skryabin theorem enables the systematic use of 
techniques from categorical representation theory in the study of 
representations of W-algebras. As an example, one can see  \cite{locw}, which discusses localization theory for highest weight modules for W-algebras. 
\end{re}

\subsubsection{} \label{ss:affhc}In the final row, we meet the endomorphisms of Kac--Moody 
representations. These are identified with affine Harish-Chandra bimodules at 
critical 
level as follows. 

First, the category $\fgc\mod$ is canonically self-dual as a 
$\Dmod(LG)$ representation and, in particular, as an abstract dg-category. 

To see this, it is orienting to begin with a simpler case. For a finite 
dimensional Lie algebra $\mathfrak{a}$ with a central 
extension 
\[
0 \rightarrow k \cdot \mathbf{1} \rightarrow \widetilde{\mathfrak{a}}_c 
\rightarrow \mathfrak{a} \rightarrow 0,
\]
consider the category $\widetilde{\mathfrak{a}}_c\mod$ of 
$\widetilde{\mathfrak{a}}_c$-modules on which $\mathbf{1}$ acts via the 
identity. If we write $\widetilde{\mathfrak{a}}_{-c}$ for the opposite central 
extension, e.g., given by the additive inverse of  any representing 2-cocycle, one has a perfect 
pairing 
\[
\widetilde{\mathfrak{a}}_c\mod \otimes \hspace{.7mm} 
\widetilde{\mathfrak{a}}_{-c}\mod \rightarrow \mathfrak{a}\mod \rightarrow 
\Vect,
\] 	
where the first arrow tensors together representations to cancel the 
extensions, and the second is Lie 
algebra homology. 

In the affine case, it was shown by Arkhipov--Gaitsgory \cite{ag}, building on previous 
work of Frenkel--Gaitsgory \cite{fg06loc}, that for any $\kappa$ 
the composition
\[
\fgk\mod \otimes \hspace{.7mm} \widehat{\fg}_{-\kappa + 2\kappa_c}\mod 
\rightarrow \widehat{\fg}_{2 \kappa_c}\mod \rightarrow \Vect
\]
is a perfect pairing, where the first arrow tensors representations, and the 
second arrow is the functor of semi-infinite cohomology. Briefly, this is the natural 
homology theory for representations of infinite dimensional Lie algebras like 
$\gk$ (formally, Tate Lie algebras), but it requires a central extension of level $2\kappa_c$. One has such a duality statement for general Tate Lie algebras, cf. \cite{mys}, \cite{lpw}. 

In particular, one has a self-duality
\[
    \fgc\mod \otimes \hspace{.7mm} \fgc\mod \rightarrow \Vect,
\]
which is another sense in which the critical level is the midpoint of the 
space of levels. Moreover, as in Section \ref{ss:hchbimod}, this identifies the 
$\Dmod(LG)$ equivariant endofunctors of $\fgc\mod$ with the category of critical level affine 
Harish-Chandra bimodules 
\[
      \fgc \oplus \fgc \mod^{LG}.
\]
This is a beautiful, if somewhat elusive, category. For example, all of its 
objects lie in cohomological degree $-\infty$ with respect to the natural 
$t$-structure on $\fgc \oplus \fgc\mod$, i.e., it consists of homological 
phantoms introduced through renormalization. 

\begin{re}The affine Harish-Chandra bimodules at other levels 
are also subtle, with some of the basic works and conjectures 
due to I. 
Frenkel--Malikov \cite{frenkelmalikov1} and Gaitsgory \cite{quantum-langlands-summary}. 
\end{re}

For these bimodules, local geometric Langlands predicts a monoidal equivalence 
with the spectral convolution algebra 
\[ 
      \fgc \oplus \fgc \mod^{LG} \simeq \QCoh( \Op \underset \LS \times \Op). 
\]
This is a striking and essentially wide open conjecture. To our knowledge,
even the spherical case and the case of tame ramification with unipotent 
monodromy are open, though we expect they are within reach using some of the 
aforementioned work of Frenkel--Gaitsgory, cf. Sections \ref{ss:sphericalkmmod} and \ref{ss:iwahorikmmod}, combined with the localization theorem at 
critical level, which we now turn to.

\subsection{Localization theory}
	\label{ss:locaff}
	\subsubsection{}  We would like to discuss an important 
	application of ideas and techniques from local geometric Langlands to a 
	problem purely within the representation theory of affine Lie algebras. We 
	begin by setting up the context.

	 \subsubsection{} Recall that in Section \ref{s:bbloc} we discussed localization 
	 theory for 
	 reductive Lie algebras and sketched how it was used to determine the 
	 characters of the simple highest weight modules.

	 For affine Lie algebras, an analogous program was developed at noncritical levels by Kashiwara--Tanisaki in a series of works \cite{kashiwara90}, \cite{kashiwaratanisakiposlevelii}, \cite{kashiwaratanisakinegativeintegrallevel}, \cite{kashiwaratanisakinegativerationallevel}, \cite{kashiwaratanisakipositiverational}, \cite{kashiwaratanisakialllevels}. Briefly, 
	 at negative level they identified D-modules on the affine flag manifold 
	 with certain Kac--Moody representations. At positive level, they did the 
	 same, but instead with D-modules on the thick affine flag manifold, a 
	 scheme of infinite type parametrizing $G$-bundles on $\mathbb{P}^1$ with 
	 full level structure at zero and a Borel reduction at infinity. As a 
	 consequence, they obtained the characters for simple highest weight modules at noncritical levels.

	After their work, it was understood that to make inroads on the remaining case, namely the representation theory at critical level, one would need further ideas and tools. We would like to now 
	sketch what localization theory looks like here, following work of Beilinson--Drinfeld, Frenkel--Gaitsgory, and 
	Raskin--Yang. 
	
	\subsubsection{}  The basic new complication at critical level is that the geometric side of the localization theorem is not simply D-modules on a space, but a tweak thereof. Pleasantly, one can arrive at the correct statement using only what we have discussed so far about local geometric Langlands. Let us see how.

	To begin, recall that for any $\cG$ oper $\chi$, 
	with underlying local system $\sigma$, local geometric Langlands predicts 
	that the categorical loop group representation $\sC_{\sigma}$ corresponding 
	to $\Vect \simeq \QCoh(\pt)$ under the map $$\sigma: \pt \rightarrow \LS$$ 
	may be 
	obtained by renormalizing 
	\[
	  \fgk\mod \underset{\QCoh(\Op)} \otimes \Vect,
	\]
	i.e., Kac--Moody representations with central character $\chi$.

   \begin{re} Now that we have discussed ind-coherent sheaves of categories, 
   this renormalization may be understood as a consequence of the formal smoothness of 
   the map $$\Op 
   \rightarrow \LS,$$which, in particular, implies that the sub-2-category 
   generated by $\QCoh(\Op)$ consists solely of quasi-coherent sheaves of categories. 
   That is, there is no room for nontrivial singular support.  
   
   On the Langlands dual side, this predicts that $\fgc\mod$ is tempered, i.e., 
   that the canonical map 
   \[
        \Dmod(LG/ LN, \psi ) \underset {\Dmod(LN, \psi \backslash LG / LN, 
        \psi)} \otimes \mathscr{W}_{\kappa_c}\mod \rightarrow \fgc\mod
   \]
   is an equivalence. This is a remarkable conjecture about Kac--Moody 
   representations at critical level. 
   \end{re}

	\subsubsection{} Let us specialize the previous discussion, which applies 
	to any local system $\sigma$, to  the trivial connection. For this, we may 
	tautologically factor $\sigma$ as
	\[
	  \pt \rightarrow \pt/\cG \rightarrow \LS.
	\]

	We have already discussed in some detail what $\QCoh(\pt/\cG)$ is 
	meant to 
	correspond to under local geometric Langlands, namely D-modules on the 
	affine Grassmannian. This gives an alternative presentation of 
	$\sC_{\sigma}$, as follows. Since $\pt / \cG \rightarrow \LS$ is as far 
	from smooth as possible, we have that $\Vect$ lies in the category 
	generated by $\QCoh(\pt / \cG)$. That is, we find that

	\begin{align*}
	     \Vect &\simeq  \QCoh( \pt / \cG) \underset {\Hom_{\TwoIC(\LS)}( 
	     \QCoh(\pt / \cG), \QCoh(\pt / \cG))} \otimes \Hom_{\TwoIC(\LS)}( 
	     \QCoh(\pt / 
	     \cG), \QCoh(\pt)) \\ &\simeq 
	     \QCoh( \pt / \cG) \underset {\ICn(\pt / \cG \underset \LS \times \pt / 
	     \cG)} \otimes \ICn( \pt \underset \LS \times \pt ). 
	\end{align*}
Therefore, on the automorphic side, we should have that
\begin{align}
   \sC_{\sigma} &\simeq \Dmod(\Gr_G)\underset {\ICn(\pt / \cG \underset \LS 
   \times \pt / 
   	\cG)} \otimes \ICn( \pt \underset \LS \times \pt ).
\intertext{As tensoring from $\ICn( \pt / \cG \underset \LS \times \pt / \cG)$ 
to 
$\ICn(\pt 
\underset  \LS \times \pt)$ is simply forgetting equivariance for 
the 
action of $\cG$, we may rewrite this more plainly as }
 \sC_{\sigma} & \simeq \Dmod(\Gr_G) \underset{\QCoh(\pt / \cG)} \otimes 
 \QCoh(\pt).\label{e:heckeeig}
	\end{align}
	
Explicitly, this is the category of Hecke eigensheaves on the affine 
Grassmannian. 
That is, recall that $\QCoh( \pt / \cG)$, i.e., $\on{Rep}(\cG)$, is identified 
with the abelian Satake category $$\Dmod(L^+G \backslash LG / 
L^+G)^\heartsuit.$$This category acts via right 
convolution on D-modules on the affine Grassmannian, which we denote by 
\[
    - \overset {L^+G} \star -: \Dmod(LG/L^+G) \otimes \Dmod(L^+G \backslash LG 
    / L^+G) \rightarrow \Dmod(LG / L^+G).
\]
Our desired category consists of eigenobjects for this action, i.e., sheaves 
$\sF$ on 
$\Dmod(\Gr_G)$ equipped with compatible isomorphisms 
\begin{equation} \label{e:heckeeigen}
    \sF \overset {L^+G }\star \on{Sat}(V) \simeq \sF \otimes \on{Oblv}(V), 
    \quad \text{for } V \in \on{Rep}(\cG),
\end{equation}	
where $\on{Oblv}(V)$ denotes the underlying vector space of $V$. A basic 
observation to make about such an $\sF$, for which convolution with any Satake 
sheaves returns a direct sum of copies of $\sF$, is that $\sF$ must have 
infinite dimensional support and, in particular, cannot be a compact object of 
$\Dmod(\Gr_G)$. 

Before comparing the two descriptions of $\sC_\sigma$, let us make a couple 
more orienting comments. 

\begin{re}  One arrives at the more down to earth formulation \eqref{e:heckeeigen} of Hecke eigensheaves 
by dragging the Satake 
sheaf 
through 
the tensor product \eqref{e:heckeeig}, i.e., by considering 
\[
     (\sF \overset {L^+G} \star V) \otimes k \simeq \sF \otimes (V \cdot k) 
     \simeq 
     \sF \otimes \on{Oblv}(V), 
\]
where $\cdot$ denotes the action of $\QCoh(\pt/\cG)$ on $\QCoh(\pt)$. 
	\end{re}

	\begin{re}The reader meeting these ideas for the first time may find some 
	of the details of the preceding discussion somewhat involved. It is therefore worth emphasizing 
	that one can also predict this by analogy. Namely, in the $p$-adic case, to 
	extract an irreducible unramified representation $\pi_{\eta}$ of $G(\Fqt)$ 
	from 
	\[
	  \Fun( G(\Fqt)/G(\Fq[[t]])),
	\]
	one passes to eigenvectors for the right action of the spherical Hecke 
	algebra with eigenvalue $\eta$, i.e., the corresponding homomorphism
	\[   \eta: \Fun(G(\Fq[[t]]) \backslash G(\Fqt)/G(\Fq[[t]]))) \rightarrow k.\]
	 In the geometric setting, there is only one possible eigenvalue by Tannaka 
	 duality,  namely the
	 monoidal functor
	 \[ \on{Oblv}: \on{Rep}(\cG) \rightarrow \Vect,\]
	  and we are again building the irreducible unramified 
	  representation from the universal one by passing to eigenvectors. In particular, the observation we made about the support of Hecke eigensheaves above lifts an analogous assertion about the support of Hecke eigenvectors. 
	\end{re}

	\subsubsection{} By comparing the two constructions of $\sC_{\sigma}$ in 
	the unramified case, one obtains the following prediction. Namely, for any 
	oper structure $\chi$ on the trivial connection, there should be an 
	equivalence
	\begin{equation} \label{e:locaff}
	      \Dmod(\Gr_G) \underset {\QCoh(\pt / \cG)} \otimes \QCoh(\pt) \simeq 
	      \fgc\mod_\chi,
	\end{equation}
	where the right hand side denotes an appropriate renormalization of the 
	category of Kac--Moody representations with central character $\chi$. 
	
	This is in fact an important theorem,\footnote{More carefully, it is a 
	theorem at least under the assumption 
	that $\chi$ corresponds to an oper on the formal punctured disk, cf. 
	Example \ref{e:regop}. However, it is 
	expected that the general case may be handled similarly.} whose development 
	involved significant contributions from multiple groups of authors over 
	several decades. A brief history is as 
	follows. In their 
	work on quantization of the Hitchin system \cite{bdh}, Beilinson--Drinfeld showed 
	that, under the convolution action 
	\[
	     \Dmod(L^+G \backslash LG / L^+G) \otimes \fgc\mod^{L^+G} \rightarrow 
	     \fgc\mod^{L^+G},
	\] 
	any central quotient of the vacuum module is an eigenobject for the action 
	of the abelian Satake category. This theorem 
	was the major local input into their construction of Hecke eigensheaves in  global geometric Langlands 
	corresponding to global opers. 
	
	\subsubsection{} The implications of Beilinson--Drinfeld's work for localization theory at critical level were realized by Frenkel--Gaitsgory, and substantially developed in \cite{frenkelgaitsgorylocgrassmannian04}, \cite{fg06loc}, and \cite{fgloc09}. To begin, the Hecke eigen-property established by Beilinson--Drinfeld  was used by Frenkel--Gaitsgory to construct a 
	functor
	\begin{equation} \label{e:globgr}
	      \Dmod(\Gr_G)\underset {\QCoh(\pt / \cG)} \otimes \QCoh(\pt) 
	      \rightarrow 
	      \fgc\mod_\chi.
	\end{equation}
	Roughly, this amounted to showing that the global sections of any D-module on 
	$\Gr_G$ (twisted by the critical line bundle, so that one obtains a 
	representation at critical level, and with the central quotient corresponding to $\chi$ imposed) is naturally a Hecke eigenobject. For 
	delta D-modules, this is exactly the previous theorem of 
	Beilinson--Drinfeld, and the general case may be reduced to this one. 
	
	 In addition to conjecturing that \eqref{e:globgr} should be a $t$-exact equivalence, Frenkel--Gaitsgory established many of its fundamental properties. Namely, they 
	proved that it was right $t$-exact and a fully faithful embedding on bounded 
	below derived categories.\footnote{The fully faithfulness in general and, 
	in particular, the correct renormalization of the category of Kac--Moody 
	representations, was later obtained by Raskin \cite{RaskinLocCrit}.} They moreover proved  
	that it was 
	an 
	equivalence on $I$-equivariant objects, i.e., highest 
	weight modules. 
	
	\begin{re} With a little work, one obtains from these results the analog of the Kazhdan--Lusztig conjecture for the characters of simple highest weight modules at critical level. This was conjectured by Feigin--Frenkel, see \cite{arakawafiebig}, and is proven in \cite{DhillonYang}. 
	\end{re}

\subsubsection{} After the work of Frenkel--Gaitsgory, the remaining pieces of their conjecture were that \eqref{e:globgr} was $t$-exact and essentially surjective. The latter is profitably reinterpreted as a generation statement under the categorical action of the loop group as follows. The left hand side of \eqref{e:globgr} is by definition generated by 
	its $L^+G$-equivariant objects. Therefore, it remained to show the same held 
	for the right hand side. 
	
	On the one hand, this is a fairly natural sounding 
	statement. Indeed, the central character imposed is of unramified nature, 
	corresponding to an oper on the formal non-punctured disk, and so it is 
	rather plausible that the category itself should be unramified. However, 
	this proved to be a difficult assertion, in large part due to how little 
	was, and largely still is, explicitly known about 
	Kac--Moody representations beyond tame 
	ramification.

	\subsubsection{} Over a decade later, significant progress was made by 
	Raskin \cite{RaskinLocCrit}, who proved the conjecture in the first nontrivial case of $GL_2$ 
	and, in fact, for any group of semisimple rank one. To do so, he 
	proved a general assertion about categorical representations of $LGL_2$, 
	namely their generation by Whittaker and Iwahori invariants, which 
	parallels an earlier and similarly fundamental assertion in the $p$-adic 
	theory. As \eqref{e:globgr} is an equivalence on such objects, the full equivalence then follows.  
	
	The general case of the conjecture was then settled by Raskin--Yang \cite{RaskinYang}. The 
	major new ingredient was the adaptation of Moy--Prasad theory from 
	the $p$-adic setting by Yang \cite{yangmoyprasad}; see also the earlier \cite{chenkamgarpour}. Crucially, besides the filtration of loop group representations by depth, Moy--Prasad theory moreover 
	provided an explicit list of compact generators for each successive 
	quotient, corresponding to unrefined minimal $K$-types in the $p$-adic theory. Their proof then showed by a nontrivial analysis that having a 
	nonzero map from any such generator beyond depth zero is incompatible 
	with an unramified central character, which yielded the desired essential 
	surjectivity. 
	
	\begin{re}In fact, a crucial technical point, which we have swept under the rug in the preceding discussion, is that the category $\fgc\mod_\chi$ does not carry an a priori categorical action of $LG$  due to its renormalization. This is one of the basic difficulties surmounted in \cite{RaskinYang}.
	\end{re}

	\begin{re} We would like to make a remark on Moy--Prasad theory in geometric Langlands, which is otherwise independent of our discussion of localization. Namely, an interesting feature of working with $\cG$-connections as opposed to Galois representation, due to Bremer--Sage \cite{sage1}, \cite{sage2}, \cite{sage3}, is that one also has a theory of minimal $K$-types for $\cG$-connections. The interaction with local geometric Langlands and particularly Moy--Prasad theory on the automorphic side is not yet fully understood. 	\end{re}

	\subsubsection{} Let us finish our discussion of localization theory at critical level with two comments. 
	First, as the reader may imagine, the above admits an extension from the 
	unramified case to the case of tame ramification, e.g.,  with unipotent 
	monodromy. This was developed by Frenkel and Gaitsgory in \cite{fg06loc} and \cite{fglocaff}, roughly by passing from the affine Grassmannian to the affine 
	flag variety.

	Second, it is rather arresting that, given the remarkable nature of the localization theorem at critical level we have discussed, it is in some sense the tip of the iceberg. That	is, it is expected that many of the known constructions of irreducible representations of $p$-adic 
	groups of higher depth should admit geometric versions, leading to equivalences and,  in 
	particular, localization theorems for  critical level representations 
	with 
	wildly ramified central characters; see \cite{fg06loc} and \cite{kamgarpourschedler} for some discussion and conjectures in this direction. 
	
	In particular, the surveyor expects that, for an oper whose underlying local system admits a $\cB$-reduction, the corresponding category of Kac--Moody modules localizes as a category of twisted D-modules on the semi-infinite flag manifold. For ramified characters, such a localization is no longer possible also on the affine flag variety as in \cite{fglocaff}, but instead on the quotient by a smaller compact open subgroup first introduced in the $p$-adic setting by Roche \cite{roche}.

	\subsection{Siblings} \label{ss:variants}
	In this section, we have discussed the local geometric Langlands 
	conjecture, and hopefully the reader has gained a basic feeling for its 
	statement and some of its concrete predictions in representation theory and 
	geometry, both known and conjectural.

	Before we finish, we should mention some important variants or, rather, sibling conjectures in nearby areas, as well as some 
	other emerging connections. 
	
	\subsubsection{Quasi-split groups} The first comment we should make is that we have worked throughout this survey with loop groups, i.e., split reductive groups over $k(\!(t)\!)$. However, even for the study of their categorical representations, as in number theory, it is important to also consider the case of twisted loop groups, i.e., the analogues  of quasi-split groups in the $p$-adic setting. Note this allows for a similar picture of the representation theory of twisted affine Lie algebras, and the forms of local geometric Langlands duality discussed below again are expected to extend to the quasi-split case.

	\subsubsection{Quantum local geometric Langlands} Recall that $\fgk\mod$ 
	is  
	a categorical representation of the loop group only for integral 
	levels. For other levels, it instead carries an action of the category of 
	twisted D-modules $\Dmod_\kappa(LG)$. One can ask if there is a form of 
	local Langlands duality for such twisted categorical representations of 
	loop groups. Indeed there is, and this is the subject of quantum geometric Langlands.   
	
    To proceed, let us for simplicity assume 
    that $\kappa - \kappa_c$ is a nondegenerate bilinear form. In this case, 
    one may attach a dual form $\check{\kappa}$ for $\cG$, as follows. If we 
    write $\kappa_{\cfg, c}$ for the critical level for $\cfg$, then 
    $\check{\kappa}$ is defined by the property 
    that $$\kappa - \kappa_c \text{ and } \check{\kappa} - \kappa_{\cfg, c}$$ 
    are dual 
    bilinear forms on $\ft$ and $\check{\ft}$, respectively. That is, up to the 
    critical shift we have repeatedly met throughout, one asks that $\kappa$ 
    and $\check{\kappa}$ be dual forms.  
    
    The local quantum Langlands conjecture, due to Gaitsgory \cite{quantum-langlands-summary}, is that one has an 
    equivalence of 
    $(\infty,2)$-categories 
    $$\Dmod_\kappa(LG)\mod \simeq \Dmod_{-\check{\kappa}}(L\cG)\mod.$$
    This builds on earlier work in the global setting, due to 
    Beilinson--Drinfeld, Feigin--Frenkel, and Stoyanovsky \cite{stoyanovsky}. 
    
     We confine ourselves here to only a few comments; see however \cite{winter-langlands-summary}, \cite{paris-notes} for more details.

    \begin{re} A basic thing to notice about the quantum case is that the two 
    sides look more symmetric, as they only concern twisted categorical loop 
    group 
    representations. Moreover, in a precise sense, one recovers the usual 
    statement \eqref{e:loclang} in the limit as $\kappa \rightarrow \kappa_c$, 
    and hence $\check{\kappa} \rightarrow \infty$, and one recovers its analogue with the 
    roles of $G$ and $\cG$ reversed in the limit as $\check{\kappa} \rightarrow 
    \kappa_{\cfg, c}$, and hence $\kappa \rightarrow \infty$. 	
    	\end{re}

    \begin{re} The adjective `quantum' arises as follows. Recall the 
    Feigin--Frenkel isomorphism 
    \[      \mathscr{W}_{\fg, \kappa_c} \simeq \OO_{\Op},
    \]
    which plays a basic role in the story at critical level. In fact, 
    Feigin--Frenkel 
    showed in \cite{feiginfrenkel91} that
    this is the quasi-classical limit of a family of isomorphisms of affine 
    W-algebras
    \[   \mathscr{W}_{\fg, \kappa} \simeq \mathscr{W}_{\cfg, \check{\kappa}}. \]
    These isomorphisms, now between non-commutative algebras, 
    deformation quantize the previous isomorphism and play a basic role in the 
    quantum theory. 
    
    Relatedly, while at critical level we saw that  the representations of the 
    dual group $\on{Rep}(\cG)$ play a basic role in 
    the unramified cases of geometric Langlands, the analogous role is played 
    in the quantum setting by representations of the quantum group 
    $\on{Rep}_q(\cG)$.  
    \end{re}

   \subsubsection{Connections with physics} So far, we have only touched on 
   some of the relations between geometric Langlands and physics and mostly on ties to two-dimensional 
   conformal field theory. However, there is an extremely rich connection with three- and four-dimensional supersymmetric gauge theories, starting with the works \cite{kapustinwitten}, \cite{wittenwild}, \cite{gukovwitten}, \cite{frenkelwitten}. This is the 
   source of many interesting constructions and expected compatibilities in 
   (quantum) geometric Langlands, both local and global. We refer the reader to \cite{gaiotto16s-duality}, \cite{gaiottorapcak}, \cite{bfn}, \cite{gaiottotwistcompact}, \cite{costellogaiotto}, \cite{creutziggaiotto}, \cite{creutziggaiottolinshaw}, \cite{frenkelgaiotto}, \cite{bfgt}, \cite{hilburnraskin} for a partial 
   indication of recent work, as well as the forthcoming \cite{gaiottoyu}.

 \begin{re} 
   One exciting aspect of some of the above work is its close relation to new developments in the arithmetic Langlands program. Namely, many predictions of Langlands dual categories coming from
    S-duality of boundary conditions for four-dimensional  
   $\mathscr{N} = 4$ Yang--Mills, as studied by Gaiotto--Witten \cite{gaiottowitten}, match analogous predictions in arithmetic, notably from the relative Langlands program of Sakellaridis--Venkatesh \cite{sakellaridisvenkatesh}. This connection is 
   developed in ongoing work of Ben-Zvi--Sakellaridis--Venkatesh. 
\end{re}

	\subsubsection{}  We have arrived at our final point, 
	which we state first somewhat informally. When working with constructible 
	sheaves on algebraic varieties, there are different categories of sheaves 
	one can choose that behave similarly. Plainly, one typically works  with 
	either (i) D-modules, for which the corresponding cohomology theory, i.e., 
	global sections of the `constant sheaf,' is algebraic de Rham cohomology; (ii) $\ell$-adic sheaves, for which the corresponding cohomology theory is 
	\'etale cohomology; and (iii) when working over the complex 
	numbers, constructible sheaves in the analytic topology, for which the 
	corresponding cohomology theory is singular, i.e., Betti cohomology. 
	 
	In the version of (quantum) geometric Langlands we have discussed, the 
	constructible parts of the correspondence are of de Rham type. That is, on 
	the automorphic side we study D-modules on the moduli of $G$-bundles, and on 
	the spectral side we study coherent sheaves on the moduli of 
	$\cG$-bundles with flat connections, i.e., local systems in the de Rham 
	sense. Relatedly, for a smooth projective curve $X$ over $k$, the global conjecture, in its modern formulation by Arinkin--Gaitsgory \cite{arinkingaitssingsupp}, posits an equivalence
	\[
	\Dmod(\Bun_G(X)) \simeq \ICn(\LS(X)). 
	\]
	As we now describe, there are also versions of the theory of \'etale  
	and Betti flavors.

	\subsubsection{Local geometric Langlands with restricted variation} An 
	important recent development in local arithmetic Langlands has been its 
	formulation in families of representations, i.e. for whole categories and not only irreducible by 
	irreducible \cite{emertonhelm}, \cite{hellmann}, \cite{zhufamilies}, \cite{farguescholze}. That is, its modern formulation is roughly as 
	an equivalence between all
	representations of $G(\Fqt)$ and coherent sheaves on an appropriate moduli 
	space of Galois representations.
	
	\begin{re} In fact, in the spirit of the work of Zelevinsky, Lusztig, and Vogan \cite{voganloclang}, the automorphic side is subtler and is instead also glued with the categories of representations of other groups, e.g., certain inner forms of $G$.  
	\end{re} 
	
    In the series of works \cite{restrictedvar1}, \cite{restrictedvar2}, \cite{restrictedvar3}, Arinkin--Gaitsgory--Kazhdan--Raskin--Rozenblyum--Varshavsky 
    have produced a global categorical geometric Langlands conjecture that 
    makes sense both in characteric zero and over function fields. It takes the 
    form 
    \begin{equation} \label{e:etalelang}
           \on{Shv}^{\on{\acute{e}tale}}_{\on{nilp}}(\on{Bun}_G(X)) \simeq 
           \on{IndCoh}_{\on{nilp}}(\on{LocSys}_{\cG}^{\on{restr}}(X)).  
    \end{equation}
	On the left hand side, one is considering ind-constructible $\ell$-adic sheaves on the moduli of $G$-bundles on a global curve $X$ with singular support in the global nilpotent cone of $T^* \Bun_G(X)$.  On the right hand side, one meets ind-coherent sheaves with nilpotent singular support 
	on an appropriately defined moduli space of $\ell$-adic $\cG$ local systems on $X$. 
	
	\begin{re}
	Very roughly speaking, the basic picture is that, over the complex numbers,	the spectral side looks like a `torn apart' version of the Betti moduli 
	stack of $\cG$-local systems, i.e., the moduli space parametrizing representations of the topological fundamental group. That is, one decomposes the Betti stack into many 
	pieces 
	based on the semisimplification of the local system and 
	takes the disjoint union of their formal neighborhoods. For example, the formal 
	neighborhood of an irreducible 
	local system is isolated in  $\on{LocSys}_{\cG}^{\on{restr}}$ and 
	identifies with the corresponding formal completion in the Betti or de Rham moduli 
	spaces. Relatedly, the superscript `restr' stands for restricted variation. 
	\end{re}

	In the function field setting, upon taking the categorical trace of 
	Frobenius, this was shown to recover (compactly supported) automorphic 
	forms and, in particular, the automorphic-to-Galois direction of unramified 
	global Langlands for cuspidal representations over function fields. 
	
	In unpublished work, Gaitsgory has formulated 
	a local version of \eqref{e:etalelang}. This roughly asserts an equivalence
	\[
	     \on{Shv}^{\on{\acute{e}tale}}(G(\overline{\mathbb{F}}_q(\!(t)\!))\mod \simeq 
	     \TwoIC(\LS^{\on{restr}}),
	\] 
	where on the left hand side one really works with a subtler notion of 
	categorical representation to correct the non-essential surjectivity of 
	the fully faithful embedding
	$$\on{Shv}^{\on{\acute{e}tale}}(X) \otimes \on{Shv}^{\on{\acute{e}tale}}(Y) \hookrightarrow \on{Shv}^{\on{\acute{e}tale}}(X \times Y)$$
	and passes further to a full subcategory of `spectrally finite' objects. 
	
	Via a categorical trace of Frobenius, this is expected in the function 
	field setting to recover a form of the  arithmetic local Langlands conjecture in its 
	formulation by Zhu \cite{zhufamilies}.

	\subsubsection{Local Betti geometric Langlands} Finally, a version of 
	geometric Langlands for Betti sheaves has been formulated by 
	Ben-Zvi--Nadler \cite{benzvinadlerbetti}. 
	Here, the global equivalence, for $X$ a projective algebraic 
	curve over the complex numbers, takes the form 
	\[
	  \on{Shv}^{\on{Betti}}_{\on{nilp}}(\Bun_G(X)) \simeq 
	  \ICn(\LS^{\on{Betti}}(X) ).
	\] 
	On the left hand side, one is considering Betti sheaves, not necessarily 
	ind-constructible, on the analytification of $\Bun_G(X)$ with singular 
	support in the global nilpotent cone. On the right hand 
	side, one is considering ind-coherent sheaves on the Betti moduli space of 
	local systems, i.e., the natural enhancement of the character variety of the 
	underlying topological space of $X$ to a derived stack.

	The local case is not particularly well documented, but let us sketch some 
	basic parts.  First, the portion with regular singularities is established 
	by the equivalence
	\begin{equation} \label{e:bezfam}
	      \on{Shv}^{\on{Betti}}_{\on{nilp}}( \mathring{I}\backslash LG / 
	      \mathring{I}) \simeq \ICn( \widetilde{\cG} / \cG \underset {\cG / 
	      \cG} \times \widetilde{\cG} / \cG).
	\end{equation}
	Informally speaking, this is a version of Bezrukavnikov's equivalence 
	\eqref{e:2realz} and its variants for nonunipotent monodromy in families, where the monodromy varies through the adjoint quotient $\check{G} / \check{G}$.\footnote{This equivalence may not be available in the literature. However, the surveyor has an argument, joint with H. Chen, which hopefully will be recorded soon.}
	
	By construction, the modules for the right hand side of \eqref{e:bezfam} 
	identify with ind-coherent sheaves of categories with nilpotent singular 
	support on $\LS^{\on{Betti}}$, i.e., $\cG / \cG$. The full assertion, i.e., 
	including the counterpart of wild ramification, 
	should relate Betti categorical representations of the 
	loop group, which is to first approximation 
	$\on{Shv}^{\on{Betti}}(LG)\mod$, 
	with ind-coherent sheaves of categories with nilpotent singular support on 
	the moduli of Stokes data.

	\subsubsection{} At this point, we have arrived at the end of the survey. We hope the reader has gained a feel for some of the basic ideas, results, and conjectures in this area, as well as some appreciation for their beauty. The interested reader is encouraged to wade into the literature we have toured, and, as is particularly feasible in an area with so many attractive open problems, further learn by doing.

	\appendix

	\section{From functions on $X(\Fq)$ to D-modules on $X(\CC)$}
	\label{s:frmfun}
	
	\subsection{Overview}

	\subsubsection{} The goal of this appendix, which is a supplement to Section \ref{s:bbloc}, is to explain in what sense passing from functions on the $\Fq$-points of the flag variety to D-modules on the $\CC$-points of the flag variety is a natural thing to do. 
	
    The short answer is that one thinks of the latter as some sort of categorification of the former. Any explanation of this necessarily involves three topics that play an important role in geometric representation theory -- D-modules, the functions-sheaves correspondence, and the Riemann--Hilbert correspondence. 
    
    We hope the following presentation of these topics gives enough context and intuition for nonspecialists to appreciate the discussion in the main body of the text. We will defer the suggestion of references to the end of each subsection. 
	
	\subsection{Algebraic D-modules}

	\subsubsection{}  In this subsection we would like to discuss the basics of D-modules. We begin with the affine case. 
	
	So, suppose $X$ is a smooth and affine algebraic variety 
	over a field $k$, which for simplicity we take to be of characteristic zero. Informally, a D-module on $X$ is something in which it 
	makes sense to scale by functions on $X$ and to differentiate by vector 
	fields on $X$. 
	
	More formally, let us denote by $O_X$ the algebra of regular functions on $X$, and 
	by $T_X$ the $O_X$-module of regular vector fields on $X$. A D-module $M$ on $X$ 
	is then an $O_X$-module equipped with a flat connection. That is, $M$ is a 
	vector space equipped with an action of the algebra of functions 
	\[
	O_X \otimes M \rightarrow M, \quad f \otimes m \mapsto f \cdot m,  
	\]
	as well an action by the Lie algebra of vector fields
	\[
	T_X \otimes M \rightarrow M, \quad \xi \otimes m \mapsto \xi \cdot m.
	\]
	That the above is a Lie algebra action is exactly the assumption of flatness, i.e., vanishing curvature. These two operations are asked to be compatible with scaling vector fields by functions 
	\[
	f \cdot (\xi \cdot m) = (f \cdot \xi) \cdot m
	\]
	and with differentiating functions by vector fields, i.e., to satisfy the `Leibnitz rule'
	\[
	\xi \cdot (f \cdot m) = (\xi \cdot f) \cdot m  + f \cdot (\xi \cdot m),
	\]
	where $\xi \cdot f$ denotes the derivative of $f$ with respect to $\xi$.

	Equivalently, one can form the algebra $D_X$ of polynomial linear 
	differential operators on $X$, i.e., the subalgebra of all endomorphisms of $O_X$ generated by multiplication by functions and differentiation by vector fields. With this, a D-module is simply a left 
	$D_X$-module.

	\begin{ex} For the affine space $X = \mathbb{A}^n_k$, the regular functions 
	and vector fields take the form 
		\[
		O_X \simeq k[x_1, \ldots, x_n] \quad \text{and} \quad T_X \simeq 
		\underset{i}{\oplus} \hspace{1mm} O_X \cdot \partial_i, 
		\]
		where the $\partial_i$, $1 \leqslant i \leqslant n$, denote the standard 
		coordinate vector fields associated to the coordinates $x_i$. In this case, $D_X$ 
		is generated as a $k$-algebra by $x_i, \partial_j$, $1 \leqslant i,j 
		\leqslant n$, subject to the standard commutation relations 
		\[
		[x_i, x_j] = [\partial_i, \partial_j] = 0, \quad \quad [\partial_i, 
		x_j] = \begin{cases} 1 & \text{if } i = j, \\ 0 & \text{otherwise}. 
		\end{cases}
		\]
		Explicitly, $D_X$ has as a basis the ordered monomials 
		\[
		x_1^{a_1}    \cdots  x_n^{a_n} \cdot \partial_1^{b_1} \cdots 
		\partial_n^{b_n}, \quad \quad \text{for } a_i, b_j \in 
		\mathbb{Z}^{\geqslant 0}, 
		\]
		i.e., consists of linear differential operators with polynomial 
		coefficients. 
	\end{ex}
	
   We write $\on{D-mod}(X)^\heartsuit$ for the abelian category 
	of D-modules on $X$ and denote its unbounded derived category by 
	$\on{D-mod}(X)$. Let us discuss three basic objects of either category, i.e., examples of D-modules. 
	
	\begin{ex} The algebra of differential operators ${D}_X$ is naturally a 
	left module over itself, i.e., a D-module. 
	\end{ex}
	
	Next, we turn to a slightly less tautological example of `global nature.' 
	
	\begin{ex} For any $X$, the ring of regular functions $O_X$ carries a 
	tautological structure of a D-module. Namely, one multiplies and 
	differentiates functions in the standard way.

		This module admits the following presentation via generators and relations. 
		Acting on the constant function 1 yields an isomorphism of D-modules
		\[
		D_X/ D_X \cdot T_X \simeq O_X. 
		\]
		In particular, giving a map of D-modules $O_X \rightarrow M$ is the 
		same as specifying a flat section of $M$, i.e., an element $m$ 
		annihilated by $T_X$.

		More generally, any vector bundle $E$ equipped with a flat connection 
		$\nabla$ is a D-module, and giving a map of D-modules $E \rightarrow 
		M$ is the same as specifying a flat section of the tensor product$$M 
		\underset{\phantom{ii} O_X} \otimes E^\vee$$equipped with its natural flat connection. \end{ex}
	
	Finally, we would like to describe a class of examples of `local nature.' 
	\begin{ex} For any closed point $x$ of $X$, one has an associated delta 
	function D-module. Namely, if we write $\fm_x \subset O_X$ for the maximal 
	ideal of functions vanishing at $x$, we set 
		\[
		\delta'_x := D_X / (D_X \cdot \fm_x). 
		\]
		In particular, giving a map of D-modules $\delta'_x \rightarrow M$ is 
		the same as specifying a section of $M$ annihilated by $\fm_x$, i.e., 
		scheme-theoretically supported at $x$.  
		
		The D-module $\delta'_x$ is an algebro-geometric version of the vector 
		space of distributions supported at $x$.\footnote{This parallel is 
		improved by passing to the corresponding right D-module, but the reader 
		may wish to ignore this point.} Relatedly, if $\xi_1, \ldots, \xi_n$ 
		are a collection of vector fields in $T_X$ whose image in the tangent 
		space $T_x X$ at $x$ forms a basis, then the images of the monomials 
		\begin{equation} \label{PBW}
		\xi_1^{a_1} \cdots \xi_n^{a_n}, \quad \text{for } a_i \in 
		\mathbb{Z}^{\geqslant 0}
		\end{equation}
		form a basis for $\delta'_x$. One thinks of the image of 1 from $D_X$ 
		as the usual delta distribution, and the entire module is then spanned by 
		its derivatives. 
		
		Here is a somewhat technical point, which the reader may safely ignore 
		at first pass. It turns out to be better to twist the underlying vector 
		space of $\delta'_x$ by a line, namely the determinant of the tangent 
		space at $x$:
		\[
		\delta_x := \delta'_x \underset{k} \otimes \det(T_xX).
		\]
		The reason is as follows. Given a map $f: X \rightarrow Y$ of smooth affine varieties, one has an associated pushforward of D-modules $$f_*: 
		\on{D-mod}(X) \rightarrow \on{D-mod}(Y),$$cf. Section \ref{ss:pushpull} below. When we incorporate the 
		twist, we then have canonical isomorphisms
		\[
		f_*( \delta_{x} ) \simeq \delta_{f(x)}, \quad \quad \text{for } x \in X.
		\]
		That is, delta functions pushforward to delta functions, as one expects 
		from the distributional analogue. \end{ex}

	\begin{re} \label{r:hwdelta}Here is an alternative construction of 
	$\delta_x$, which may 
	help elucidate the appearance of the determinant line. 
		
		Consider the Lie algebra of vector fields vanishing at $x$, i.e., $\fm_x \cdot T_X \subset T_X.$ This acts naturally on the quotient $T_X/ 
		\hspace{.1mm} \fm_x \cdot T_X \simeq T_x X.$ Explicitly, the action 
		factors through 
		\[
		\fm_x \cdot T_X / \hspace{.1mm} \fm_x^2 \cdot T_X \simeq T^*_x X 
		\otimes T_x X \simeq \on{End}(T^*_xX), 
		\]
		and the action on $T_x X$ is simply the natural action, 
		for any vector space $W$, of its Lie algebra of endomorphisms 
		$\on{End}(W)$ on the dual vector space $W^*$. 
		
		Consider the induced action of $\fm_x \cdot T_X$ on the determinant line $\det(T_xX)$. This line also  
		carries a natural action of $O_X$, and these two actions define 
		an action of the subalgebra $D_{X, x}$ of $D_X$ generated by $O_X$ and 
		$\fm_x \cdot T_X$. A point here is that, before passing to the 
		determinant, the relevant compatibilities would not be satisfied. 
		
		Finally, the induced D-module is equipped with a canonical isomorphism 
		\[
		D_X \underset{D_{X,x}} \otimes \det(T_xX) \simeq \delta_x,   
		\]
		which acts as the identity on the generating determinant lines, i.e., 
		we have described the action of $D_{X, x}$ on the span of the delta 
		distribution itself. \end{re}

	As the reader may guess, many D-modules have behavior which interpolates 
	between the last two examples. That is, they resemble an algebraic connection along a subvariety and a delta D-module along its normal directions.

	\subsubsection{} \label{ss:pushpull}Given a map $X \rightarrow Y$ of smooth affine varieties, one has natural pullback and pushforward functors between their categories of D-modules, i.e., 
	\[
	  f^!: \Dmod(Y) \rightarrow \Dmod(X) \quad \text{and} \quad \quad f_*: \Dmod(X) \rightarrow \Dmod(Y).
	\]
	As for the nature of the functors, for our purposes it suffices to bear the following in mind. On a manifold $M$, one has the algebra of smooth functions and its dual space of distributions. Given a map $M \rightarrow N$, one can pull back functions and dually push forward distributions. 
	
	Intuitively, $f^!$ records, for a function on $Y$ satisfying a linear PDE, the linear PDEs satisfied by its pullback. Similarly, $f_*$ records, for a distribution on $X$ satisfying a linear PDE, the linear PDEs satisfied by its pushforward. In particular, $f^!$ and $f_*$ are not always adjoint functors, but rather are dual. This can be made precise using duality of dg-categories as in Section \ref{ss:dualizable}.

	 We now provide a more detailed discussion of $f^!$ and $f_*$ for the reader's convenience, though they may safely skip this. The pullback $f^!$ is simply given by the quasi-coherent pullback up to a cohomological shift, i.e., 
    \[
         f^!(M) \simeq O_X \underset{\phantom{ii}O_Y} \otimes M [d_X - d_Y],
    \]
	where we denote by $d_X$ and $d_Y$ the dimensions of $X$ and $Y$, respectively. Here, $O_X$ acts in the tautological way, and $T_X$ acts via the chain rule, i.e., via the map $$T_X \rightarrow O_X \underset {O_Y} \otimes T_Y$$and the natural action of the latter on $f^!(M)$. 
	To describe $f_*$, note the previous formula for the pullback can be rewritten as tensoring over $D_Y$ with the bimodule
	\[
	    O_X \underset {\phantom{ii}O_Y} \otimes D_Y.
	\]
	The same bimodule also defines the pushforward functor $f_*$. More carefully, it yields a functor between categories of right D-modules $$f_{\on{r}, *}: \on{mod-D}(X) \rightarrow \on{mod-D}(Y).$$
	The discrepancy is bridged by the fact that left D-modules and right D-modules are canonically identified. Namely, on a smooth affine variety $Z$, the canonical bundle $K_Z$ carries a natural right action of $D_Z$. Roughly, one thinks of the canonical bundle as an algebraic substitute for distributions and the action is by integration by parts. Tensoring by it defines the desired equivalence 
	\[
	   - \underset {\phantom{ii}O_Z} \otimes K_Z: \Dmod(Z) \simeq \on{mod-D}(Z), 
	\]
	and applying this on both sides we obtain the pushforward
	\[
	  f_*: \Dmod(X) \simeq \on{mod-D}(X) \rightarrow \on{mod-D}(Y) \simeq \Dmod(Y). 
	\]
	%


	\subsubsection{} The previous definition of a D-module in the affine case 
	globalizes straightforwardly, as follows. Suppose now $X$ is a smooth 
	algebraic variety over $k$, not necessarily affine. 
	
	A D-module $\sM$ on $X$ is a quasi-coherent sheaf $\sM$ equipped with a 
	flat connection. That is, if we write $\OO_X$ for the sheaf of regular 
	functions, $\sT_X$ for the sheaf of regular fields, and $\underline{k}$ for 
	the constant sheaf associated to $k$, $\sM$ is a sheaf of $k$ vector spaces 
	equipped with maps 
	\[
	\OO_X \underset{\ck}{\otimes} \sM \rightarrow \sM \quad \text{and} \quad 
	\sT_X \underset{\ck}{\otimes} \sM \rightarrow \sM,
	\]
	satisfying compatibilities as in the previous section. We recall that 
	quasicoherence is the further condition that, for any affine open $U$ and 
	regular function $f$ on $U$ with nonvanishing locus $U_f$, the natural map 
	\[
	\Gamma(U, \sM) \underset{\phantom{ii} O_U} \otimes O_{U_f} \xrightarrow{\sim} \Gamma(U_f, \sM)
	\]
	is an isomorphism. By this assumption, if $X$ is affine, the present 
	sheaf-theoretic category of D-modules is canonically equivalent to the one 
	of the previous subsection, via the functor of taking global sections.

	So, somewhat informally, the category of D-modules on a general smooth 
	variety is glued from the previously discussed categories of D-modules on 
	its affine open subvarieties. Similarly, a D-module $\sM$ is roughly a 
	sheaf whose local sections one can scale by regular functions and 
	differentiate by regular vector fields. 
	
	The examples of the previous subsection carry over to this setting as 
	follows. 
	
	\begin{ex} One has a global sheaf $\sD_X$ of differential operators on $X$, 
	and D-modules are quasicoherent sheaves of $\sD_X$-modules.
	\end{ex}

	\begin{ex} For any $X$, the sheaf of regular functions $\OO_X$ is naturally a  D-module, as is the sheaf of sections of any 
	vector bundle with flat connection $(\sE, \nabla)$. 
	\end{ex}
	
	\begin{ex} For an closed point $x$ of $X$, one again has a delta function 
	D-module supported at $x$, which we denote by $\delta_x$. 
	\end{ex}
	
	\label{ss:dmodglob}
	\subsubsection{}  Finally, we note that the D-modules on $X$, as discussed 
	in Section \ref{ss:dmodglob}, naturally form an abelian category in such a 
	way that the forgetful functor to sheaves of abelian groups is $t$-exact. That is, one calculates kernels, cokernels, and direct sums in the usual way 
	for sheaves, and these naturally inherit structures of D-modules and 
	satisfy the correct universal properties. 
	
	\begin{defn} We write $\on{D-mod}(X)^\heartsuit$ for the abelian category 
	of D-modules on $X$. We denote its unbounded derived category by 
	$\on{D-mod}(X)$.  
	\end{defn}
	
	Given a map $f: X \rightarrow Y$, one again has pushforward and pullback functors 
	\[ f^!: \Dmod(Y) \rightarrow \Dmod(X) \quad \text{and} \quad f_*: \Dmod(X) \rightarrow \Dmod(Y), \]
	defined as in Section \ref{ss:pushpull}. We should highlight that in the non-affine case, $f_*$ is the composite of the left derived functor of tensoring with the bimodule $$\sO_X \underset{\sO_Y}{\otimes} \sD_Y$$and the right derived functor of pushforward of sheaves along $f$. As such, it is not the derived functor of a functor between abelian categories of D-modules, and requires working with derived categories from the outset. 
	
	\subsubsection{} 	For further reading on D-modules, the reader may consult, as some of many possible references, ordered roughly by increasing demands on the reader,  \cite{coutinho}, \cite{arkhipovnotes},  \cite{bernsteinnotes}, \cite{kashiwarabookd-modules}, \cite{gaitsgorynotes}, \cite{hottatakeuchitanisaki}.

	\subsection{The functions-sheaves correspondence}
	
		\subsubsection{}
	
	We have, as promised, made sense of the category of D-modules on 
	the complex flag variety
	\[
	\on{D-mod}(\sB(\C)). 
	\]
We would like to next explain the sense in which this is analogous to the vector 
	space of functions $$\Fun(\sB(\Fq)).$$ 
	\subsubsection{} For any $\Fq$-variety $X$ and field $k$ of characteristic 
	zero, recall we associated to it 
	\[
	\Fun(X(\Fq), k),
	\]
	the vector space of $k$-valued functions on the finite set of rational 
	points $X(\Fq)$. 
	
	A wonderful fact of life is that special functions in $\Fun(X(\Fq), k)$ are 
	often the shadows of special sheaves on $X$. Let us say this more 
	precisely. 
	
	\subsubsection{}\label{s:ladic} To do so, let us take $k = \Qlb$, an algebraic closure of 
	the $\ell$-adic numbers $\mathbb{Q}_\ell$, where $\ell$ is any prime number coprime to $q$. 
	
	With this choice of $k$ one has a well-behaved category of \'{e}tale local 
	systems of $k$ vector spaces on $X$. This category is 
	equivalent to the category of continuous representations of the \'{e}tale 
	fundamental group of $X$ on $k$ vector spaces.

	\begin{re} Here is an orienting remark for the reader unfamiliar with these 
	notions. Consider a smooth complex algebraic variety $X$ equipped with a 
	local system $\sL$, by which we mean a locally constant sheaf of finite dimensional $k$ vector spaces, in 
	the sense of the analytic topology.

		The basic observation to make is that one can hardly ever trivialize 
		$\sL$ on a Zariski open subset of $X$. Concretely, the reader may wish 
		to think through the case of $X = \mathbb{G}_m$ and $\sL$ a rank one 
		local system whose monodromy around the origin is a root of unity.
		
		However, if the corresponding representation of the fundamental group 
		of $X$ factors through a finite quotient, one can trivialize $\sL$ 
		after pulling back to a finite covering space $\tilde{X}$ of $X$. 
		Crucially, $\tilde{X}$ is again canonically an algebraic variety -- we 
		emphasize that this is a nonobvious fact in general. In the previous 
		example, if the monodromy is an $n^{th}$ root of unity, one can take 
		the cover to be 
		\[
		\mathbb{G}_m \rightarrow \mathbb{G}_m, \quad \quad z \mapsto z^n. 
		\]

		The covering map $\tilde{X} \rightarrow X$ is a prototypical example of 
		an \'{e}tale, i.e., a flat and unramified, map of algebraic varieties. The 
		upshot is that if one allows oneself to check local properties of 
		sheaves on not just Zariski open sets, but also \'{e}tale covers, one 
		obtains a robust theory of local systems which makes sense for general 
		schemes, not necessarily over the complex numbers. \end{re}

	More generally, if one wants to consider local systems also on locally closed 
	subvarieties of $X$, these all naturally belong to, and in fact generate, the bounded derived 
	category of constructible \'{e}tale sheaves 
	$
	\on{Sh}(X)^b. 
	$
	\subsubsection{} Let us explain what this category looks like in the 
	simplest case of $X$ being a point, i.e., the spectrum of a field $\kappa$. 
	
	If $\kappa$ is separably closed, its constructible derived category of sheaves of $k$ vector spaces is 
	simply the bounded derived category of $k$ vector spaces, i.e. 
	\[
	\on{Sh}(\on{Spec }\kappa)^b \simeq \on{Vect}_k^b.
	\]
	Plainly, this is because $\on{Spec }\kappa$ does not have interesting 
	\'{e}tale covers, and this should match the reader's intuition for what 
	happens for usual sheaves on a point in topology or complex geometry. We emphasize that $k$ plays the role of the coefficients for the sheaves, and $\kappa$ provides the geometric object $\Spec \kappa$ on which the sheaves live. In particular, there need be no relation between $\kappa$ and $k$.

	However, for a general field $\kappa$, $\on{Spec } \kappa$ {does} have 
	interesting \'{e}tale covers. For example, given a separable extension of 
	fields $\kappa \rightarrow \kappa'$, the map 
	\[
	\Sp \kappa' \rightarrow \Sp \kappa
	\]
	is \'{e}tale, and any finite \'{e}tale cover is a disjoint union of such 
	maps. For this reason, writing $\kappa_s$ for a separable closure of 
	$\kappa$ and $H := \on{Gal}(\kappa_s/\kappa)$ for its Galois group, one 
	has a canonical equivalence between constructible sheaves on $\Sp \kappa$ and 
	the bounded derived category of finite dimensional continuous representations of $H$ 
	on $k$ vector spaces, i.e.
	\[
	\on{Sh}(\Sp \kappa)^b \simeq \on{Rep}_k(H)^b. 
	\]
	In particular, the category of \'{e}tale sheaves captures arithmetic 
	information about the field $\kappa$. An analogue to $\Sp \kappa$ in usual 
	topology would not be a point, but rather the classifying space $\on{pt}/H$ 
	of a (profinite) group $H$. In this analogy, one thinks of the map $\Sp 
	\kappa_s \rightarrow \Sp \kappa$ as akin to the map $\on{pt} \rightarrow 
	\on{pt}/H$. 
	
	\label{ss:shvpt}

	\subsubsection{} Given a map $f: X \rightarrow Y$ of schemes, one has 
	pushforward and pullback functors\footnote{One further has a 
	$!$-pushforward and a $*$-pullback, but we will not need these.} 
	\[
	f_*: \on{Sh}(X)^b \rightarrow \on{Sh}(Y)^b \quad \text{and} \quad f^!: 
	\on{Sh}(Y)^b \rightarrow \on{Sh}(X)^b. 
	\]
	Here $f_*$ is the usual pushforward of sheaves, and $f^!$ is the so-called 
	exceptional inverse image. 
	
	If $f$ is a proper map (recall this roughly means that the fibres of the map are 
	projective varieties), $f_*$ is left adjoint to $f^!$. In particular, given 
	a closed point $x$ of $X$, with residue field $\kappa_x$, one has an 
	adjunction
	\[
	i_{x, *}: \on{Sh}(\on{Spec }\kappa_x)^b \rightleftarrows \on{Sh}(X)^b: 
	i_x^!.  
	\]
	Concretely, $i_{x, *}$ produces skyscraper sheaves supported at $x$, and 
	its adjoint $i_{x}^!$ sends a sheaf to its (derived) sections supported at 
	$x$. 
	
	\subsubsection{} Let $X$ again be a variety over $\Fq$. We are now 
	ready to describe the promised map
	\[
	\on{FF}: \on{Sh}(X)^b \rightarrow \Fun(X(\Fq)), 
	\]
	which is known as Grothendieck's functions-sheaves correspondence. 
	For each rational point $x: \Sp \Fq \rightarrow X$, we therefore need 
	to produce a map 
	\[
	\on{Sh}(X)^b \rightarrow k,
	\]
	i.e., the values of the associated functions at $x$. We remind the reader that,  as in Section \ref{s:ladic},  $k$ denotes the algebraic closure $\overline{\mathbb{Q}}_\ell$ of the $\ell$-adic numbers. To define the desired evaluation, we first pass 
	to the $!$-stalk at $x$, i.e., consider the composition 
	\[
	\on{Sh}(X)^b \xrightarrow{i_x^!} \on{Sh}(\Sp \Fq)^b \simeq 
	\on{Rep}_k(\on{Gal}(\overline{\mathbb{F}}_q / \Fq))^b. 
	\]
	From here, we recall that $\on{Gal}(\overline{\mathbb{F}}_q / \Fq)$ is 
	simply the profinite completion $\widehat{\mathbb{Z}}$ of the integers, 
	with topological generator the Frobenius automorphism $\on{Fr}$. 
	
	In particular, by taking the trace of Frobenius, we obtain the desired map 
	\[
	\on{Sh}(X)^b \xrightarrow{i^!_x} 
	\on{Rep}_k(\on{Gal}(\overline{\mathbb{F}}_q/\Fq))^b 
	\xrightarrow{\on{tr}(\on{Fr}, -)} k. 
	\]
	Explicitly, this sends a sheaf $\sS$ to the number 
	\[
	\underset{j}{\Sigma} \hspace{1mm} (-1)^j \cdot \tr(\on{Fr}, H^j( i^!_x 
	\sS)). 
	\]

	\begin{re} In case this seems remote to the reader, let us pursue the analogy 
	of Section \ref{ss:shvpt}. Replacing the completion $\widehat{\mathbb{Z}}$ 
	by the usual integers $\mathbb{Z}$, note that $\on{pt} / \mathbb{Z}$ is 
	simply the circle $S^1$. That is, one can think of $\Sp \Fq$ as something 
	like a circle.  
		
		One can then think of an $n$-dimensional variety $X \rightarrow \Fq$ as 
		akin to a $2n+1$ dimensional real manifold fibred over a circle $\pi: M 
		\rightarrow S^1$ and its rational points as a collection of embedded 
		circles in $M$ providing sections of $\pi$. The analogue of the trace of Frobenius would send a 
		local system on $M$ to the traces of its monodromies about these 
		circles. \end{re}

	\subsubsection{} A basic property of the functions-sheaves correspondence 
	is that it intertwines the $*$-pushforward of functions with the $*$-pushforward of 
	sheaves and similarly for $!$-pullbacks. For pullbacks, this is by 
	definition, but, for pushforwards, one needs the Lefschetz fixed point 
	formula in \'{e}tale cohomology. 
	
	In addition, since $i^!$ sends distinguished triangles to distinguished 
	triangles,\footnote{Recall that this simply means a sequence $A \rightarrow B 
	\rightarrow C$ of maps in $\on{Sh}(X)^b$, which up to isomorphism comes 
	from a short exact sequence of complexes. More intrinsically,  
	$A$ is the homotopy kernel of $B \rightarrow C$, and, equivalently, $C$ is the 
	homotopy cokernel, i.e., cone, of $A \rightarrow B$.} we obtain a map of 
	$k$ vector spaces 
	\[
	\on{FF}: K( \on{Sh}(X)^b) \underset{\mathbb{Z}}{\otimes} k \rightarrow 
	\Fun(X(\Fq)),
	\]
	where $K(-)$ denotes the Grothendieck group, so that one may reasonably 
	think of $\on{Sh}(X)^b$ as a `categorification' of the space of functions. 
	This is a useful point of view, although we mildly caution that the above 
	arrow is typically not an isomorphism. 
	
	\subsubsection{} 
	
	Let us describe two basic examples of the correspondence. 
	
	\begin{ex} Let us denote by $\underline{k}$ the constant sheaf on $X$. If 
	$X$ is smooth and $n$-dimensional, then $\on{FF}(\underline{k})$ is the 
	constant function on $X(\Fq)$ with value $q^{n}$. \end{ex}
	
	\begin{re} Because we work with $!$-stalks, the sheaf which `lifts' the 
	constant function 1 on any variety $X$ is not the constant sheaf but the 
	dualizing sheaf $\omega_X$. This is by definition obtained by pulling back along the projection $\pi: X \rightarrow \on{Spec} \mathbb{F}_q$ the trivial representation $k$ of $\on{Gal}(\overline{\mathbb{F}}_q / \mathbb{F}_q)$, i.e. $$\omega_X \simeq \pi^!(k).$$ The reader should be aware that many authors work instead with the Verdier dual convention, i.e., with $*$-stalks, when setting up the functions-sheaves correspondence, in which case the constant sheaf, and not the dualizing sheaf, would correspond to the constant function 1. \end{re}

	\begin{ex} Given a rational point $x: \Sp \Fq \rightarrow X$, consider the 
	skyscraper sheaf $i_{x, *} k$, where we view $k$ as a trivial 
	$\on{Gal}(\overline{\mathbb{F}}_q / \Fq)$-module, concentrated in 
	cohomological degree zero. Then we have 
		\[
		\on{FF}( i_{x, *} k) = \delta_x,
		\]
		i.e., the skyscraper sheaf `lifts' the indicator function of the point. 
	\end{ex}
	
	Although these simple examples hopefully  convey the basic idea 
	to the reader, we should mention that much more interesting numbers and 
	functions, e.g., Gauss sums, Kloosterman sums, and characters of finite groups 
	of Lie type, have sheaf-theoretic counterparts. Another rich 	example, which gave birth to the geometric Langlands theory, is the lift of 
	unramified Hecke eigenfunctions over function fields to Hecke eigensheaves on 
	moduli spaces of $G$-bundles on curves over finite fields. 
	
	\subsubsection{} For further details the reader may wish to consult \cite{milnenotes}, \cite{milneetalebook}, \cite{freitagbook}.

	\subsection{The Riemann--Hilbert correspondence}
	
	\subsubsection{}  Having passed from $\Fun(X(\Fq))$ to $\on{Sh}(X)^b$, we 
	can make sense of the latter for $X$ over any base field and not just 
	$\Fq$. 
	
	Over the complex numbers, it will be convenient to work instead with the 
	closely related category of constructible analytic sheaves 
	\[
	\on{Sh}(X_{\on{an}})^b.
	\]
	Here $X_{\on{an}}$ denotes the analytification of $X$, i.e., its set of closed points $X(\C)$ equipped with their metric topology. Explicitly, $\on{Sh}(X_{\on{an}})^b$ is the full subcategory of the derived category of sheaves 
	of $k$ vector spaces generated by sheaves that are locally constant along 
	a stratification of $X$ into Zariski locally closed subsets. 
	
	\begin{re} The category $\on{Sh}(X_{\on{an}})^b$ is an enlargement of the 
	analogous category of \'{e}tale sheaves, essentially by the Artin 
	comparison theorem. For example, if $X = \mathbb{G}_m$, a rank one local 
	system in the \'{e}tale category must have monodromy in the units of the
	algebraic integers, i.e.,  $\overline{\mathbb{Z}}_\ell^{\times}$, whereas in the 
	analytic category  the monodromy may be any element of 
	$\overline{\mathbb{Q}}_\ell^\times$. 
		
	\end{re}

	\subsubsection{} Finally, let us obtain the desired connection to 
	D-modules. First, fix an isomorphism $\overline{\mathbb{Q}}_\ell \simeq 
	\mathbb{C}$ so that we can identify 
	$
	\on{Sh}(X_{\on{an}})^b
	$
	with its analogue with coefficients in complex vector spaces. 
	
	Recall that $\on{Sh}(X_{\on{an}})^b$ contains objects like local systems 
	supported on subvarieties of $X$. These may be converted into a full 
	subcategory of D-modules by an elaboration of the operation of passing 
	between systems of differential equations and their local solutions, as 
	follows.

	\subsubsection{} Consider first the subcategory of the constructible 
	derived category consisting of local systems on $X$. Let us assume that $X$ 
	is smooth. In this case, one has an equivalence of categories 
	\begin{equation} \label{e:vivaflatsections}
	\{ \text{finite rank local systems on $X_{\on{an}}$} \} \simeq \{ 
	\text{flat holomorphic connections on $X_{\on{an}}$} \}, 
	\end{equation}
	where the right hand side denotes the category of holomorphic vector 
	bundles $\sE$ equipped with a holomorphic flat connection $\nabla$,  in the 
	sense of complex geometry. 
	
	Let us review how this equivalence works. To go from right to left in \eqref{e:vivaflatsections}, one 
	sends a connection $(\sE, \nabla)$ to its sheaf of flat sections 
	$\sE^{\nabla}$. The point is that for any $x \in X$, any vector 
	in the fibre $\sE_x$ extends in a small analytic ball to a unique flat 
	section. Here one uses both that the connection is flat and that one is 
	allowing analytic, and not merely polynomial solutions, to essentially 
	reduce to the case of ordinary differential equations. To go from left to 
	right, one sends a local system $\sL$ to the holomorphic bundle 
	\[
	\sL \underset{\ck}{\otimes} \OO_{X_{\on{an}}},
	\]
	with its connection induced by the tautological one on $\OO_{X_{\on{an}}}$. 
	
	\subsubsection{} 
	To make contact with algebraic D-modules on $X$, note that applying $- 
	\underset{\phantom{ii}\OO_{X}}{\otimes} \OO_{X_{\on{an}}}$, i.e., allowing not only polynomial, but also holomorphic sections, yields a natural 
	analytification functor 
	\begin{equation} \label{e:gaga}
	\{ \text{flat algebraic connections on $X$} \} \rightarrow \{  \text{flat 
	holomorphic connections on $X_{\on{an}}$} \}. 
	\end{equation}
	We emphasize that, on the left, one is considering objects of algebro-geometric nature and, on the right, objects of a complex analytic nature. 
	
	If $X$ is proper then the functor \eqref{e:gaga} is an equivalence. This 
	is a consequence of the so-called GAGA (i.e., g\'{e}om\'etrie alg\'ebrique 
	et g\'eom\'etrie analytique) theorem.
	
	However, if $X$ is non-proper, the functor \eqref{e:gaga} is no longer an equivalence. Namely, it is still essentially surjective, but it typically sends many 
	distinct flat algebraic connections to the same holomorphic connection. The 
	reason is simple and already visible for rank one connections -- there are 
	typically many more holomorphic functions than polynomial functions (unlike 
	in the case of $X$ proper), and the differential equations satisfied by them lead to identifications after analytification between distinct algebraic connections. 
	
	\begin{ex}Suppose $X = \mathbb{A}^1$ with coordinate $z$, and consider the 
	left algebraic D-modules 
		\[
		D_X / (D_X \cdot \partial_z) \quad \text{and} \quad  D_X / (D_X \cdot 
		(\partial_z - 1)).
		\]
		Explicitly, the underlying $O_X$-modules both carry canonical identifications 
		with $O_X$, given by acting on the image of 1 in $D_X$, and carry the 
		connections 
		\[
		\nabla =  \partial_z \quad \text{and} \quad \nabla = \partial_z - 1. 
		\]
		As algebraic D-modules, they are not isomorphic, but their 
		analytifications are. Indeed, this follows from the fact that $e^z$ is 
		a holomorphic, but not polynomial, function and is, up to scalars, the 
		unique solution to 
		\[
		\partial_z u = u, 
		\]
		i.e., a flat section of the second connection. 
	\end{ex}
	
	What rectifies the many-to-oneness of \eqref{e:gaga} is the following 
	nontrivial observation. For a given flat holomorphic bundle 
	$\sE_{\on{an}}$, among the various algebraic connections $\sE$ mapping to 
	it, there will be a unique one whose flat sections do not have `essential 
	singularities at infinity.' More precisely, an algebraic connection $(\sE, 
	\nabla)$ on $X$ is said to be {\em regular} if, for a 
	smooth 
	compactification $$j: X \rightarrow \overline{X}$$ of $X$, with boundary 
	$\overline{X} \setminus X$ a simple normal crossings divisor,
	the connection form $j_* \nabla$ on $j_* \sE$ roughly has at most simple poles on 
	the boundary divisor, cf. \cite{hottatakeuchitanisaki} for a more careful formulation. It turns out this is independent of the choice of $\overline{X}$. 
	
	\begin{ex} In the previous example, we may take the compactification 
		\[
		j: \mathbb{A}^1 \rightarrow \mathbb{P}^1, 
		\]
		i.e., we add the point at infinity. With respect to the coordinate $w = 
		z^{-1}$ near infinity, the two connections take the form 
		\[
		\nabla = \partial_w \quad \text{and} \quad \nabla = \partial_w + w^{-2},
		\]
		respectively. In particular, the second connection has a double pole 
		and so is not regular. This is equivalent to the fact that the function 
		$e^z$ has an essential singularity at $\infty$. 
	\end{ex}
	
	When one restricts to regular connections, the analytification functor 
	\[
	\{ \text{regular algebraic connections on $X$} \} \rightarrow \{  
	\text{flat holomorphic connections on $X_{\on{an}}$} \} 
	\]
	is an equivalence. The resulting composite equivalence with local systems
	\begin{equation} \label{e:rhconn}
	\nabla: \{ \text{regular algebraic connections on $X$} \} \simeq 
	\{\text{finite rank local systems on $X_{\on{an}}$} \}
	\end{equation}
	is explicitly given by taking the flat analytic sections. This is known as the 
	Riemann--Hilbert correspondence for connections.

	\subsubsection{} Finally, let us describe the desired extension of 
	\eqref{e:rhconn} to all constructible sheaves on $X$. Intuitively, we would 
	like to interchange now not only local systems and regular connections, but 
	also those supported on subvarieties. 
	
	A first guess would be that one has an equivalence of abelian categories 
	between constructible sheaves and a suitable category of D-modules. After 
	all, the previous correspondence for connections matched the abelian 
	categories. For example, one might expect to exchange the constant sheaf 
	$\CC_X$ and the regular functions $\OO_X$, and, for any closed point $x$ of 
	$X$, the rank one skyscraper sheaf $\CC_x$ and the delta D-module 
	$\delta_x$.  
	\label{s:cannotlit}
	However, this cannot literally hold. Indeed, if $X$ is of dimension $d > 
	0$, and one looks at the homomorphisms in the abelian category, one has 
	\begin{align*}
	\on{Hom}_{\on{Sh}(X_{\on{an}})^\heartsuit}( \CC_X, \CC_x) \simeq \CC \quad 
	&\text{but} \quad \on{Hom}_{\on{D-mod}(X)^\heartsuit}( \OO_X, \delta_x) 
	\simeq 0.  
	\intertext{To correct for this, observe that when passing to derived 
	categories, we have }
	\on{Hom}_{\on{Sh}(X_{\on{an}})}( \CC_X, \CC_x) \simeq \CC \quad &\text{and} 
	\quad \on{Hom}_{\on{D-mod}(X)}( \OO_X, \delta_x) \simeq \CC[-d].  
	\end{align*}
	The reader is encouraged to check both assertions, e.g., for the affine line $X = 
	\mathbb{A}^1$. In any case, we see the desired extension will need to 
	involve cohomological shifts and, in particular, must be formulated in terms 
	of derived categories. 
	
	\subsubsection{}   With this in mind, we are led to the following. The 
	category of regular holonomic D-modules, denoted by 
	$\on{D-mod}^{\on{rh}}(X)$, is the full triangulated subcategory of 
	$\on{D-mod}(X)$ generated by $*$-extensions of regular algebraic 
	connections from smooth subvarieties $Z \subset X$.  
	
	One has a functor of (derived) flat analytic sections, defined on the 
	category of all D-modules 
	\begin{equation} \label{e:fllat}
	\nabla:  \on{D-mod}(X) \rightarrow \{ \text{sheaves of $\CC$-vector spaces} 
	\},
	\end{equation}
	which we shall describe in more detail momentarily. However, the main 
	result here is that $\eqref{e:fllat}$, when restricted to regular holonomic 
	D-modules, factors through the constructible derived category and yields 
	an equivalence 
	\begin{equation}
	\nabla: \on{D-mod}^{\on{rh}}(X) \simeq \on{Sh}(X_{\on{an}})^b. 
	\end{equation}
	This is known as the Riemann--Hilbert correspondence for D-modules.

	\subsubsection{}  Invariantly, \eqref{e:fllat} is defined as follows. The 
	holomorphic canonical bundle $\omega_{X_{\on{an}}}$ is naturally a right 
	D-module, and \eqref{e:fllat} is obtained by tensoring with it:
	\[
	\nabla := \omega_{X_{\on{an}}} \underset{\phantom{ii}\sD_X}{\otimes} -. 
	\]
	As with Lie algebra homology, a canonical projective resolution of 
	$\omega_{X_{\on{an}}}$ yields the following presentation. Let us denote by 
	$\Omega^i_{\on{an}}$, for $1 \leqslant i \leqslant d$, the sheaves of 
	holomorphic $i$-forms on $X$. In particular, note that $$\Omega^d_{\on{an}} 
	= \omega_{X_{\on{an}}}.$$ With this, \eqref{e:fllat} sends a D-module 
	$\sM$, concentrated in cohomological degree zero, to the de Rham complex of 
	sheaves
	\begin{equation} \label{e:flatsec}
	\sM \rightarrow \sM \underset{\phantom{ii}\OO_X} \otimes \Omega^1_{\on{an}} \rightarrow 
	\cdots \rightarrow \sM \underset{\phantom{ii}\OO_X}\otimes \Omega^d_{\on{an}}
	\end{equation}
	where $\sM$ is placed in cohomological degree $-d$. Explicitly, for a local 
	section $m$ of $\sM$, and local coordinates $z_1, \ldots, z_d$ on $X$, the 
	first differential sends 
	\[
	m \mapsto \underset{i}{\Sigma} \hspace{1mm} (\partial_i \cdot m) \otimes 
	dz_i,
	\]
	and similarly for higher differentials. A complex of D-modules is then sent 
	to the totalization of a corresponding de Rham bicomplex. This construction 
	visibly respects quasi-isomorphisms and yields the desired functor 
	\eqref{e:fllat}.

	\subsubsection{} Let us conclude this discussion with some basic examples 
	and properties. 
	
	\begin{ex} The Riemann--Hilbert correspondence sends the structure sheaf to 
	a shift of the constant sheaf, namely
		\begin{align*}
		\nabla(\OO_X) &\simeq \CC_X[d_X].
		\intertext{It sends a delta D-module $\delta_x$ to the corresponding 
		skyscraper sheaf, i.e., }
		\nabla(\delta_x) &\simeq \CC_x. 
		\intertext{The reader may wish to check these directly and note in 
		particular their compatibility with the discussion of Section 
		\ref{s:cannotlit}.}
		\intertext{\hspace{1em} More generally, for a smooth subvariety $i: Z 
		\rightarrow X$ of dimension $d_Z$ and a regular connection $\sE$ on 
		$Z$ with flat analytic sections $\sL$, we have }
		\nabla( i_{*} \sE) &\simeq i_* \sL [d_Z]. 
		\end{align*}
		\label{ex:exrhs}
	\end{ex}
	
	\subsubsection{}  It is also possible to treat D-modules on a singular 
	variety $X$. For any embedding $X \hookrightarrow S$ into a smooth variety, 
	this category canonically agrees with the full subcategory of D-modules on $S$ 
	set-theoretically supported on $X$. The Riemann--Hilbert correspondence 
	again holds for singular varieties and may be deduced from the smooth 
	case.

	The formulas of Example \ref{ex:exrhs} are then a special case of the 
	following. For any map of algebraic varieties $X \rightarrow Y$, the 
	Riemann--Hilbert correspondence exchanges $*$-pushforwards and 
	$!$-pullbacks of D-modules and constructible sheaves.

	\subsubsection{} Finally, let us discuss abelian categories. The 
	subcategory $\on{D-mod}^{\on{rh}}(X)$ is closed under the standard 
	truncation functors on $\on{D-mod}(X)$ and canonically identifies with the 
	derived category of the corresponding abelian category of regular holonomic 
	D-modules. From Example \ref{ex:exrhs}, we see that this abelian category is not exchanged with 
	the standard abelian category of constructible sheaves. Instead, it is 
	identified with the abelian category of perverse sheaves.

	\subsubsection{} For further details, the reader may wish to look at \cite{arkhipovnotes}, \cite{bernsteinnotes}, \cite{hottatakeuchitanisaki}.

		\vspace{.5cm}
	
	\noindent {\bf Conflict of interest.} There is no conflict of interest.

	\bibliographystyle{amsalpha}
	\bibliography{samplez}

\newcommand{\etalchar}[1]{$^{#1}$}
\providecommand{\bysame}{\leavevmode\hbox to3em{\hrulefill}\thinspace}
\providecommand{\MR}{\relax\ifhmode\unskip\space\fi MR }
\providecommand{\MRhref}[2]{%
  \href{http://www.ams.org/mathscinet-getitem?mr=#1}{#2}
}
\providecommand{\href}[2]{#2}
\begin{thebibliography}{AGK{\etalchar{+}}20b}

\bibitem[AB00]{avramovbuchweitz}
Luchezar Avramov and Ragnar-Olaf Buchweitz, \emph{Support varieties and
  cohomology over complete intersections}, Invent. Math. \textbf{142} (2000),
  no.~2, 285--318. \MR{1794064}

\bibitem[AB09]{ab}
Sergey Arkhipov and Roman Bezrukavnikov, \emph{Perverse sheaves on affine flags
  and {L}anglands dual group}, Israel J. Math. \textbf{170} (2009), 135--183,
  With an appendix by Bezrukavnikov and Ivan Mirkovi\'{c}. \MR{2506322}

\bibitem[ABC{\etalchar{+}}18]{paris-notes}
Dima Arinkin, Dario Beraldo, Justin Campbell, Lin Chen, Yuchen Fu, Dennis
  Gaitsgory, Quoc Ho, Sergey Lysenko, Sam Raskin, Simon Riche, Nick Rozenblyum,
  James Tao, David Yang, and Yifei Zhao, \emph{{Notes from the winter school on
  local geometric {L}anglands}}, 2018, Available at
  \url{https://sites.google.com/site/winterlanglands2018/notes-of-talks}.

\bibitem[ABG04]{abg}
Sergey Arkhipov, Roman Bezrukavnikov, and Victor Ginzburg, \emph{Quantum
  groups, the loop {G}rassmannian, and the {S}pringer resolution}, J. Amer.
  Math. Soc. \textbf{17} (2004), no.~3, 595--678. \MR{2053952}

\bibitem[AF12]{arakawafiebig}
Tomoyuki Arakawa and Peter Fiebig, \emph{On the restricted {V}erma modules at
  the critical level}, Trans. Amer. Math. Soc. \textbf{364} (2012), no.~9,
  4683--4712. \MR{2922606}

\bibitem[AG15a]{arinkingaitssingsupp}
Dima Arinkin and Dennis Gaitsgory, \emph{Singular support of coherent sheaves
  and the geometric {L}anglands conjecture}, Selecta Math. (N.S.) \textbf{21}
  (2015), no.~1, 1--199. \MR{3300415}

\bibitem[AG15b]{ag}
Sergey Arkhipov and Dennis Gaitsgory, \emph{Localization and the long
  intertwining operator for representations of affine {K}ac-{M}oody algebras},
  Preprint \textbf{\url{http://www. math. harvard. edu/~gaitsgde/GL/Arkh. pdf}}
  (2015).

\bibitem[AGK{\etalchar{+}}20a]{restrictedvar2}
Dima Arinkin, Dennis Gaitsgory, David Kazhdan, Sam Raskin, Nick Rozenblyum, and
  Yakov Varshavsky, \emph{Duality for automorphic sheaves with nilpotent
  singular support}, Preprint, arXiv:2012.07665 (2020).

\bibitem[AGK{\etalchar{+}}20b]{restrictedvar1}
\bysame, \emph{The stack of local systems with restricted variation and
  geometric {L}anglands theory with nilpotent singular support}, Preprint,
  arXiv:2010.01906 (2020).

\bibitem[AGK{\etalchar{+}}21]{restrictedvar3}
\bysame, \emph{Automorphic functions as the trace of {F}robenius}, Preprint,
  arXiv:2102.07906 (2021).

\bibitem[AN12]{arkhipovnotes}
Sergey Arkhipov and Nikita Nikolaev, \emph{D-modules}, Notes, available at
  \url{https://www.yumpu.com/en/document/view/19420256/d-modules-lecture-notes-wiki-university-of-toronto}
  (2012).

\bibitem[AR15]{acharrider}
Pramod Achar and Laura Rider, \emph{Parity sheaves on the affine {G}rassmannian
  and the {M}irkovi\'{c}-{V}ilonen conjecture}, Acta Math. \textbf{215} (2015),
  no.~2, 183--216. \MR{3455233}

\bibitem[Ara17]{araintro}
Tomoyuki Arakawa, \emph{Introduction to {W}-algebras and their representation
  theory}, Perspectives in {L}ie theory, Springer INdAM Ser., vol.~19,
  Springer, Cham, 2017, pp.~179--250. \MR{3751125}

\bibitem[Ari18]{arinkinindcohcats}
Dima Arinkin, \emph{Day {III}, {T}alk 3: {S}pectral side in the classical
  case},
  \url{https://lysenko.perso.math.cnrs.fr/Notes_talks_winter2018/GL-6(Dima).pdf}
  (2018).

\bibitem[Ark04]{ark}
Sergey Arkhipov, \emph{Algebraic construction of contragradient quasi-{V}erma
  modules in positive characteristic}, Representation theory of algebraic
  groups and quantum groups, Adv. Stud. Pure Math., vol.~40, Math. Soc. Japan,
  Tokyo, 2004, pp.~27--68. \MR{2074588}

\bibitem[BB81]{beilinson-bernstein81}
Alexandre Be\u{\i}linson and Joseph Bernstein, \emph{Localisation de
  {$g$}-modules}, C. R. Acad. Sci. Paris S\'{e}r. I Math. \textbf{292} (1981),
  no.~1, 15--18. \MR{610137}

\bibitem[BB83]{bb82}
Alexander Be\u{\i}linson and Joseph Bernstein, \emph{A generalization of
  {C}asselman's submodule theorem}, Representation theory of reductive groups
  ({P}ark {C}ity, {U}tah, 1982), Progr. Math., vol.~40, Birkh\"{a}user Boston,
  Boston, MA, 1983, pp.~35--52. \MR{733805}

\bibitem[BD]{bdh}
Alexander Be\u{\i}linson and Vladimir Drinfeld, \emph{Quantization of
  {H}itchin's integrable system and {H}ecke eigensheaves}, Preprint
  \url{http://math.uchicago.edu/~drinfeld/langlands/QuantizationHitchin.pdf}.

\bibitem[BD04]{BD}
\bysame, \emph{Chiral algebras}, American Mathematical Society Colloquium
  Publications, vol.~51, American Mathematical Society, Providence, RI, 2004.
  \MR{2058353}

\bibitem[Ber82]{bernsteinnotes}
Joseph Bernstein, \emph{Algebraic theory of {D}-modules}, Notes, available at
  \url{http://www.math.tau.ac.il/~bernstei/Unpublished_texts/unpublished_texts/Bernstein82-dmod-from-chicago.pdf}
  (1982).

\bibitem[Ber17a]{Bearaldo}
Dario Beraldo, \emph{Loop group actions on categories and {W}hittaker
  invariants}, Advances in Mathematics \textbf{322} (2017), 565--636.

\bibitem[Ber17b]{beraldo}
\bysame, \emph{Loop group actions on categories and {W}hittaker invariants},
  Adv. Math. \textbf{322} (2017), 565--636. \MR{3720805}

\bibitem[Bez16]{bez2realz}
Roman Bezrukavnikov, \emph{On two geometric realizations of an affine {H}ecke
  algebra}, Publ. Math. Inst. Hautes \'{E}tudes Sci. \textbf{123} (2016),
  1--67. \MR{3502096}

\bibitem[BF08]{bezfinksat}
Roman Bezrukavnikov and Michael Finkelberg, \emph{Equivariant {S}atake category
  and {K}ostant-{W}hittaker reduction}, Mosc. Math. J. \textbf{8} (2008),
  no.~1, 39--72, 183. \MR{2422266}

\bibitem[BF19]{bravermanfinkelbergconjecturecoloumb}
Alexander Braverman and Michael Finkelberg, \emph{Coulomb branches of
  3-dimensional gauge theories and related structures}, Geometric
  representation theory and gauge theory, Lecture Notes in Math., vol. 2248,
  Springer, Cham, [2019] \copyright 2019, pp.~1--52. \MR{4286060}

\bibitem[BFGT21]{bfgt}
Alexander Braverman, Michael Finkelberg, Victor Ginzburg, and Roman Travkin,
  \emph{Mirabolic {S}atake equivalence and supergroups}, Compos. Math.
  \textbf{157} (2021), no.~8, 1724--1765. \MR{4292176}

\bibitem[BFN18]{bfn}
Alexander Braverman, Michael Finkelberg, and Hiraku Nakajima, \emph{Towards a
  mathematical definition of {C}oulomb branches of 3-dimensional
  {$\mathcal{N}=4$} gauge theories, {II}}, Adv. Theor. Math. Phys. \textbf{22}
  (2018), no.~5, 1071--1147. \MR{3952347}

\bibitem[BG80]{bege1980}
Joseph Bernstein and Sergey Gel'fand, \emph{Tensor products of finite- and
  infinite-dimensional representations of semisimple {L}ie algebras},
  Compositio Math. \textbf{41} (1980), no.~2, 245--285. \MR{581584}

\bibitem[BG99]{beilinsonginzburg}
Alexander Be\u{\i}linson and Victor Ginzburg, \emph{Wall-crossing functors and
  {D}-modules}, Representation Theory of the American Mathematical Society
  \textbf{3} (1999), no.~1, 1--31.

\bibitem[BGM{\etalchar{+}}19]{bezgaitsmiriwahotiwhittaker}
Roman Bezrukavnikov, Dennis Gaitsgory, Ivan Mirkovi\'{c}, Simon Riche, and
  Laura Rider, \emph{An {I}wahori-{W}hittaker model for the {S}atake category},
  J. \'{E}c. polytech. Math. \textbf{6} (2019), 707--735. \MR{4014634}

\bibitem[BIK08]{bensoniyengarkrause}
Dave Benson, Srikanth Iyengar, and Henning Krause, \emph{Local cohomology and
  support for triangulated categories}, Ann. Sci. \'{E}c. Norm. Sup\'{e}r. (4)
  \textbf{41} (2008), no.~4, 573--619. \MR{2489634}

\bibitem[BK81]{brylin81}
Jean-Luc Brylinski and Masaki Kashiwara, \emph{Kazhdan-{L}usztig conjecture and
  holonomic systems}, Invent. Math. \textbf{64} (1981), no.~3, 387--410.
  \MR{632980}

\bibitem[BK15]{krem}
Erik Backelin and Kobi Kremnizer, \emph{Singular localization of
  $\mathfrak{g}$-modules and applications to representation theory}, Journal of
  the European Mathematical Society \textbf{17} (2015), no.~11, 2763--2787.

\bibitem[BMR06]{bmr}
Roman Bezrukavnikov, Ivan Mirkovi\'{c}, and Dmitriy Rumynin, \emph{Singular
  localization and intertwining functors for reductive {L}ie algebras in prime
  characteristic}, Nagoya Math. J. \textbf{184} (2006), 1--55. \MR{2285230}

\bibitem[BR18]{bezrichetop}
Roman Bezrukavnikov and Simon Riche, \emph{A topological approach to {S}oergel
  theory}, Preprint, arXiv:1807.07614 (2018).

\bibitem[BRar]{bezricheprinblock}
\bysame, \emph{Hecke action on the principal block}, Compositio Mathematica
  (2021, to appear).

\bibitem[BRR20]{bezricheriderregcentralizer}
Roman Bezrukavnikov, Simon Riche, and Laura Rider, \emph{Modular affine {H}ecke
  category and regular unipotent centralizer, {I}}, Preprint, arXiv:2005.05583
  (2020).

\bibitem[BS13]{sage1}
Christopher Bremer and Daniel Sage, \emph{Moduli spaces of irregular singular
  connections}, Int. Math. Res. Not. IMRN (2013), no.~8, 1800--1872.
  \MR{3047490}

\bibitem[BS18]{sage3}
\bysame, \emph{A theory of minimal {$K$}-types for flat {$G$}-bundles}, Int.
  Math. Res. Not. IMRN (2018), no.~11, 3507--3555. \MR{3810225}

\bibitem[BZFN10]{benzvifrancisnadler}
David Ben-Zvi, John Francis, and David Nadler, \emph{Integral transforms and
  {D}rinfeld centers in derived algebraic geometry}, J. Amer. Math. Soc.
  \textbf{23} (2010), no.~4, 909--966. \MR{2669705}

\bibitem[BZGO18]{BZO}
David Ben-Zvi, Sam Gunningham, and Hendrik Orem, \emph{{Highest Weights for
  Categorical Representations}}, International Mathematics Research Notices
  (2018), rny258.

\bibitem[BZN09]{bznadlerchartheory}
David Ben-Zvi and David Nadler, \emph{The character theory of a complex group},
  arXiv preprint arXiv:0904.1247 (2009).

\bibitem[BZN18]{benzvinadlerbetti}
\bysame, \emph{Betti geometric {L}anglands}, Algebraic geometry: {S}alt {L}ake
  {C}ity 2015, Proc. Sympos. Pure Math., vol.~97, Amer. Math. Soc., Providence,
  RI, 2018, pp.~3--41. \MR{3821166}

\bibitem[Cas80]{cass80spherical}
William Casselman, \emph{The unramified principal series of
  {${\mathfrak{p}}$}-adic groups. {I}. {T}he spherical function}, Compositio
  Math. \textbf{40} (1980), no.~3, 387--406. \MR{571057}

\bibitem[CD21]{ahc}
Justin Campbell and Gurbir Dhillon, \emph{Steinberg--{W}hittaker localization
  and affine {H}arish--{C}handra bimodules}, Preprint, arXiv:2108.02806 (2021).

\bibitem[CG10]{chrissginzburg}
Neil Chriss and Victor Ginzburg, \emph{Representation theory and complex
  geometry}, Modern Birkh\"{a}user Classics, Birkh\"{a}user Boston, Ltd.,
  Boston, MA, 2010, Reprint of the 1997 edition. \MR{2838836}

\bibitem[CG19]{costellogaiotto}
Kevin Costello and Davide Gaiotto, \emph{Vertex {O}perator {A}lgebras and 3d
  {$\mathcal N=4$} gauge theories}, J. High Energy Phys. (2019), no.~5, 018,
  37. \MR{3976868}

\bibitem[CG20]{creutziggaiotto}
Thomas Creutzig and Davide Gaiotto, \emph{Vertex algebras for {S}-duality},
  Comm. Math. Phys. \textbf{379} (2020), no.~3, 785--845. \MR{4163353}

\bibitem[CGL20]{creutziggaiottolinshaw}
Thomas Creutzig, Davide Gaiotto, and Andrew Linshaw, \emph{S-duality for the
  large {$N=4$} superconformal algebra}, Comm. Math. Phys. \textbf{374} (2020),
  no.~3, 1787--1808. \MR{4076088}

\bibitem[CK17]{chenkamgarpour}
Tsao-Hsien Chen and Masoud Kamgarpour, \emph{Preservation of depth in the local
  geometric {L}anglands correspondence}, Trans. Amer. Math. Soc. \textbf{369}
  (2017), no.~2, 1345--1364. \MR{3572276}

\bibitem[Cou95]{coutinho}
Severino~Collier Coutinho, \emph{A primer of algebraic {$D$}-modules}, London
  Mathematical Society Student Texts, vol.~33, Cambridge University Press,
  Cambridge, 1995. \MR{1356713}

\bibitem[CS80]{cassshalika}
William Casselman and Joseph Shalika, \emph{The unramified principal series of
  {$p$}-adic groups. {II}. {T}he {W}hittaker function}, Compositio Math.
  \textbf{41} (1980), no.~2, 207--231. \MR{581582}

\bibitem[CW21]{williamsonsurvey2}
Joshua Ciappara and Geordie Williamson, \emph{Lectures on the geometry and
  modular representation theory of algebraic groups}, J. Aust. Math. Soc.
  \textbf{110} (2021), no.~1, 1--47. \MR{4211934}

\bibitem[Dhi21]{lpw}
Gurbir Dhillon, \emph{Semi-infinite cohomology and the linkage principle for
  {$\mathscr{W}$}-algebras}, Advances in Mathematics \textbf{381} (2021),
  107625, 62. \MR{4215746}

\bibitem[DR20]{locw}
Gurbir Dhillon and Sam Raskin, \emph{Localization for affine
  {$\mathscr{W}$}-algebras}, Preprint, arXiv:2010.11434 (2020).

\bibitem[DS84]{drinfeldsokolov}
Vladimir Drinfel'd and Vladimir Sokolov, \emph{Lie algebras and equations of
  {K}orteweg-de {V}ries type}, Current problems in mathematics, {V}ol. 24,
  Itogi Nauki i Tekhniki, Akad. Nauk SSSR, Vsesoyuz. Inst. Nauchn. i Tekhn.
  Inform., Moscow, 1984, pp.~81--180. \MR{760998}

\bibitem[DY]{DhillonYang}
Gurbir Dhillon and David Yang, \emph{Characters of the critical level highest
  weight modules for affine {L}ie algebras}, In preparation.

\bibitem[EH14]{emertonhelm}
Matthew Emerton and David Helm, \emph{The local {L}anglands correspondence for
  {${\rm GL}_n$} in families}, Ann. Sci. \'{E}c. Norm. Sup\'{e}r. (4)
  \textbf{47} (2014), no.~4, 655--722. \MR{3250061}

\bibitem[Eis80]{eisenbud}
David Eisenbud, \emph{Homological algebra on a complete intersection, with an
  application to group representations}, Trans. Amer. Math. Soc. \textbf{260}
  (1980), no.~1, 35--64. \MR{570778}

\bibitem[Enr79]{enright79}
Thomas Enright, \emph{On the fundamental series of a real semisimple {L}ie
  algebra: their irreducibility, resolutions and multiplicity formulae}, Ann.
  of Math. (2) \textbf{110} (1979), no.~1, 1--82. \MR{541329}

\bibitem[FBZ04]{fbz}
Edward Frenkel and David Ben-Zvi, \emph{Vertex algebras and algebraic curves},
  second ed., Mathematical Surveys and Monographs, vol.~88, American
  Mathematical Society, Providence, RI, 2004. \MR{2082709}

\bibitem[FF91]{feiginfrenkel91}
Boris Feigin and Edward Frenkel, \emph{Duality in {$W$}-algebras}, Internat.
  Math. Res. Notices (1991), no.~6, 75--82. \MR{1136408}

\bibitem[FG04]{frenkelgaitsgorylocgrassmannian04}
Edward Frenkel and Dennis Gaitsgory, \emph{{$D$}-modules on the affine
  {G}rassmannian and representations of affine {K}ac-{M}oody algebras}, Duke
  Math. J. \textbf{125} (2004), no.~2, 279--327. \MR{2096675}

\bibitem[FG06a]{frenkgaitslocallangds}
\bysame, \emph{Local geometric {L}anglands correspondence and affine
  {K}ac-{M}oody algebras}, Algebraic geometry and number theory, Progr. Math.,
  vol. 253, Birkh\"{a}user Boston, Boston, MA, 2006, pp.~69--260. \MR{2263193}

\bibitem[FG06b]{fg06loc}
\bysame, \emph{Local geometric {L}anglands correspondence and affine
  {K}ac-{M}oody algebras}, Algebraic geometry and number theory, Progr. Math.,
  vol. 253, Birkh\"{a}user Boston, Boston, MA, 2006, pp.~69--260. \MR{2263193}

\bibitem[FG09a]{fglocaff}
\bysame, \emph{{$D$}-modules on the affine flag variety and representations of
  affine {K}ac-{M}oody algebras}, Represent. Theory \textbf{13} (2009),
  470--608. \MR{2558786}

\bibitem[FG09b]{fg-sph}
\bysame, \emph{Local geometric {L}anglands correspondence: the spherical case},
  Algebraic analysis and around, Adv. Stud. Pure Math., vol.~54, Math. Soc.
  Japan, Tokyo, 2009, pp.~167--186.

\bibitem[FG09c]{fgloc09}
\bysame, \emph{Localization of {$\mathfrak{g}$}-modules on the affine
  {G}rassmannian}, Ann. of Math. (2) \textbf{170} (2009), no.~3, 1339--1381.
  \MR{2600875}

\bibitem[FG20]{frenkelgaiotto}
Edward Frenkel and Davide Gaiotto, \emph{Quantum {L}anglands dualities of
  boundary conditions, {$D$}-modules, and conformal blocks}, Commun. Number
  Theory Phys. \textbf{14} (2020), no.~2, 199--313. \MR{4084137}

\bibitem[FGV01]{fgv}
Edward Frenkel, Dennis Gaitsgory, and Kari Vilonen, \emph{Whittaker patterns in
  the geometry of moduli spaces of bundles on curves}, Ann. of Math. (2)
  \textbf{153} (2001), no.~3, 699--748. \MR{1836286}

\bibitem[FK88]{freitagbook}
Eberhard Freitag and Reinhardt Kiehl, \emph{\'{E}tale cohomology and the {W}eil
  conjecture}, Ergebnisse der Mathematik und ihrer Grenzgebiete (3) [Results in
  Mathematics and Related Areas (3)], vol.~13, Springer-Verlag, Berlin, 1988,
  Translated from the German by Betty S. Waterhouse and William C. Waterhouse,
  With an historical introduction by J. A. Dieudonn\'{e}. \MR{926276}

\bibitem[FM97]{frenkelmalikov1}
Igor~B. Frenkel and Feodor~G. Malikov, \emph{Kazhdan-{L}usztig tensoring and
  {H}arish-{C}handra categories}, arXiv preprint q-alg/9703010 (1997).

\bibitem[FM99]{finkelbergmirkovicconjecture}
Michael Finkelberg and Ivan Mirkovi\'{c}, \emph{Semi-infinite flags. {I}.
  {C}ase of global curve {$\bold P^1$}}, Differential topology,
  infinite-dimensional {L}ie algebras, and applications, Amer. Math. Soc.
  Transl. Ser. 2, vol. 194, Amer. Math. Soc., Providence, RI, 1999,
  pp.~81--112. \MR{1729360}

\bibitem[Fre07a]{frenkelloop}
Edward Frenkel, \emph{Langlands correspondence for loop groups}, Cambridge
  Studies in Advanced Mathematics, vol. 103, Cambridge University Press,
  Cambridge, 2007. \MR{2332156}

\bibitem[Fre07b]{frenkelglobalnotes}
\bysame, \emph{Lectures on the {L}anglands program and conformal field theory},
  Frontiers in number theory, physics, and geometry. {II}, Springer, Berlin,
  2007, pp.~387--533. \MR{2290768}

\bibitem[Fre08]{frenkellocalnotes}
\bysame, \emph{Ramifications of the geometric {L}anglands program},
  Representation theory and complex analysis, Lecture Notes in Math., vol.
  1931, Springer, Berlin, 2008, pp.~51--135. \MR{2409698}

\bibitem[FS21]{farguescholze}
Laurent Fargues and Peter Scholze, \emph{Geometrization of the local
  {L}anglands correspondence}, Preprint, arXiv:2102.13459 (2021).

\bibitem[FW08]{frenkelwitten}
Edward Frenkel and Edward Witten, \emph{Geometric endoscopy and mirror
  symmetry}, Commun. Number Theory Phys. \textbf{2} (2008), no.~1, 113--283.
  \MR{2417848}

\bibitem[FZ10]{frenkelzhuoper}
Edward Frenkel and Xinwen Zhu, \emph{Any flat bundle on a punctured disc has an
  oper structure}, Math. Res. Lett. \textbf{17} (2010), no.~1, 27--37.
  \MR{2592725}

\bibitem[Gai05]{gaitsgorynotes}
Dennis Gaitsgory, \emph{Geometric representation theory}, Notes, available at
  \url{https://people.math.harvard.edu/~gaitsgde/267y/catO.pdf} (2005).

\bibitem[Gai07]{quantum-langlands-summary}
\bysame, \emph{Quantum {L}anglands correspondence}, arXiv preprint
  arXiv:1601.05279 (2007).

\bibitem[Gai13]{gaitsgoryindcoh}
\bysame, \emph{{I}nd-coherent sheaves}, Mosc. Math. J. \textbf{13} (2013),
  no.~3, 399--528, 553. \MR{3136100}

\bibitem[Gai15]{1affine}
\bysame, \emph{Sheaves of categories and the notion of 1-affineness}, Stacks
  and categories in geometry, topology, and algebra, Contemp. Math., vol. 643,
  Amer. Math. Soc., Providence, RI, 2015, pp.~127--225. \MR{3381473}

\bibitem[Gai16]{Gaitsgoryinventshtukas}
\bysame, \emph{From geometric to function-theoretic {L}anglands (or how to
  invent shtukas)}, Preprint, arXiv:1606.09608 (2016).

\bibitem[Gai18a]{gaiotto16s-duality}
Davide Gaiotto, \emph{S-duality and boundary conditions and the geometric
  {L}anglands program}, String-{M}ath 2016, Proc. Sympos. Pure Math., vol.~98,
  Amer. Math. Soc., Providence, RI, 2018, pp.~139--179. \MR{3821752}

\bibitem[Gai18b]{winter-langlands-summary}
Dennis Gaitsgory, \emph{Winter school on local geometric {L}anglands theory:
  {P}rogram}, Notes, available at
  \url{http://www.iecl.univ-lorraine.fr/~Sergey.Lysenko/program_1.pdf} (2018).

\bibitem[Gai19]{gaiottotwistcompact}
Davide Gaiotto, \emph{Twisted compactifications of {$3{\rm d}\ \mathcal N=4$}
  theories and conformal blocks}, J. High Energy Phys. (2019), no.~2, 061,
  front matter+41. \MR{3933138}

\bibitem[Gai20]{whitlocglob}
Dennis Gaitsgory, \emph{The local and global versions of the {W}hittaker
  category}, Pure Appl. Math. Q. \textbf{16} (2020), no.~3, 775--904.
  \MR{4176538}

\bibitem[Gin95]{ginzburgsatake}
Victor Ginzburg, \emph{Perverse sheaves on a {L}oop group and {L}anglands'
  duality}, Preprint, arXiv:9511007 (1995).

\bibitem[GKRV21]{gkrvtoyshtuka}
Dennis Gaitsgory, David Kazhdan, Nick Rozenblyum, and Yakov Varshavsky, \emph{A
  toy model for the {D}rinfeld--{L}afforgue shtuka construction}, Preprint,
  arXiv:1908.05420 (2021).

\bibitem[GR17]{gaitsroz}
Dennis Gaitsgory and Nick Rozenblyum, \emph{A study in derived algebraic
  geometry. {V}ol. {I}. {C}orrespondences and duality}, Mathematical Surveys
  and Monographs, vol. 221, American Mathematical Society, Providence, RI,
  2017. \MR{3701352}

\bibitem[GR19]{gaiottorapcak}
Davide Gaiotto and Miroslav Rap\v{c}\'{a}k, \emph{Vertex algebras at the
  corner}, J. High Energy Phys. (2019), no.~1, 160, front matter+85.
  \MR{3919335}

\bibitem[Gul74]{gulliksen}
Tor Gulliksen, \emph{A change of ring theorem with applications to
  {P}oincar\'{e} series and intersection multiplicity}, Math. Scand.
  \textbf{34} (1974), 167--183. \MR{364232}

\bibitem[GW08]{gukovwitten}
Sergei Gukov and Edward Witten, \emph{Gauge theory, ramification, and the
  geometric {L}anglands program}, Current developments in mathematics, 2006,
  Int. Press, Somerville, MA, 2008, pp.~35--180. \MR{2459305}

\bibitem[GW09]{gaiottowitten}
Davide Gaiotto and Edward Witten, \emph{{$S$}-duality of boundary conditions in
  {$\mathscr N=4$} super {Y}ang-{M}ills theory}, Adv. Theor. Math. Phys.
  \textbf{13} (2009), no.~3, 721--896. \MR{2610576}

\bibitem[GY]{gaiottoyu}
Davide Gaiotto and Philsang Yoo, \emph{Gauge theory and the local {G}eometric
  {L}anglands program}, In preparation.

\bibitem[Hel21]{hellmann}
Eugene Hellmann, \emph{On the derived category of the {I}wahori--{H}ecke
  algebra}, Preprint, arXiv:2006.03013 (2021).

\bibitem[HR21]{hilburnraskin}
Justin Hilburn and Sam Raskin, \emph{Tate's thesis in the de {R}ham setting},
  Preprint, arXiv:2107.11325 (2021).

\bibitem[HTT08]{hottatakeuchitanisaki}
Ryoshi Hotta, Kiyoshi Takeuchi, and Toshiyuki Tanisaki, \emph{{$D$}-modules,
  perverse sheaves, and representation theory}, Progress in Mathematics, vol.
  236, Birkh\"{a}user Boston, Inc., Boston, MA, 2008, Translated from the 1995
  Japanese edition by Takeuchi. \MR{2357361}

\bibitem[Jan79]{jantz79}
Jens~Carsten Jantzen, \emph{Moduln mit einem h\"{o}chsten {G}ewicht}, Lecture
  Notes in Mathematics, vol. 750, Springer, Berlin, 1979. \MR{552943}

\bibitem[Kas89]{kashloc}
Masaki Kashiwara, \emph{Representation theory and {$D$}-modules on flag
  varieties}, no. 173-174, 1989, Orbites unipotentes et repr\'{e}sentations,
  III, pp.~9, 55--109. \MR{1021510}

\bibitem[Kas90]{kashiwara90}
\bysame, \emph{Kazhdan-{L}usztig conjecture for a symmetrizable {K}ac-{M}oody
  {L}ie algebra}, The {G}rothendieck {F}estschrift, {V}ol. {II}, Progr. Math.,
  vol.~87, Birkh\"{a}user Boston, Boston, MA, 1990, pp.~407--433. \MR{1106905}

\bibitem[Kas03]{kashiwarabookd-modules}
\bysame, \emph{{$D$}-modules and microlocal calculus}, Translations of
  Mathematical Monographs, vol. 217, American Mathematical Society, Providence,
  RI, 2003, Translated from the 2000 Japanese original by Mutsumi Saito,
  Iwanami Series in Modern Mathematics. \MR{1943036}

\bibitem[Kel94]{kellerderiving}
Bernhard Keller, \emph{Deriving {DG} categories}, Ann. Sci. \'{E}cole Norm.
  Sup. (4) \textbf{27} (1994), no.~1, 63--102. \MR{1258406}

\bibitem[KL79]{kl79}
David Kazhdan and George Lusztig, \emph{Representations of {C}oxeter groups and
  {H}ecke algebras}, Invent. Math. \textbf{53} (1979), no.~2, 165--184.
  \MR{560412}

\bibitem[KL87]{kazhdanlusztig87}
\bysame, \emph{Proof of the {D}eligne-{L}anglands conjecture for {H}ecke
  algebras}, Invent. Math. \textbf{87} (1987), no.~1, 153--215. \MR{862716}

\bibitem[Kos63]{kostantslice}
Bertram Kostant, \emph{Lie group representations on polynomial rings}, Amer. J.
  Math. \textbf{85} (1963), 327--404. \MR{158024}

\bibitem[Kra05]{krause}
Henning Krause, \emph{The stable derived category of a {N}oetherian scheme},
  Compos. Math. \textbf{141} (2005), no.~5, 1128--1162. \MR{2157133}

\bibitem[KS15]{kamgarpourschedler}
Masoud Kamgarpour and Travis Schedler, \emph{Geometrization of principal series
  representations of reductive groups}, Ann. Inst. Fourier (Grenoble)
  \textbf{65} (2015), no.~5, 2273--2330. \MR{3449212}

\bibitem[KT90]{kashiwaratanisakiposlevelii}
Masaki Kashiwara and Toshiyuki Tanisaki, \emph{Kazhdan-{L}usztig conjecture for
  symmetrizable {K}ac-{M}oody {L}ie algebra. {II}. {I}ntersection cohomologies
  of {S}chubert varieties}, Operator algebras, unitary representations,
  enveloping algebras, and invariant theory ({P}aris, 1989), Progr. Math.,
  vol.~92, Birkh\"{a}user Boston, Boston, MA, 1990, pp.~159--195. \MR{1103590}

\bibitem[KT95]{kashiwaratanisakinegativeintegrallevel}
\bysame, \emph{Kazhdan-{L}usztig conjecture for affine {L}ie algebras with
  negative level}, Duke Math. J. \textbf{77} (1995), no.~1, 21--62.
  \MR{1317626}

\bibitem[KT96]{kashiwaratanisakinegativerationallevel}
\bysame, \emph{Kazhdan-{L}usztig conjecture for affine {L}ie algebras with
  negative level. {II}. {N}onintegral case}, Duke Math. J. \textbf{84} (1996),
  no.~3, 771--813. \MR{1408544}

\bibitem[KT98]{kashiwaratanisakipositiverational}
\bysame, \emph{Kazhdan-{L}usztig conjecture for symmetrizable {K}ac-{M}oody
  {L}ie algebras. {III}. {P}ositive rational case}, Asian J. Math. \textbf{2}
  (1998), no.~4, 779--832, Mikio Sato: a great Japanese mathematician of the
  twentieth century. \MR{1734129}

\bibitem[KT00]{kashiwaratanisakialllevels}
\bysame, \emph{Characters of irreducible modules with non-critical highest
  weights over affine {L}ie algebras}, Representations and quantizations
  ({S}hanghai, 1998), China High. Educ. Press, Beijing, 2000, pp.~275--296.
  \MR{1802178}

\bibitem[KW07]{kapustinwitten}
Anton Kapustin and Edward Witten, \emph{Electric-magnetic duality and the
  geometric {L}anglands program}, Commun. Number Theory Phys. \textbf{1}
  (2007), no.~1, 1--236. \MR{2306566}

\bibitem[Lan93]{langlands1985harish}
Robert Langlands, \emph{Harish-{C}handra (11 {O}ctober 1923--16 {O}ctober
  1983)}, Current Sci. \textbf{65} (1993), no.~12, 922--936. \MR{1263347}

\bibitem[Lus83]{lusztigqchar}
George Lusztig, \emph{Singularities, character formulas, and a {$q$}-analog of
  weight multiplicities}, Analysis and topology on singular spaces, {II}, {III}
  ({L}uminy, 1981), Ast\'{e}risque, vol. 101, Soc. Math. France, Paris, 1983,
  pp.~208--229. \MR{737932}

\bibitem[Lus85a]{luschar1}
\bysame, \emph{Character sheaves. {I}}, Adv. in Math. \textbf{56} (1985),
  no.~3, 193--237. \MR{792706}

\bibitem[Lus85b]{luschar23}
\bysame, \emph{Character sheaves. {II}, {III}}, Adv. in Math. \textbf{57}
  (1985), no.~3, 226--265, 266--315. \MR{806210}

\bibitem[Lus86a]{luschar4}
\bysame, \emph{Character sheaves. {IV}}, Adv. in Math. \textbf{59} (1986),
  no.~1, 1--63. \MR{825086}

\bibitem[Lus86b]{luschar5}
\bysame, \emph{Character sheaves. {V}}, Adv. in Math. \textbf{61} (1986),
  no.~2, 103--155. \MR{849848}

\bibitem[Mil80]{milneetalebook}
James Milne, \emph{\'{E}tale cohomology}, Princeton Mathematical Series, No.
  33, Princeton University Press, Princeton, N.J., 1980. \MR{559531}

\bibitem[Mil89]{milnenotes}
\bysame, \emph{Lectures on \'{e}tale cohomology}, Notes, available at
  \url{https://www.jmilne.org/math/CourseNotes/LEC.pdf} (1989).

\bibitem[MP94]{moyprasad}
Allen Moy and Gopal Prasad, \emph{Unrefined minimal {$K$}-types for {$p$}-adic
  groups}, Invent. Math. \textbf{116} (1994), no.~1-3, 393--408. \MR{1253198}

\bibitem[MV07]{mirkovic-vilonen}
Ivan Mirkovi{\'c} and Kari Vilonen, \emph{Geometric {L}anglands duality and
  representations of algebraic groups over commutative rings}, Ann. of Math.
  (2) \textbf{166} (2007), no.~1, 95--143.

\bibitem[Orl04]{orlovsingularity}
Dmitri Orlov, \emph{Triangulated categories of singularities and {D}-branes in
  {L}andau-{G}inzburg models}, Tr. Mat. Inst. Steklova \textbf{246} (2004),
  no.~Algebr. Geom. Metody, Svyazi i Prilozh., 240--262. \MR{2101296}

\bibitem[Ras15a]{rdm}
Sam Raskin, \emph{{$D$}-modules on infinite dimensional varieties}, Preprint
  \url{https://web.ma.utexas.edu/users/sraskin/dmod.pdf} (2015).

\bibitem[Ras15b]{raskinspectraldecomp}
\bysame, \emph{On the notion of spectral decomposition in local geometric
  {L}anglands}, Preprint, arXiv:1511.01378, current version available at
  \url{https://web.ma.utexas.edu/users/sraskin/locsys.pdf} (2015).

\bibitem[Ras18]{raskin2ndadjoint}
\bysame, \emph{Second adjointness for loop groups},
  \url{https://lysenko.perso.math.cnrs.fr/Notes_talks_winter2018/Ja-1(Raskin).pdf}
  (2018).

\bibitem[Ras20a]{RaskinLocCrit}
\bysame, \emph{Affine {B}eilinson--{B}ernstein localization at the critical
  level for {GL}(2)}, Annals of Math, to appear (2020).

\bibitem[Ras20b]{mys}
\bysame, \emph{{H}omological methods in semi-infinite contexts}, Preprint,
  arXiv:2002.01395 (2020).

\bibitem[Ras21]{whit}
\bysame, \emph{{$\mathscr{W}$}-algebras and {W}hittaker categories}, Selecta
  Math. (N.S.) \textbf{27} (2021), no.~3, Paper No. 46, 114. \MR{4273643}

\bibitem[Roc98]{roche}
Alan Roche, \emph{Types and {H}ecke algebras for principal series
  representations of split reductive {$p$}-adic groups}, Ann. Sci. \'{E}cole
  Norm. Sup. (4) \textbf{31} (1998), no.~3, 361--413. \MR{1621409}

\bibitem[Rod75]{rodier}
Fran\c{c}ois Rodier, \emph{Mod\`ele de {W}hittaker et caract\`eres de
  repr\'{e}sentations}, Non-commutative harmonic analysis ({A}ctes {C}olloq.,
  {M}arseille-{L}uminy, 1974), 1975, pp.~151--171. Lecture Notes in Math., Vol.
  466. \MR{0393355}

\bibitem[RW18]{richewilliamsontilting}
Simon Riche and Geordie Williamson, \emph{Tilting modules and the
  {$p$}-canonical basis}, Ast\'{e}risque (2018), no.~397, ix+184. \MR{3805034}

\bibitem[RY]{RaskinYang}
Sam Raskin and David Yang, \emph{Affine {B}eilinson--{B}ernstein localization
  at the critical level}, In preparation.

\bibitem[Sag17]{sage2}
Daniel Sage, \emph{Regular strata and moduli spaces of irregular singular
  connections}, New trends in analysis and interdisciplinary applications,
  Trends Math. Res. Perspect., Birkh\"{a}user/Springer, Cham, 2017, pp.~69--75.
  \MR{3695627}

\bibitem[Sat63]{satake63}
Ichir\^{o} Satake, \emph{Theory of spherical functions on reductive algebraic
  groups over {${\mathfrak p}$}-adic fields}, Inst. Hautes \'{E}tudes Sci.
  Publ. Math. (1963), no.~18, 5--69. \MR{195863}

\bibitem[Sha90]{shahidi90}
Freydoon Shahidi, \emph{A proof of {L}anglands' conjecture on {P}lancherel
  measures; complementary series for {$p$}-adic groups}, Ann. of Math. (2)
  \textbf{132} (1990), no.~2, 273--330. \MR{1070599}

\bibitem[Spa88]{spaltenstein}
Nicolas Spaltenstein, \emph{Resolutions of unbounded complexes}, Compositio
  Math. \textbf{65} (1988), no.~2, 121--154. \MR{932640}

\bibitem[Sto06]{stoyanovsky}
Alexander Stoyanovsky, \emph{Quantum {L}anglands duality and conformal field
  theory}, Preprint, arXiv:0610974 (2006).

\bibitem[SV17]{sakellaridisvenkatesh}
Yiannis Sakellaridis and Akshay Venkatesh, \emph{Periods and harmonic analysis
  on spherical varieties}, Ast\'{e}risque (2017), no.~396, viii+360.
  \MR{3764130}

\bibitem[To{\"{e}}07]{toenpaper}
Bertrand To{\"{e}}n, \emph{The homotopy theory of {$dg$}-categories and derived
  {M}orita theory}, Invent. Math. \textbf{167} (2007), no.~3, 615--667.
  \MR{2276263}

\bibitem[To{\"{e}}11]{toenlectures}
\bysame, \emph{Lectures on dg-categories}, Topics in algebraic and topological
  {$K$}-theory, Lecture Notes in Math., vol. 2008, Springer, Berlin, 2011,
  pp.~243--302. \MR{2762557}

\bibitem[TT90]{thomasontrobaugh}
Robert Thomason and Thomas Trobaugh, \emph{Higher algebraic {$K$}-theory of
  schemes and of derived categories}, The {G}rothendieck {F}estschrift, {V}ol.
  {III}, Progr. Math., vol.~88, Birkh\"{a}user Boston, Boston, MA, 1990,
  pp.~247--435. \MR{1106918}

\bibitem[Ver66]{verma66}
Daya-Nand Verma, \emph{Structure of certain induced representations of complex
  semisimple {L}ie algebras}, ProQuest LLC, Ann Arbor, MI, 1966, Thesis
  (Ph.D.)--Yale University. \MR{2615829}

\bibitem[Vog93]{voganloclang}
David Vogan, \emph{The local {L}anglands conjecture}, Representation theory of
  groups and algebras, Contemp. Math., vol. 145, Amer. Math. Soc., Providence,
  RI, 1993, pp.~305--379. \MR{1216197}

\bibitem[Wil17]{willimasonsurvey1}
Geordie Williamson, \emph{Algebraic representations and constructible sheaves},
  Jpn. J. Math. \textbf{12} (2017), no.~2, 211--259. \MR{3694932}

\bibitem[Wit08]{wittenwild}
Edward Witten, \emph{Gauge theory and wild ramification}, Anal. Appl. (Singap.)
  \textbf{6} (2008), no.~4, 429--501. \MR{2459120}

\bibitem[Yan21]{yangmoyprasad}
David Yang, \emph{Categorical {M}oy--{P}rasad theory}, Preprint,
  arXiv:2104.12917 (2021).

\bibitem[Zhu21]{zhufamilies}
Xinwen Zhu, \emph{Coherent sheaves on the stack of {L}anglands parameters},
  Preprint, arXiv:2008.02998 (2021).

\end{thebibliography}

\end{document}